\documentclass[a4paper,10pt,reqno]{amsart}

\usepackage[utf8]{inputenc}
\usepackage[T1]{fontenc}
\usepackage{lmodern}
\usepackage[english]{babel}
\usepackage{microtype}

\usepackage{amsmath,amssymb,amsfonts,amsthm,comment}
\usepackage{mathtools,accents}
\usepackage{mathrsfs}
\usepackage{aliascnt}
\usepackage{braket}
\usepackage{bm}
\usepackage{esint}
\usepackage{amsfonts}
\usepackage{csquotes}

\usepackage[a4paper,twoside,top=0.85in, bottom=0.85in, left=0.7in, right=0.7in]{geometry} 
\usepackage{hyperref}
\usepackage[english,capitalize]{cleveref}

\usepackage{enumerate}
\usepackage{xcolor}
\usepackage{comment}

\usepackage{imakeidx} %Analytical index
\usepackage{xparse} %In order to use ExplSyntax (in older latex versions)
\usepackage{mathtools}
\usepackage{float}
\usepackage{comment}
\usepackage{nicefrac}
\usepackage{graphicx}
\usepackage{enumerate}
\usepackage{enumitem}

\definecolor{newblue}{rgb}{0.2, 0.3, 0.85}
\hypersetup{colorlinks=true, linkcolor=newblue, citecolor=newblue, urlcolor  = newblue}

\usepackage[maxbibnames=99,backend=biber, sorting=nyt, doi=false, url=false, isbn=false]{biblatex}
\makeatletter
\g@addto@macro\@floatboxreset\centering
\makeatother

\makeatletter
\def\newaliasedtheorem#1[#2]#3{
  \newaliascnt{#1@alt}{#2}
  \newtheorem{#1}[#1@alt]{#3}
  \expandafter\newcommand\csname #1@altname\endcsname{#3}
}
\makeatother

\numberwithin{equation}{section}

\newtheoremstyle{slanted}{\topsep}{\topsep}{\slshape}{}{\bfseries}{.}{.5em}{}

\theoremstyle{plain}
\newtheorem{theorem}{Theorem}[section]
\newaliasedtheorem{proposition}[theorem]{Proposition}
\newaliasedtheorem{question}[theorem]{Question}
\newaliasedtheorem{problem}[theorem]{Problem}
\newaliasedtheorem{lemma}[theorem]{Lemma}
\newaliasedtheorem{corollary}[theorem]{Corollary}
\newaliasedtheorem{counterexample}[theorem]{Counterexample}

\theoremstyle{definition}
\newaliasedtheorem{definition}[theorem]{Definition}
\newaliasedtheorem{example}[theorem]{Example}
\newaliasedtheorem{conjecture}[theorem]{Conjecture}

\theoremstyle{remark}
\newaliasedtheorem{remark}[theorem]{Remark}
\newcommand{\setN}{\mathbb{N}}
\newcommand{\N}{\mathbb{N}}

\newcommand{\setR}{\mathbb{R}}
\newcommand{\R}{\mathbb{R}}
\renewcommand{\S}{\mathbb{S}}

\newcommand{\X}{{\rm X}}
\newcommand{\lip}{{\rm lip \,}}
\newcommand{\nchi}{{\raise.3ex\hbox{\(\chi\)}}}
\newcommand{\eps}{\varepsilon}

\let\phi\varphi

\newcommand{\di}{\mathop{}\!\mathrm{d}}

\newcommand{\res}{\mathop{\hbox{\vrule height 7pt width .5pt depth 0pt
\vrule height .5pt width 6pt depth 0pt}}\nolimits}

\newcommand{\Ch}{{\sf Ch}}
\newcommand{\st}{\ensuremath{\ :\ }} % Such that in formulas
\newcommand{\eqdef}{\ensuremath{\vcentcolon=}}

\newcommand{\haus}{\mathcal{H}}

\newcommand{\dist}{\mathsf{d}}

\newcommand{\meas}{\mathfrak{m}}

\newcommand{\diam}{\mathrm{diam}}
\DeclareMathOperator{\CD}{CD}
\DeclareMathOperator{\RCD}{RCD}
\DeclareMathOperator{\CBB}{CBB}

\DeclareMathOperator{\Ric}{Ric}
\DeclareMathOperator{\Per}{Per}

\newfont{\tmpf}{cmsy10 scaled 2500}

\newcommand{\de}{\ensuremath{\,\mathrm d}} % Il de degli integrali (contiene lo spazio da mettere tra integranda e misura)

\author{Gioacchino Antonelli}
\address{Courant Institute Of Mathematical Sciences (NYU), 251 Mercer Street, 10012, New York, USA}
\email{ga2434@nyu.edu}

\author{Marco Pozzetta}
\address{Dipartimento di Matematica e Applicazioni, Universit\`a di Napoli Federico II, Via Cintia, Monte S. Angelo, 80126 Napoli, Italy.}
\email{marco.pozzetta@unina.it}


\begin{document}
\title[Isoperimetric problem and structure at infinity on $\CBB(0)$ spaces]{Isoperimetric problem and structure at infinity on Alexandrov spaces with nonnegative curvature}

\maketitle

\begin{abstract}
In this paper we consider nonnegatively curved finite dimensional Alexandrov spaces with a non-collapsing condition, i.e.,
such that unit balls have volumes uniformly bounded from below away from zero.
We study the relation between the isoperimetric profile, the existence of isoperimetric sets, and the asymptotic structure at infinity of such spaces.

In this setting, we prove that the following conditions are equivalent: the space has linear volume growth; it is Gromov--Hausdorff asymptotic to one cylinder at infinity; it has uniformly bounded isoperimetric profile; the entire space is a tubular neighborhood of either a line or a ray.

Moreover, on a space satisfying any of the previous conditions, we prove existence of isoperimetric sets for sufficiently large volumes, and we characterize the geometric rigidity at the level of the isoperimetric profile. 

Specializing our study to the $2$-dimensional case, we prove that unit balls have always volumes uniformly bounded from below away from zero, and we prove existence of isoperimetric sets for every volume, characterizing also their topology when the space has no boundary. 

The proofs exploit a variational approach, and in particular apply to Riemannian manifolds with nonnegative sectional curvature and to Euclidean convex bodies. Up to the authors' knowledge, most of the results are new even in these smooth cases.
\end{abstract}

\tableofcontents

\section{Introduction}

We consider Alexandrov spaces $(X,\dist)$ with curvature bounded from below by $0$, shortly $\CBB(0)$, of (Hausdorff) dimension $N\in\N$ with $N\ge 2$, see \cref{def:CBB}. For an introduction to the theory of such spaces we refer to \cite{BuragoBuragoIvanovBook, AKP, PetruninSurvey}, and to the foundational works \cite{Perelman91, BuragoGromovPerelman}. The class of $\CBB(0)$ spaces includes complete Riemannian manifolds with nonnegative sectional curvature (possibly with convex boundary), convex bodies, i.e., closure of open convex subsets of $\R^N$, topological boundaries of open convex sets of subsets of $\R^N$ (endowed with corresponding intrinsic distance), and cones over closed manifolds of sectional curvature greater or equal to $1$ (see \cref{sec:Cones}). Moreover the $\CBB$ condition is stable with respect to quotients by isometries, gluing over boundaries, {and pointed Gromov--Hausdorff convergence}. 

In the following, an $N$-dimensional $\CBB(0)$ space $(X,\dist)$ is endowed with the $N$-dimensional Hausdorff measure $\haus^N$, which makes the triple $(X,\dist,\haus^N)$ a metric measure space. The basic theory of $BV$-functions on metric measure spaces \cite{Miranda03, Ambrosio02} allows to define the perimeter $\Per (E)$ of a measurable subset $E\subset X$, see \cref{def:BVperimetro}. Thus we consider the \emph{isoperimetric problem}: given an $N$-dimensional $\CBB(0)$ space $(X,\dist)$, for $V \in (0,\haus^N(X))$ one aims at studying the minimization
\[
\inf \{ \Per(E) \st E\subset X, \,\,\haus^N(E)=V \}.
\]
The value of the infimum as a function $I_X:(0,\haus^N(X))\to [0,+\infty]$ of the volume is called \emph{isoperimetric profile} of $X$. A set $E\subset X$ such that $\Per (E)=I_X(\haus^N(E))$, is called \emph{isoperimetric set}, or \emph{region}. We shall drop the subscript $X$ when there is no risk of confusion.

\medskip

In the last fifty years, a great interest has been devoted to the study of several aspects of the isoperimetric problem on spaces with (classical or synthetic) notions of lower bounds on the curvature. We mention \cite{BerardMeyer82, BerardBessonGallotiso, BavardPansu86, Gromovmetric, MorganJohnson00, Milmanconvexity, Milmanmodels, CavallettiMaggiMondino} concerning sharp, rigid and quantitative isoperimetric inequalities on compact spaces with Ricci lower bounds, \cite{AgostinianiFogagnoloMazzieri, BrendleFigo, BaloghKristaly, Johne, APPSb, CavallettiManini, CavallettiManiniRigidity} about sharp and rigid isoperimetric inequalities on noncompact spaces with nonnegative Ricci curvature and Euclidean volume growth, see \eqref{eq:DefAVR} below, \cite{BavardPansu86, SternbergZumbrun, MorganJohnson00, KuwertIsop, Bayle03, Bayle04, BayleRosales, MondinoNardulli16, APPSa} for what concerns differential properties of the isoperimetric profile of spaces with Ricci lower bounds, and \cite{RitoreExistenceSurfaces01, RitRosales04, MorganRitore02, Nar14, MondinoNardulli16, ChodoshEichmairVolkmann17, AFP21, AntBruFogPoz, AntonelliNardulliPozzetta} studying existence of isoperimetric sets on noncompact spaces with lower curvature bounds.\\
Specializing to the case of convex bodies in the Euclidean space, a number of results about existence of minimizers, isoperimetric inequalities, stability and topology of isoperimetric sets have been obtained in \cite{LionsPacella, SternbergZumbrun, RitoreVernadakis15, LeonardiRitore}.\\
We refer the reader to \cite{RitoreBook} for a comprehensive account on the isoperimetric problem.

\medskip

On a compact space, precompactness {of sets with uniformly bounded perimeter} and lower semicontinuity of the perimeter functional imply existence of isoperimetric sets of any volume $V \in (0,\haus^N(X))$. On the contrary, the problem of existence of minimizers for some assigned volume is highly nontrivial on noncompact spaces. In fact, it is proved in \cite{AntonelliGlaudo} that for every dimension $N\geq 3$ there exists a convex body $X\subset \R^N$ with infinite volume such that no isoperimetric sets exist for volumes $<1$\footnote{In this case the perimeter of a set $E\subset X$ is the relative perimeter in the convex body. Indeed, the perimeter measure $\mathrm{Per}(E,\cdot)$ of a set $E$ never charges the boundary $\partial X$, see \cref{sec:RCD}.}, thus answering in the negative a question in \cite{LeonardiRitore}. A similar example can be constructed where $(X,\dist)$ is a smooth complete noncompact Riemannian manifold with dimension $N\geq 3$ and strictly positive sectional curvature. A reason why isoperimetric sets of some volumes do not exist on some noncompact space is that minimizing sequences may escape from any compact region, see examples in \cite{Rit01NonExistence, AFP21}. It follows that there is a strong relation between the existence of isoperimetric sets on noncompact spaces, and the asymptotic structure at infinity.

{
The link between the problem and the structure at infinity,} together with recent developments on the Geometric Analysis on metric measure spaces with curvature bounded below (\cref{sec:RCD}), has led to several contributions which exploit the direct method of the Calculus of Variations in order to find an isoperimetric region. 
Building on the fundamental concentration-compactness principle \cite{Lions84I} and the results in \cite{RitRosales04, Nar14, MondinoNardulli16}, such direct method has been established in \cite{AFP21, AntonelliNardulliPozzetta}, see also \cite{NobiliVioloStability, Reinaldo2020, NovagaClusters} for related applications, leading to new existence results of minimizers under assumptions on the asymptotic structure of the ambient {space}.

In \cite{AntBruFogPoz} it is considered the case of Riemannian manifolds with nonnegative Ricci curvature with \emph{Euclidean volume growth}, i.e., such that
\begin{equation}\label{eq:DefAVR}
{\rm AVR}(X) \eqdef \lim_{r\to+\infty} \frac{\haus^N(B_r(o))}{\omega_N r^N}  \in (0,1].
\end{equation}
By Bishop--Gromov monotonicity (\cref{sec:RCD}), the previous assumption implies that the volume of balls grows with maximal rate. Under a suitable assumption on the geometry at infinity of the space, in \cite[Theorem 1.1]{AntBruFogPoz} it is established existence of isoperimetric sets for any sufficiently large volume. Furthermore, if $(X,\dist,\haus^N)$ is $\CBB(0)$ and ${\rm AVR}(X)>0$, then isoperimetric sets exist for any sufficiently large volume without further assumptions \cite[Theorem 1.3]{AntBruFogPoz}, and \cite[Theorem 1.2]{APPSb}. The results in \cite{AntBruFogPoz} are generalized in the $\RCD$ setting in \cite{APPSb}.

\medskip

In this work, we focus on $\CBB(0)$ spaces whose balls have \emph{minimal volume growth} instead. By a classical result attributed to Calabi--Yau, see \cite[Theorem 4.1]{SchoenYauLectures}, \cite[Corollary 5.15]{Gigli12}, and also \cite{CheegerGromovTaylor} where explicit lower bounds on the Ricci curvature are related to the sharp rate of volume growth of balls, the volume growth of balls in a $\CBB(0)$ space $(X,\dist)$ is at least linear, more precisely $\haus^N(B_r(o)) \ge C(N, \haus^N(B_1(o))) r$ for every $o \in X$, $r>1$.

In the class of $\CBB(0)$ spaces, our first main result yields a full equivalence between linear volume growth and: the cylindrical structure at infinity of the space, an upper bound on the isoperimetric profile, and the fact that the space is a tubular neighborhood of either a line or a ray. Moreover, a rigidity statement holds in the sense that the space is a cylinder or it has a cylindrical end - necessarily, over isoperimetric boundaries - if and only if the isoperimetric profile is eventually constant.

\begin{theorem}[{Isoperimetry and structural rigidity, cf. \cref{thm:raggiosselimiteainfinito}, \cref{cor:Equivalenze}, \cref{prop:ProfiloCostiffSPlitCilindro}}]\label{thm:IntroEquivalenzeStruttura}
Let $N\geq 2$, and let $(X,\dist)$ be a noncompact $N$-dimensional $\CBB(0)$ metric space such that $\mathcal{H}^N(B_1(x))\geq v_0$ for some $v_0>0$, and for every $x\in X$. Let $o \in X$. Then the following are equivalent. 
\begin{enumerate}[label=(\roman*)]
    \item $\limsup_{r\to+\infty} \haus^N(B_r(o)) /r <+\infty$;
    
    \item For any diverging sequence of points $p_i\in X$ such that $(X,p_i)$ pGH-converges to a limit $(Y,y)$, it holds that $Y$ is isometric to a product $\R\times K$, where $K$ is compact;
    
    \item There exists a diverging sequence of points $p_i\in X$ such that $(X,p_i)$ pGH-converges to a limit $\R\times K$, where $K$ is compact;

    \item There exists a constant $C>0$ such that $I(V)\leq C$ for every $V>0$.
    
    \item Either $X$ is isometric to $\R\times K$, or $X$ has one end and there exists a ray $\gamma:[0,+\infty)\to X$ and a constant $D>0$ such that $\dist(p,\gamma)\leq D$ for any $p\in X$.
\end{enumerate}
{If any of the previous items holds, then there exists a unique compact $\CBB(0)$ space $K$ such that any pGH limit of $X$ along a sequence of diverging points is isometric to $\R\times K$.}

Moreover, the following are equivalent.
    \begin{enumerate}
        \item There exists $V_0>0$ such that the isoperimetric profile $I$ is constant on $[V_0,+\infty)$;
        \item Either $X$ is isometric to $\mathbb R\times K$ for some compact $K$, or there exists an isoperimetric set $A$ such that $(\partial A, \dist|_{\partial A})$ is $\CBB(0)$ and $(X \setminus A, \dist|_{X\setminus A})$ is isometric to $(\partial A\times [0,+\infty), \dist|_{\partial A}\otimes\dist_{\rm eu})$;
        \item Either $X$ is isometric to $\mathbb R\times K$ for some compact $K$, or there exists an open bounded set $A$ such that $(X \setminus A, \dist|_{X\setminus A})$ is $\CBB(0)$ and it is isometric to $(\partial A\times [0,+\infty), \dist|_{\partial A}\otimes\dist_{\rm eu})$.
    \end{enumerate}
Moreover, if (1), (2), or (3) holds, either $X$ is isometric to $\mathbb R\times K$ for some compact $K$ and there exists $V_0>0$ such that for every $V\geq V_0$ isoperimetric regions of volume $V$ exist and coincide with $(0,a)\times K$ for some $a>0$, up to translations; or there exists $V_0>0$ such that for every $V\geq V_0$ isoperimetric regions of volume $V$ exist, are unique, and coincide with $\{x \st \dist(x,A) < b_V\}$ for some $b_V\geq 0$.
\end{theorem}

The most challenging part of the previous \cref{thm:IntroEquivalenzeStruttura} is the equivalence of the conditions in (i)-(ii)-(iii)-(v) with (iv), i.e., the connection between the structure of the space and the boundedness of the isoperimetric profile. {
The most difficult implications, in which the $\CBB(0)$ condition comes into play, are (i)$\Rightarrow$(ii), and (iii)$\Rightarrow$(v), from which we can infer the global structure of the space just by having one limit at infinity that is of the form $\mathbb R\times K$}. The $\CBB(0)$ condition comes into play also in the proof of (i)$\Leftrightarrow$(iv), see \cref{problemequivalencelineargrowth}.

Up to the authors' knowledge the equivalence of items (i)-(v) and the rigidity in items (1)-(3) in \cref{thm:IntroEquivalenzeStruttura} are new both in the case of Euclidean convex bodies, and in the class of Riemannian manifolds with nonnegative sectional curvature. Moreover, since the proof exploits the analysis of pointed Gromov--Hausdorff limits at infinity, which are a priori non-smooth, the use of non-smooth geometry appears to be necessary even in the smooth realm.

In the case of manifolds with nonnegative Ricci curvature, or even $\RCD(0,N)$ spaces $(X,\dist,\haus^N)$ with linear volume growth, the picture depicted in \cref{thm:IntroEquivalenzeStruttura} might be more involved, see \cref{rem:SormaniResults} and \cref{rem:CompletenessBusemann} below. Moreover, there exist examples of manifolds with nonnegative Ricci and linear volume growth having limits at infinity of the form $\R\times K$ and $\R\times K'$ where $K$ and $K'$ are non-isometric compact manifolds \cite[Example 27]{SormaniJDG98}. The investigation of the relations between conditions (i)-(v) in \cref{thm:IntroEquivalenzeStruttura} at such level of generality goes beyond the scope of the present note. 

\begin{problem}\label{problemequivalencelineargrowth}
Investigate the relation between conditions (i)--(v) in \cref{thm:IntroEquivalenzeStruttura} in the class of non-collapsed noncompact $\RCD(0,N)$ metric measure spaces $(X,\dist,\mathcal{H}^N)$. In particular, investigate whether the implications {(i)$\Rightarrow$(ii)$\Rightarrow$(v)} and (i)$\Leftrightarrow$(iv) hold in the previous setting.
\end{problem}

{
\textbf{Addendum}. 
After this work was posted on arXiv, X. Zhu \cite{XingyuZhu} addressed \cref{problemequivalencelineargrowth}, proving that on a noncompact manifold $(M^N,g)$ with $\Ric\ge0$ and $\inf_{x\in M} \haus^N(B_1(x))>0$, conditions (i), (ii), (v) in \cref{thm:IntroEquivalenzeStruttura} are equivalent \cite[Theorem 1.6]{XingyuZhu} It remains open to understand the full equivalence with the uniform boundedness of the isoperimetric profile, that is, condition (iv) in \cref{thm:IntroEquivalenzeStruttura}.
}
\smallskip

Exploiting the structural information coming from \cref{thm:IntroEquivalenzeStruttura}, we prove that on a $\CBB(0)$ space with minimal volume growth there exist isoperimetric sets for any sufficiently large volume.

\begin{theorem}[{Existence and structural rigidity, cf. \cref{thm:AsintoticaProfilo}}]\label{thm:IntroExistence}
Let $N\geq 2$, and let $(X,\dist)$ be a noncompact $N$-dimensional $\CBB(0)$ metric space such that $\mathcal{H}^N(B_1(x))\geq v_0>0$ for some $v_0>0$ for every $x\in X$. Assume any of the equivalent conditions (i)-(v) in \cref{thm:IntroEquivalenzeStruttura} and assume that $X$ has one end\footnote{If it has at least two ends, then it is isometric to $\mathbb R\times K$ by the splitting theorem, see \cref{thm:SplittingTheorem}.}. Then the following hold.
\begin{enumerate}
\item All the pGH limits at infinity of $X$ {are isometric to $\mathbb R\times K$, for a fixed compact $\CBB(0)$ space $K$. Denote $\mathfrak{V}\eqdef \mathcal{H}^{N-1}(K)$;}
\item There exists a constant {$V_0:=V_0(v_0,N,\mathfrak{V})$} such that isoperimetric regions exist for every volume $V\geq V_0$;
\item For every $V>0$ we have that 
    \[
    I(V)\leq \mathfrak{V}.
    \]
    If for some $V'$ the equality is attained, then it is attained for all $V\geq V'$ and there exists an isoperimetric set $\Omega$ of volume $\haus^N(\Omega)\ge V'$ such that $(\partial \Omega, \dist|_{\partial \Omega})$ is $\CBB(0)$ and $(X \setminus\Omega, \dist|_{X\setminus\Omega})$ is isometric to $(\partial\Omega\times [0,+\infty), \dist|_{\partial \Omega}\otimes\dist_{\rm eu})$.
\item {If $E_i$ is a sequence of isoperimetric sets with $\haus^N(E_i)\to+\infty$, $x_i \in \partial E_i$, and $(X,\dist,x_i)$ pGH-converges to $(\R\times K,\dist_{\rm eu} \otimes \dist_{K}, (0,k_0))$, then $E_i$ converges to $(-\infty,0)\times K$ in $L^1_{\rm loc}$ and $\partial E_i$ converges to $(-\infty,0)\times K$ in Hausdorff distance, in a realization of the pGH-convergence.} Moreover there exists a sequence of isoperimetric sets $\Omega_i\Subset \Omega_{i+1}$ such that $\cup_i \Omega_i = X$, and there holds
    \[
    \lim_{V\to+\infty} I(V) =  \mathfrak{V}.
    \]
\end{enumerate} 
\end{theorem}

Let us mention that \cref{thm:IntroEquivalenzeStruttura} and \cref{thm:IntroExistence} recover several results in the context of Euclidean convex bodies proved in \cite{RVcylindricallyBounded, LeonardiRitore}, and extend the theory to the larger class of $\CBB(0)$ spaces. Indeed, every convex body in $\mathbb R^n$ that satisfies item (v) in \cref{thm:IntroEquivalenzeStruttura} is {\em cylindrically bounded} in the terminology of \cite{RVcylindricallyBounded}, see also \cite[Theorem 6.20]{LeonardiRitore}. We mention that in order to prove \cref{thm:IntroExistence} we need some auxiliary results about the characterization of isoperimetric sets with large volumes in $\RCD$ cylinders over compact spaces which might be of independent interest, see \cref{AppendixA}.

\medskip

We can further specialize our approach to the case of $2$-dimensional $\CBB(0)$ spaces $(X,\dist)$. In this case we prove both that the volume non-collapsing assumption $\inf_{x \in X}\haus^2(B_1(x))>0$ holds true, and that isoperimetric sets exist for every volume. We further investigate the topology of isoperimetric sets in case the space has no boundary, see \cref{sec:Cones}. We denote by $\mathbb S^1(\rho)$ the circle of length $\rho$.

\begin{theorem}[{$2$-dimensional $\CBB(0)$ spaces, cf. \cref{thm:Noncollapsed}, \cref{thm:Existence2D}, \cref{thm:IsopBoundaryConnected}}]\label{thm:IntroDimension2}
Let $(X,\dist)$ be a noncompact $2$-dimensional $\CBB(0)$ metric space. Then the following hold.
\begin{enumerate}
    \item There exists $v_0>0$ such that $\mathcal{H}^2(B_1(x))\geq v_0$ for any $x\in X$.
    
    \item For any $V>0$ there exists an isoperimetric region $\Omega\subset X$ of volume $\haus^2(\Omega)=V$.
\end{enumerate}
Assume further that $\partial X=\emptyset$, and let $E\subset X$ be an isoperimetric set. Then one of the following alternatives holds.
\begin{enumerate}[label=(\alph*)]
\item $\partial E$ is connected.
\item There exist $b,\rho>0$ such that $X$ is isometric to a right cylinder $\mathbb R\times \mathbb S^1(\rho)$ and $E=(0,b)\times\mathbb S^1(\rho)$, up to translation along the factor $\mathbb R$.
\end{enumerate}
In particular, in case $\partial X=\emptyset$ and $X$ is not isometric to a Euclidean $2$-dimensional right cylinder, the boundary of every isoperimetric set is connected and it is parametrized by a Lipschitz curve homeomorphic to $\mathbb S^1$.
\end{theorem}

By a classification argument, \cref{thm:IntroDimension2} implies the existence of isoperimetric sets on any $\RCD(0,2)$ space $(X,\dist,\meas)$, see \cref{cor:RCDK2}. We mention that an alternative argument for the proof of (1) in \cref{thm:IntroDimension2} has been recently proposed in \cite[Lemma 2.12, Remark 2.13]{PositiveScalarRicci}.

We observe that \cref{thm:IntroDimension2} refutes the existence of foliations by isoperimetric boundaries even at the level of $2$-dimensional $\CBB(0)$ spaces, and despite existence of isoperimetric sets for any volume, see \cref{rem:Foliations}.

We remark that (1) in \cref{thm:IntroDimension2} recovers the result known on $2$-dimensional Riemannian manifolds from \cite[Theorem A]{CrokeKarcher}. We stress that property (1) is false on higher dimensional $\CBB(0)$ spaces, see \cite[Example 1]{CrokeKarcher}, or consider the convex hull of the union of the unit ball in $\R^3$ with a parabola $\{(x,x^2,0)\st x\in \R\}$. Moreover (2) in \cref{thm:IntroDimension2} recovers the result known on $2$-dimensional Riemannian manifolds from \cite{RitoreExistenceSurfaces01}, and (a) in \cref{thm:IntroDimension2} recovers the topological conclusion in \cite[Theorem 1.1]{RitoreExistenceSurfaces01}.\\
We remark that (1), (2) in \cref{thm:IntroDimension2} here are proved by applying the above mentioned direct method from \cite{AFP21, AntonelliNardulliPozzetta}, thus unifying these results as corollaries of the method.

\medskip

It remains open the following fundamental question related to the existence part in \cref{thm:IntroExistence}. We remark that, up to the authors' knowledge, the forthcoming question has not a definitive answer both in the class of smooth complete Riemannian manifolds with $\mathrm{Sect}\geq 0$ and in the class of Euclidean convex bodies.

\begin{question}
Let $N\ge 2$. Let $(X,\dist)$ be an $N$-dimensional $\CBB(0)$ noncompact metric space such that $\inf_{x\in X}\mathcal{H}^N(B_1(x))\geq v_0$ for some $v_0>0$. Do there exist isoperimetric regions for sufficiently large volumes?
\end{question}

Concerning the topological result in \cref{thm:IntroDimension2}, we mention the following question. 

\begin{question}
Let $(X,\dist)$ be a noncompact $2$-dimensional $\CBB(0)$ metric space. If the space is not a $2$-dimensional cylinder, is it true that isoperimetric sets are convex and homeomorphic to a disk?
\end{question}

The question whether isoperimetric sets are homeomorphic to disks has a positive answer on smooth complete Riemannian surfaces with $\mathrm{Sect}\geq 0$, see \cite{RitoreExistenceSurfaces01}.

\medskip

Concerning the strategy to carry out the analysis of the structure at infinity of a $\CBB(0)$ space, a first key observation is the fact that any pGH limit at infinity of a $\CBB(0)$ space splits a line, i.e., it is isometric to a product $\R\times Y$, see \cref{lem:EveryLimitAtInfinitySplits}. Later, we found that this property had already observed in a work of Kleiner--Lott \cite[Page 2852, Appendix G]{KleinerLott}. Then for a $\CBB(0)$ space having at least a limit at infinity of the form $\R\times K$ with $K$ compact, it is possible to derive the asymptotic behavior of the Busemann function (\cref{def:Busemann}) of a given ray, see \cref {thm:raggiosselimiteainfinito}, and \cref{prop:StimaFq}. Hence, exploiting also recent rigidity results from \cite{Ketterer21, KettererKitabeppuLakzian}, we derive the following theorem on the asymptotic behavior of the sublevel sets of the Busemann function.

\begin{theorem}[cf. \cref{thm:ConvergenzaBusemann}]\label{thm:IntroConvergenzaBusemann}
Let $N\geq 2$, and let $(X,\dist)$ be a noncompact $N$-dimensional $\CBB(0)$ metric space such that $\mathcal{H}^N(B_1(x))\geq v_0$ for some $v_0>0$, and for every $x\in X$.
Assume that there exists a ray $\gamma:[0,+\infty)\to X$ such that
\[
\dist(p,\gamma)\leq C, \qquad \text{for all $p\in X$},
\]
for some constant $C>0$. Let $F$ be the Busemann function of $\gamma$, see \cref{def:Busemann}.
Let $t_j\to +\infty$ and assume that $(X,\gamma(t_j))$ pGH-converges to $(\mathbb R\times K,(0,k_0))$ for some $(N-1)$-dimensional compact $\CBB(0)$ space $K$.

Let us denote by $f:\mathbb R\times K\to \R$ the coordinate function $f(t,k)=t$. Let $F_j:=F-t_j$.
Then the following hold.
\begin{itemize}
    \item For every $j\in\mathbb N$ and every $T\geq 0$, the set $\{F_j\leq T\}$ is compact.
    
    \item $F_j$ converges locally uniformly to $f$ along the pointed sequence of converging metric spaces.
    
    \item For any $s \in \R$, the set $\{F<  s\}$ has finite perimeter and
    \begin{equation*}
    \mathrm{Per}(\{F< s_1\})\leq \mathrm{Per}(\{F< s_2\}),\qquad \text{$\forall\, s_1, s_2 \in \R$ with $s_1<s_2$,}
    \end{equation*}
    \begin{equation*}
        \lim_j \Per(\{F_j<t\}) =\lim_{s\to +\infty}\mathrm{Per}(\{F< s\})=\mathcal{H}^{N-1}(K) 
        \qquad \forall\, t \in \R.
    \end{equation*}
    Moreover, if $\Per (\{F< \bar s\})=\mathcal{H}^{N-1}(K)$ for some $\bar s \in \R$, then there exists $s_0\ge \bar s$ such that $(\{F=s_0\}, \dist|_{\{F=s_0\}})$ is isometric to $K$, and $(X\setminus \{F<s_0\}, \dist|_{X\setminus \{F<s_0\}})$ is isometric to $(\{F=s_0\}\times [0,+\infty),\dist|_{\{F=s_0\}} \otimes \dist_{\rm eu})$, which is isometric to $K\times[0,+\infty)$.
    
    \item {If $q_j$ is an arbitrary sequence of diverging points on $X$, then up to subsequences $(X,q_j)$ pGH-converges to $\mathbb R\times K$.}
\end{itemize}
\end{theorem}

Eventually, \cref{thm:IntroEquivalenzeStruttura}, \cref{thm:IntroExistence} and \cref{thm:IntroDimension2} follow by variational arguments exploiting \cref{thm:IntroConvergenzaBusemann} and the approach by direct method developed in \cite{AFP21, AntonelliNardulliPozzetta}.

\medskip
\noindent\textbf{Organization.} In \cref{sec:Preliminari} we recall definitions and basic results on $\CBB(0)$ spaces, $\RCD$ spaces, and sets of finite perimeter. We also derive a series of technical lemmas on the asymptotic geometry of $\CBB(0)$ spaces, see \cref{sec:AsympGeometry}. In \cref{sec:Structure} we introduce Busemann functions, we prove \cref{thm:IntroConvergenzaBusemann} and we derive the equivalence of items (i)-(v) in \cref{thm:IntroEquivalenzeStruttura}. In \cref{sec:Dimension2} we prove the volume non-collapsedness of balls in $2$-dimensional $\CBB(0)$ spaces arguing by direct method. In \cref{sec:IsoperimetricProblem} we consider the isoperimetric problem on $\CBB(0)$ metric spaces, completing the proof of the results stated in the Introduction.\\
For the convenience of the reader, appendices collect statements and proofs of some results that may be well-known to experts in the field. \cref{AppendixA} characterizes the isoperimetric sets for large volumes in $\RCD$ cylinders over compact spaces, \cref{AppendixB} contains the characterization of rays in product spaces and cylinders, while in \cref{AppendixC} we provide a proof of the fact that the distance from a line is a concave function on $2$-dimensional $\CBB(0)$ spaces.

\medskip
\noindent\textbf{Acknowledgements.} The first author acknowledges the financial support of the Courant Institute, and the AMS-Simons Travel grant. The second author is partially supported by the INdAM - GNAMPA Project 2022 CUP\_E55F22000270001 ``Isoperimetric problems: variational and geometric aspects''. The authors thank Christian Ketterer for useful discussions on the results in \cite{KettererKitabeppuLakzian, Ketterer21}. They also thank Elia Bruè, Mattia Fogagnolo, Federico Glaudo, Andrea Mondino, Daniele Semola, Christina Sormani, and Guofang Wei for many fruitful discussions on the topic of this work during various stages of the writing of this paper. The authors would also like to thank Alexander Lytchak for several comments that improved the exposition and the results of this note. The authors wish to thank the anonymous referee for several comments that improved the readability of the paper.

\section{Preliminaries}\label{sec:Preliminari}

\subsection{Notation and general facts about $\mathrm{CBB}(0)$ metric spaces}

In this section we discuss general facts about $\CBB(0)$ metric spaces.
Let $(X,\dist)$ be a geodesic metric space. A geodesic $\gamma:[a,b]\to X$ is a curve with finite length $L(\gamma)$ such that $L(\gamma)=\dist(\gamma(a),\gamma(b))$. {Except otherwise stated, we shall assume that geodesics are parametrized with metric speed equal to $1$ almost everywhere.} 
For any \(k\in[0,\infty)\), we denote by \(\haus^k\) the \emph{\(k\)-dimensional Hausdorff measure} on \((X,\dist)\).\\
Let $p,x,y \in X$. A given geodesic from $x$ to $y$ parametrized with metric speed equal to $1$ at a.e. time is denoted by $[xy]$. With a little abuse of notation, we shall denote by $\gamma$ or $[xy]$ both the parametrization and the image of the curve $\gamma$ or $[xy]$. A given hinge determined by two geodesics $[px]$ and $[py]$ is denoted by $p^x_y$. A triangle with vertices $p,x,y$ determined by three given geodesics is denoted by $\triangle pxy$.

For $\bar p, \bar x, \bar y \in \R^2$, the angle at $\bar p$ of the triangle with vertices $\bar p, \bar x, \bar y$ is denoted by $\angle^0 \bar x \bar p \bar y$, that is, $\angle^0 \bar x \bar p\bar y = \arccos \frac{|\bar p \bar x|^2 + |\bar p \bar y|^2 - |\bar x \bar y|^2}{2|\bar p \bar x||\bar p \bar y|}$.

A triangle $\triangle \bar p \bar x \bar y \subset \R^2$ is a comparison triangle for $\triangle pxy$ if $|\bar p \bar x|=\dist(p,x)$, $|\bar p \bar y|=\dist(p,y)$, $|\bar x \bar y|=\dist(x,y)$.

Denoting $\alpha=[px]$ and $\beta=[py]$, the hinge $p^x_y$ in $X$ has angle at $p$ if there exists the limit 
\[
\angle p^x_y  \eqdef \lim_{s,t\to0^+} \angle^0 \bar \alpha(s) \bar p \bar \beta(t),
\]
where $\triangle \bar \alpha(s) \bar p \bar \beta(t) \in \R^2$ is a comparison triangle for $\triangle \alpha(s) p  \beta(t)$. In our notation it is therefore consistent to write $\angle^0 \bar x \bar p \bar y = \angle \bar p^{\bar x}_{\bar y}$ for any $\bar p, \bar x, \bar y \in \R^2$. In order to avoid confusion, we will denote $\angle^0 \bar p^{\bar x}_{\bar y}\eqdef \angle \bar p^{\bar x}_{\bar y}$ the angle of a hinge in $\R^2$.

{ Given $\alpha=[px],\beta=[py]$ two geodesics emanating from the same point $p$ in $X$, we define
\begin{equation}\label{eqn:DefAngolo}
\angle \alpha'(0)\beta'(0):=\angle p^{\alpha(t)}_{\beta(s)},
\end{equation}
for any $t \in (0,|px|)$, $s \in (0,|py|)$. Notice that the definition is well-posed, i.e., it does not depend on the choice of $t,s$, by definition of angle.
}
\begin{definition}\label{def:CBB}
	Let $(X,\dist)$ be a geodesic proper, and thus complete, metric space. We say that $(X,\dist)$ is $\CBB(0)$, i.e., it has curvature bounded below by $0$, if any of the following equivalent (see \cite[Theorem 4.3.5, Proposition 10.1.1, Theorem 10.3.1]{BuragoBuragoIvanovBook}, \cite[Theorem 8.14]{AKP}) properties hold.
	\begin{enumerate}
		\item \emph{Four point comparison.} For any $p, x_1, x_2, x_3 \in X$ there holds
		\[
		\angle^0 \bar p_{ \bar x_1}^{\bar x_2} + \angle^0 \bar p_{\bar x_2 }^{\bar x_3} + \angle^0  \bar p_{\bar x_1}^{\bar x_3} \le 2\pi,
		\]
		where $ \bar p_{\bar x_i}^{\bar x_j}$ is the hinge in a comparison triangle $\triangle \bar x_i \bar p \bar x_j$ for $\triangle x_i p x_j$, for any $i,j=1,2,3$.
		
		\item \emph{Monotonicity condition.} For any $p, x, y \in X$, denoting $\alpha=[px]$ and $\beta=[py]$, the function
		\[
		(0,\dist(p,x)]\times (0,\dist(p,y)] \ni (s,t) \mapsto \frac{s^2+t^2-\dist(\alpha(s), \beta(t))^2}{2ts},
		\]
		is increasing with respect to each variable. Equivalently,
		\[
		(0,\dist(p,x)]\times (0,\dist(p,y)] \ni (s,t) \mapsto \arccos \frac{s^2+t^2-\dist(\alpha(s), \beta(t))^2}{2ts},
		\]
		is decreasing with respect to each variable.
		
		\item \emph{Triangle (or distance) condition.} For any triangle $\triangle p x y \subset  X$, if $ \triangle \bar p \bar x \bar y \in \R^2$ is a comparison triangle, then
		\[
		\dist(p,z)\ge |\bar p \bar z| \qquad \forall\, z \in [xy], \, \bar z \in [\bar x \bar y] \st \dist(x,z)=|\bar x \bar z|.
		\]
		
		\item\label{it:AngleCondition} \emph{Angle condition.} For any triangle $\triangle p x y \subset  X$ the angles $\angle p^x_y$, $\angle x^p_y$, $\angle y^p_x$ are well defined and, if $ \triangle \bar p \bar x \bar y \in \R^2$ is a comparison triangle, then
		\[
		\angle p^x_y \ge \angle^0 \bar x \bar p \bar y,
		\qquad
		\angle x^p_y \ge \angle^0\bar p \bar x \bar y,
		\qquad  
		\angle y^p_x \ge \angle^0\bar p \bar y \bar x .
		\]
		Moreover for any geodesics $[xy]$, $[pz]$, where $z \in [xy]\setminus\{x,y\}$, there holds $\angle z^p_x + \angle z^p_y = \pi$.
	\end{enumerate}
\end{definition}

If $(X,\dist)$ is $\CBB(0)$, the following further comparison condition holds.\footnote{To our knowledge, it is not known whether the hinge condition implies $\CBB(0)$, see also \cite[\text{8.51}]{AKP} and \cite[p. 108]{BuragoBuragoIvanovBook}.}
\begin{enumerate}
    \item[(5)] \emph{Hinge condition.}  For any hinge $p^x_y \subset X$, the angle $\angle p^x_y$ is well defined and, if $\bar p^{\bar x}_{\bar y} \subset \R^2$ is a hinge in $\R^2$ such that $\dist(p,x)=|\bar p \bar x|$, $\dist(p,y)=|\bar p \bar y|$, $\angle p^x_y \le \angle^0 \bar p^{\bar x}_{\bar y}$, then
	\[
		\dist(x,y)\le |\bar x \bar y|.
	\]
\end{enumerate}
The previous condition readily follows from the monotonicity condition. Indeed by monotonicity we have
\[
\begin{split}
    |\bar p \bar x|^2+ |\bar p \bar y|^2 -\dist(x,y)^2 & = 
\dist(p,x)^2+\dist(p,y)^2-\dist(x,y)^2 \ge 2\dist(p,x)\dist(p,y) \cos\angle p^x_y 
\\& \ge 2 |\bar p \bar x| \, |\bar p \bar y| \cos \angle^0 \bar p^{\bar x}_{\bar y},
\end{split}
\]
which implies $\dist(x,y)^2 \le  |\bar p \bar x|^2+ |\bar p \bar y|^2 -  2|\bar p \bar x| \, |\bar p \bar y| \cos \angle^0 \bar p^{\bar x}_{\bar y} = |\bar x \bar y|^2$.

\bigskip

When there is no risk of confusion, we shall denote by $|\cdot|$ the distance on a given metric space $(X,\dist)$, i.e., $|xy|\eqdef \dist(x,y)$ for $x, y \in X$.

\medskip

{We shall need the following elementary observation, whose proof readily follows by angle and monotonicity conditions.}

\begin{lemma}\label{lem:GeodCadePerpendicolare}
Let $(X,\dist)$ be $\CBB(0)$. Let $\gamma:[0,T]\to X$ be a geodesic. Let $q \in X$ and $t_q \in (0,T)$ be such that
\[
\dist(q,\gamma)=|q\gamma(t_q)|>0,
\]
and fix a geodesic $[\gamma(t_q)q]$.
Then $\angle \gamma(t_q)^q_{\gamma(0)} = \angle \gamma(t_q)^q_{\gamma(T)}=\tfrac\pi2$.
\end{lemma}

\subsubsection{Pointed measured Gromov--Hausdorff convergence}
{In the following, we shall always assume that a metric measure space (m.m.s. for short) $(X,\dist,\meas)$ is endowed with a nonnegative Borel measure that is finite on bounded subsets.}
We introduce the pointed measured Gromov--Hausdorff convergence already in a proper realization, which, in our setting, is equivalent to the standard definition of pmGH convergence, see \cite[Theorem 3.15 and Section 3.5]{GigliMondinoSavare15}. The following definition is taken from the introductory exposition of \cite{AmbrosioBrueSemola19}.

\begin{definition}[pGH and pmGH convergence]\label{def:GHconvergence}
A sequence $\{ (X_i, \dist_i, x_i) \}_{i\in \N}$ of pointed metric spaces is said to converge in the \emph{pointed Gromov--Hausdorff topology, in the $\mathrm{pGH}$ sense for short,} to a pointed metric space $ (Y, \dist_Y, y)$ if there exist a complete separable metric space $(Z, \dist_Z)$ and isometric embeddings
\[
\begin{split}
&\Psi_i:(X_i, \dist_i) \to (Z,\dist_Z), \qquad \forall\, i\in \N\, ,\\
&\Psi:(Y, \dist_Y) \to (Z,\dist_Z)\, ,
\end{split}
\]
such that for any $\eps,R>0$ there is $i_0(\varepsilon,R)\in\mathbb N$ such that
\[
\Psi_i(B_R^{X_i}(x_i)) \subset \left[ \Psi(B_R^Y(y))\right]_\eps,
\qquad
\Psi(B_R^{Y}(y)) \subset \left[ \Psi_i(B_R^{X_i}(x_i))\right]_\eps\, ,
\]
for any $i\ge i_0$, where $[A]_\eps\coloneqq \{ z\in Z \st \dist_Z(z,A)\leq \eps\}$ for any $A \subset Z$.

Let $\meas_i$ and $\mu$ be given in such a way $(X_i,\dist_i,\meas_i,x_i)$ and $(Y,\dist_Y,\mu,y)$ are m.m.s.\! If in addition to the previous requirements we also have $(\Psi_i)_\sharp\mathfrak{m}_i \rightharpoonup \Psi_\sharp \mu$ with respect to duality with continuous bounded functions on $Z$ with bounded support, then the convergence is said to hold in the \emph{pointed measured Gromov--Hausdorff topology, or in the $\mathrm{pmGH}$ sense for short}.
\end{definition}

The following definitions are given in \cite[Definition 3.1]{AmbrosioBrueSemola19}, and it is investigated in \cite{AmbrosioBrueSemola19} capitalizing on the results in \cite{AmbrosioHonda17}.

\begin{definition}[$L^1$-strong and $L^1_{\mathrm{loc}}$ convergence]\label{def:L1strong}
Let $\{ (X_i, \dist_i, \mathfrak{m}_i, x_i) \}_{i\in \N}$  be a sequence of pointed metric measure spaces converging in the pmGH sense to a pointed metric measure space $ (Y, \dist_Y, \mu, y)$ and let $(Z,\dist_Z)$ be a realization as in \cref{def:GHconvergence}.

We say that a sequence of Borel sets $E_i\subset X_i$ such that $\mathfrak{m}_i(E_i) < +\infty$ for any $i \in \N$ converges \emph{in the $L^1$-strong sense} to a Borel set $F\subset Y$ with $\mu(F) < +\infty$ if $\mathfrak{m}_i(E_i) \to \mu(F)$ and $\chi_{E_i}\mathfrak{m}_i \rightharpoonup \chi_F\mu$ with respect to the duality with continuous bounded functions with bounded support on $Z$.

We say that a sequence of Borel sets $E_i\subset X_i$ converges \emph{in the $L^1_{\mathrm{loc}}$-sense} to a Borel set $F\subset Y$ if $E_i\cap B_R(x_i)$ converges to $F\cap B_R(y)$ in $L^1$-strong for every $R>0$.
\end{definition}

\begin{definition}[Hausdorff convergence]
Let $\{ (X_i, \dist_i, \mathfrak{m}_i, x_i) \}_{i\in \N}$  be a sequence of pointed metric measure spaces converging in the pmGH sense to a pointed metric measure space $ (Y, \dist_Y, \mu, y)$.

We say that a sequence of closed sets $E_i\subset X_i$ converges \emph{in Hausdorff distance} (or \emph{in Hausdorff sense}) to a closed set $F\subset Y$ if there holds convergence in Hausdorff distance in a realization $(Z,\dist_Z)$ of the pmGH convergence as in \cref{def:GHconvergence}.
\end{definition}

\begin{remark}[Uniform convergence of curves and functions along pGH limits]
Let $\{ (X_i, \dist_i, x_i) \}_{i\in \N}$  be a sequence of pointed metric spaces converging in the pGH sense to a pointed metric space $ (Y, \dist_Y, y)$. Let $(Z,\dist_Z)$ be a realization of the pGH-convergence as in \cref{def:GHconvergence}.
Let $I$ be a compact interval, and let $\gamma_i:I\to X_i$, $\gamma:I\to X$ be continuous curves such that $\gamma_i(0)=x_i$, $\gamma(0)=x$.\\
We say that {\em $\gamma_i$ converges uniformly to $\gamma$} if the following holds. For every $\varepsilon$ there exists $i_0\in\mathbb N$ such that 
\[
\dist_Z(\Psi_i(\gamma_i(t)),\Psi(\gamma(t)))\leq \varepsilon,\qquad \forall t\in I,\quad \forall i\geq i_0.
\]

A slight modification of the Arzelà--Ascoli theorem gives the following. If the curves $\gamma_i$ as above are equicontinuous, then, up to subsequences, they converge uniformly to a continuous curve $\gamma$ as above. Notice that the equicontinuity is readily verified when the curves are equi-$L$-Lipschitz, and in this case the limit curve is $L$-Lipschtiz as well. In particular, if $\gamma_i$ are unit speed geodesics, the uniform limit curve $\gamma$ exists up to subsequences and it is a unit speed geodesic thanks to the continuity of the distance.

Analogous definition and result hold for a sequence of continuous functions. More precisely, we say that a sequence of continuous functions $f_i:X_i\to \R$ \emph{converges locally uniformly} to a continuous function $f:Y\to \R$ if for any $R,\eps>0$ there exist $i_0 \in \N$ and $\delta>0$ such that $|f(z)-f_i(x)|<\eps$ whenever $z \in B_R(y)$, $i\ge i_0$ and $\dist_Z(\Psi_i(x), \Psi(z)) < \delta$. Moreover, an equibounded equicontinuous sequence of functions has a locally uniformly converging subsequence by an Arzelà--Ascoli-type argument.
\end{remark}

\begin{remark}[Stability and projection]
From \cite[Proposition 10.7.1]{BuragoBuragoIvanovBook} we have that the any pGH limit of $\CBB(0)$ metric spaces is still a $\CBB(0)$ metric space.\\
Moreover, we have that if $\mathbb R\times X$ is a $\CBB(0)$ metric space when endowed with the product distance $\dist_{\rm eu} \otimes \dist_X\eqdef \sqrt{\dist_{\mathrm{eu}}^2+\dist_X^2}$, then $X$ is a $\CBB(0)$ metric space when endowed with the distance $\dist_X$. This is a direct outcome of the definitions.
\end{remark}

\begin{remark}[{Angle continuity \cite[Section 8.40, 8.41]{AKP}}]\label{rem:AngleContinuity}
Let $(X_n,\dist_n,o_n)$ be a sequence of $\CBB(0)$ spaces converging in GH-sense to $(X,\dist,o)$. Let $p_n, x_n, y_n \in X_n$ and hinges ${p_n}^{x_n}_{y_n}$ be converging to $p,x,y \in X$ and $p^x_y$. Then the following holds.
\begin{itemize}
	\item $\liminf_n \angle {p_n}^{x_n}_{y_n} \ge \angle p^x_y$.
	
	\item If $p, x, y$ are distinct and $[xp]\subset [xp']$ for some $p'\neq p$, then $\lim_n \angle {p_n}^{x_n}_{y_n} = \angle p^x_y$, whenever the limit exists.	
\end{itemize}
\end{remark}

\subsubsection{Tangent and asymptotic cones of $\CBB(0)$ metric spaces}\label{sec:Cones}

Let us recall here the classical Splitting Theorem. In the setting of manifolds with nonnegative Ricci curvature it has been proved by Cheeger and Gromoll in \cite{CheegrGromollSplitting}. The version for Alexandrov spaces with nonnegative sectional curvature below is due to \cite{MilkaSplitting}. Finally, we recall that a version for $\RCD$ spaces has been proved by Gigli in \cite{Gigli13}. When we say that $(X,\dist)$ is an {\em $N$-dimensional $\CBB(0)$ metric space} we mean that its Hausdorff dimension is $N$. The relations between different notions of dimensions is presented, e.g., in \cite[Section 15]{AKP}.

\begin{theorem}[Splitting Theorem, \cite{MilkaSplitting, CheegrGromollSplitting, Gigli13}]\label{thm:SplittingTheorem}
Let $(X,\dist)$ be an $N$-dimensional $\CBB(0)$ space. Let $\gamma:\R\to X$ be a line. Then $X$ is isometric to a product space $\R\times Y$ endowed with the product distance, where $(Y,\dist_Y)$ is an $(N-1)$-dimensional $\CBB(0)$ space.
\end{theorem}

{ Let $(X,\dist)$ be a metric space with $\mathrm{diam} (X)\leq \pi$. Let us define $C(X):=\left([0,+\infty)\times X\right)/\{0\}\times X$, and let us consider the \textit{cone} distance
\[
\dist_c((t,p),(s,q)):=\sqrt{t^2+s^2-2ts\cos(\dist(p,q))}
\]
on it. The space $(C(X),\dist_c)$ is said {\em metric cone} over $(X,d)$, see \cite[Section 3.6]{BuragoBuragoIvanovBook}. Sometimes we will refer to $(C(X),\dist_c)$ as the {\em Euclidean cone over $X$}.

Notice that, if $\mathrm{diam} (X)>\pi$, one can define $(C(X),\dist_c)$ simply as $(C(\overline X),\overline \dist_c)$, where $(\overline X,\overline d)$ is the metric space that is $X$ set-wise and, for every $p,q\in X$,
\[
\overline \dist(p,q):=\min\{\dist(p,q),\pi\},
\]
see \cite[Definition 3.6.16]{BuragoBuragoIvanovBook}. }

Let us quickly review the concepts of tangent cone, boundary, and asymptotic cone of $\CBB(0)$ metric spaces. For boundary and tangent cones, we refer the reader to classical references, e.g., \cite{Perelman91, BuragoGromovPerelman, BuragoBuragoIvanovBook}; the following facts hold also for arbitrary lower bounds on the curvature.\\
Let $X$ be an $N$-dimensional $\CBB(0)$ metric space, and let $x\in X$. The {\em space of directions} $\Sigma_x(X)$ is the metric completion of the space of geodesics emanating from $x$ endowed with the metric given by the angle between them (two geodesics are identified if their angle is zero). The metric space $\Sigma_x(X)$ is a compact $(N-1)$-dimensional metric space with curvature bounded below by $1$, see \cite[Sect. 4.6]{BuragoBuragoIvanovBook}. The tangent cone at $x$ is the metric cone over $\Sigma_x(X)$, which is itself an $N$-dimensional $\CBB(0)$ metric space \cite[Theorem 4.7.1]{BuragoBuragoIvanovBook}. The pGH limit of $(X,r^{-1}\dist,x)$ as $r\to 0$ is isometric to the tangent cone at $x$ \cite[Sect. 8.2.1]{BuragoBuragoIvanovBook}, and moreover balls centered at $x$ with sufficiently small radii are homeomorphic to the tangent cone at $x$, see \cite[Section 13.2, Theorem b)]{BuragoGromovPerelman} and \cite{Perelman91}.\\
The boundary $\partial X$ is inductively defined as follows. If $X$ has dimension $N\geq 2$, we define
\[
\partial X:=\{x\in X:\Sigma_x(X)\, \text{has boundary}\}.
\]
This is a well-posed definition after noticing that $1$-dimensional $\CBB$ spaces are $1$-dimensional manifolds.

\begin{definition}[Metric on the rays]\label{def:IdealBoundary}
    Let $(X,\dist)$ be a noncompact $N$-dimensional $\CBB(0)$ space. Fix $o \in X$ and let $\mathscr{R}_o$ be the set of rays emanating from $o$.
    For any $\gamma,\sigma \in \mathscr{R}_o$ , we define
    \[
    \angle_\infty (\gamma,\sigma) \eqdef \lim_{t\to+\infty} \arccos \Big( 1 - \dist^2(\gamma(t),\sigma(t))/(2t^2) \Big),
    \]
    {where the limit exists by \cref{def:CBB}(2).}
    Two rays $\gamma,\sigma \in \mathscr{R}_o$ are said to be equivalent, and we write $\gamma\sim\sigma$, if and only if $\angle_\infty(\gamma,\sigma)=0$, i.e., if and only if $\lim_{t\to+\infty} \dist(\gamma(t),\sigma(t))/t =0$. We denote by $[\gamma]$ the equivalence class of rays equivalent to $\gamma$.\\
    The \emph{ideal boundary} $X_\infty$ is the metric space $(\mathscr{R}_o/\sim, \angle_\infty)$. The distance $\angle_\infty$ on $X_\infty$ is called \emph{Tits metric}.
\end{definition}

\begin{theorem}[Asymptotic cones of $\CBB(0)$ spaces]\label{thm:AsymptoticCones}
Let $(X,\dist)$ be a noncompact $N$-dimensional $\CBB(0)$ space. Then the following hold.
\begin{enumerate}
    \item Every connected component of the ideal boundary $(X_\infty,\angle_\infty)$ is a compact metric space with curvature bounded below by $1$.
    \item $(X,r^{-1}\dist,o)$ converges to a unique metric cone $(C,\dist_C,o_C)$ as $r\to +\infty$ in the pointed Gromov--Hausdorff topology, and $C$ is isometric to the Euclidean cone over the ideal boundary $(X_\infty,\angle_\infty)$.
    \item $(X,\dist)$ splits a line if and only if $(C,\dist_C)$ splits a line.
\end{enumerate}
\end{theorem}
{
\begin{proof}
Proofs of the existence and uniqueness of the asymptotic cone $C$, together with its characterization as the Euclidean metric cone over $(X_\infty,\angle_\infty)$, can be found in \cite{BallmannGromovSchroeder85, Kasue88, GuijarroKapovitch95, ShiohamaBOOK, MashikoNaganoOtsuka}. A proof of claim (3) can be found in the second part of the proof of \cite[Theorem 4.6]{AntBruFogPoz}, which only exploits the $\CBB(0)$ condition.
\end{proof}
}

We conclude with a brief technical lemma.

\begin{lemma}\label{lem:AngoliRaggiCrossSection}
Let $(C,\dist)$ be an $N$-dimensional $\CBB(0)$ cone over a compact metric space $(X,\dist_X)$, and denote {$o=\{0\}\times X/\sim$}. Then for any ray $\gamma$ from $o$, any $y \in \partial B_1(o) \setminus\{\gamma(1)\}$ and any geodesic $[\gamma(1)y]$ there holds
\[
\angle \gamma(1)^o_y <\pi.
\]
\end{lemma}

\begin{proof}
If by contradiction $\angle \gamma(1)^o_y=\pi$, then $\angle \gamma(1)^{\gamma(s)}_y=0$ for any $s>1$. Hence the monotonicity condition implies that
\[
1 = \lim_{t\to0^+} 1 - \frac{\dist^2([\gamma(1)y](t), \gamma(1+t))}{2t^2} \le 1 - \frac{\dist^2([\gamma(1)y](t), \gamma(1+t))}{2t^2},
\]
for any $t \in (0,|\gamma(1)y|)$.
Hence $[\gamma(1)y](t)=\gamma(1+t)$ for $t \in [0, |\gamma(1)y|]$. But this implies that $1=|oy|=|o\gamma(1+|\gamma(1)y|)|=1+|\gamma(1)y|>1$.
\end{proof}

\subsection{$\mathrm{BV}$ and Sobolev functions on metric measure spaces}
In this paper, by a \emph{metric measure space} (briefly, m.m.s.) we mean a triple \((X,\dist,\meas)\), where \((X,\dist)\) is a complete and separable
metric space, while \(\meas\geq 0\) is a boundedly-finite Borel measure on \(X\).
\subsubsection{\texorpdfstring{$\rm BV$}{BV} functions and sets of finite perimeter in metric measure spaces}
Given a locally Lipschitz function $f:X\to\mathbb R$,
\[
\lip f (x) \eqdef \limsup_{y\to x} \frac{|f(y)-f(x)|}{\dist(x,y)},
\]
is the \emph{slope} of $f$ at $x$, for any accumulation point $x\in X$, and $\lip f(x)\coloneqq 0$ if $x\in X$ is isolated.

We give now the definitions of \emph{function of bounded variation} and \emph{set of finite perimeter} in the setting of m.m.s..
\begin{definition}[$\rm BV$ functions and perimeter on m.m.s.]\label{def:BVperimetro}
Let $(X,\dist,\meas)$ be a metric measure space.  Given $f\in L^1_{\mathrm{loc}}(X,\meas)$ we define
\[
|Df|(A) \eqdef \inf\left\{\liminf_i \int_A \lip f_i \de\meas \st \text{$f_i \in {\rm LIP}_{\mathrm{loc}}(A),\,f_i \to f $ in $L^1_{\mathrm{loc}}(A,\meas)$} \right\}\, ,
\]
for any open set $A\subset X$.
We declare that a function \(f\in L^1_{\mathrm{loc}}(X,\meas)\) is of \emph{local bounded variation}, briefly \(f\in{\rm BV}_{\mathrm{loc}}(X)\),
if \(|Df|(A)<+\infty\) for every \(A\subset X\) open bounded.
A function $f \in L^1(X,\meas)$ is said to belong to the space of \emph{bounded variation functions} ${\rm BV}(X)={\rm BV}(X,\dist,\meas)$ if $|Df|(X)<+\infty$. 

If $E\subset\X$ is a Borel set and $A\subset X$ is open, we  define the \emph{perimeter $\Per(E,A)$  of $E$ in $A$} by
\[
\Per(E,A) \eqdef \inf\left\{\liminf_i \int_A \lip u_i \de\meas \st \text{$u_i \in {\rm LIP}_{\mathrm{loc}}(A),\,u_i \to \nchi_E $ in $L^1_{\mathrm{loc}}(A,\meas)$} \right\}\, ,
\]
in other words \(\Per(E,A)\coloneqq|D\nchi_E|(A)\).
We say that $E$ has \emph{locally finite perimeter} if $\Per(E,A)<+\infty$ for every open bounded set $A$. We say that $E$ has \emph{finite perimeter} if $\Per(E,X)<+\infty$, and we denote $\Per(E)\eqdef \Per(E,X)$.
\end{definition}

Let us remark that when $f\in{\rm BV}_{\mathrm{loc}}(X,\dist,\meas)$ or $E$ is a set with locally finite perimeter, the set functions $|Df|, \Per(E,\cdot)$ above are restrictions to open sets of Borel measures that we still denote by $|Df|, \Per(E,\cdot)$, see \cite{AmbrosioDiMarino14}, and \cite{Miranda03}.
\medskip

In the sequel, we shall frequently make use of the following \emph{coarea formula}, proved in \cite{Miranda03}.
\begin{theorem}[Coarea formula]\label{thm:coarea}
Let \((X,\dist,\meas)\) be a locally compact metric measure space. Let \(f\in L^1_{\mathrm{loc}}(X)\) be given. Then for any open set
\(\Omega\subset X\) it holds that \(\R\ni t\mapsto \Per(\{f>t\},\Omega)\in[0,+\infty]\) is Borel measurable and satisfies
\[
|Df|(\Omega)=\int_\R \Per(\{f>t\},\Omega)\,\de t\, .
\]
In particular, if \(f\in{\rm BV}(X)\), then \(\{f>t\}\) has finite perimeter for a.e.\ \(t\in\setR\).
\end{theorem}
\begin{remark}[Semicontinuity of the total variation under $L^1_{\mathrm{loc}}$-convergence]\label{rem:SemicontPerimeter}
Let $(X,\dist,\meas)$ be a metric measure space. We recall (cf.\ \cite[Proposition 3.6]{Miranda03}) that whenever $g_i,g\in L^1_{\mathrm{loc}}(X,\meas)$
are such that $g_i\to g$ in $L^1_{\mathrm{loc}}(X,\meas)$, for every open set $\Omega$ we have 
$$
|Dg|(\Omega)\leq \liminf_{i\to +\infty}|Dg_i|(\Omega)\, .
$$
\end{remark}

\begin{remark}[Precompactness and lower semicontinuity of finite perimeter sets along pmGH converging sequences]\label{rem:SemicontPerimeterConverging}
Let $K\in\mathbb R$, $N\geq 1$, and $\{(X_i,\dist_i,\meas_i,x_i)\}_{i\in\mathbb N}$ be a sequence of $\RCD(K,N)$ metric measure spaces converging in the pmGH sense to $(Y,\dist_Y,\mu,y)$. Let $(Z,\dist_Z)$ be a realization of the convergence. Then, the following hold, compare with \cite[Proposition 3.3, Corollary 3.4, Proposition 3.6, Proposition 3.8]{AmbrosioBrueSemola19}, and \cite{AmbrosioHonda17}.
\begin{itemize}
    \item For any sequence of Borel sets $E_i\subset X_i$ with 
    $$
    \sup_{i\in\mathbb N}|D\chi_{E_i}|(B_R(x_i))<+\infty, \qquad \forall\,R>0,
    $$
    there exists a subsequence $i_k$ and a Borel set $F\subset Y$ such that $E_{i_k}\to F$ in $L^1_{\mathrm{loc}}$.
    \item Let $F\subset Y$ be a bounded set of finite perimeter. Then there exist a subsequence $i_k$, and uniformly bounded finite perimeter sets $E_{i_k}\subset X_{i_k}$ such that $E_{i_k}\to F$ in $L^1$-strong and $|D\chi_{E_{i_k}}|(X_{i_k})\to |D\chi_F|(Y)$ as $k\to+\infty$.
    \item Let $f_i\in \mathrm{BV}(X_i,\dist_i,\mathfrak{m}_i)$ converge in $L^1$-strong to $f\in L^1(Y,\mu)$. If $\sup_i|Df_i|(X_i)<+\infty$, then $f\in \mathrm{BV}(Y,\dist_Y,\mu)$, and we have
    \[
    |Df|(Y)\leq \liminf_{i\to+\infty}|Df_i|(X_i).
    \]
\end{itemize}
\end{remark}

\begin{definition}[Isoperimetric profile]
    Let $(X,\dist,\mathfrak{m})$ be a metric measure space. The \textit{isoperimetric profile} $I$ of $X$ is the function $I:(0,\mathfrak{m}(X))\to[0,+\infty]$
    \[
    I(V):=\inf\{\mathrm{Per}(E):E\subset X,\,E\,\text{Borel},\mathfrak{m}(E)=V\}.
    \]
\end{definition}

\subsubsection{Sobolev functions in metric measure spaces}
The \emph{Cheeger energy} on a metric measure space \((X,\dist,\meas)\) is defined as the \(L^2\)-relaxation of the functional
\(f\mapsto\frac{1}{2}\int\lip^2 f_n\di\meas\) (see \cite{AmbrosioGigliSavare11} after \cite{Cheeger99}). Namely, for any function \(f\in L^2(X)\) we define
\[
\Ch(f)\coloneqq\inf\bigg\{\liminf_{n\to\infty}\frac{1}{2}\int\lip^2 f_n\di\meas\;\bigg|\;(f_n)_{n\in\setN}\subset{\rm LIP}_{bs}(X),\,f_n\to f\text{ in }L^2(X)\bigg\}\, .
\]
The \emph{Sobolev space} \(H^{1,2}(X)\) is defined as the finiteness domain \(\{f\in L^2(X)\,:\,\Ch(f)<+\infty\}\) of the Cheeger energy.
The restriction of the Cheeger energy to the Sobolev space admits the integral representation \(\Ch(f)=\frac{1}{2}\int|\nabla f|^2\di\meas\),
for a uniquely determined function \(|\nabla f|\in L^2(X)\) that is called the \emph{minimal weak upper gradient} of \(f\in H^{1,2}(X)\).
The linear space \(H^{1,2}(X)\) is a Banach space if endowed with the Sobolev norm
\[
\|f\|_{H^{1,2}(X)}\coloneqq\sqrt{\|f\|_{L^2(X)}^2+2\Ch(f)}=\sqrt{\|f\|_{L^2(X)}^2+\||\nabla f|\|_{L^2(X)}^2},\quad\text{ for every }f\in H^{1,2}(X)\, .
\]
Following \cite{Gigli12}, when \(H^{1,2}(X)\) is a Hilbert space (or equivalently \(\Ch\) is a quadratic form) we say that the metric measure
space \((X,\dist,\meas)\) is \emph{infinitesimally Hilbertian}.
\medskip

We define the bilinear mapping \(H^{1,2}(X)\times H^{1,2}(X)\ni(f,g)\mapsto\nabla f\cdot\nabla g\in L^1(X)\) as
\[
\nabla f\cdot\nabla g\coloneqq\frac{|\nabla(f+g)|^2-|\nabla f|^2-|\nabla g|^2}{2},\quad\text{ for every }f,g\in H^{1,2}(X)\, .
\]

\subsection{Geometric Analysis and isoperimetric sets on \texorpdfstring{$\RCD$}{RCD} spaces}\label{sec:RCD}
We avoid giving a detailed introduction to the theory, addressing the reader to the survey \cite{AmbrosioSurvey} and references therein for the relevant background. Below we recall some of the main properties that will be relevant for our purposes.\\
The Riemannian Curvature Dimension condition $\RCD(K,\infty)$ was introduced in \cite{AmbrosioGigliSavare14} coupling the Curvature Dimension condition $\CD(K,\infty)$, previously proposed in \cite{Sturm1,Sturm2} and independently in \cite{LottVillani}, with the {linearity of the heat flow}.

The finite dimensional refinements subsequently led to the notions of $\RCD(K,N)$ and $\RCD^*(K, N)$ spaces, corresponding to $\CD(K, N)$ (resp. $\CD^*(K, N)$, see \cite{BacherSturm}) coupled with {the infinitesimally Hilbertian condition}. The class $\RCD(K,N)$ was proposed in \cite{Gigli12}. The (a priori more general) $\RCD^*(K,N)$ condition was thoroughly analysed in \cite{ErbarKuwadaSturm15} and (subsequently and independently) in \cite{AmbrosioMondinoSavare15} (see also \cite{CavallettiMilmanCD} and the recent \cite{LiGlobalization} for the equivalence between $\RCD^*$ and $\RCD$ conditions).
\medskip

The $\RCD(K,N)$ condition is compatible with the smooth notion. In particular, smooth $N$--dimensional Riemannian manifolds with Ricci curvature bounded from below by $K$ endowed with the canonical volume measure are $\RCD(K,N)$ spaces. Smooth Riemannian manifolds with smooth and convex boundary (i.e. non negative second fundamental form with respect to the interior unit normal) are also included in the theory { (we refer to \cite{Han20} for detailed results, also covering the case of weighted manifolds).}\\ 
We recall that if $(X,\dist)$ is $\CBB(0)$ as in \cref{def:CBB}, then it has a well defined notion of dimension $N \in \N$, that coincides with the Hausdorff dimension of $X$ (see \cite[Chapter 15]{AKP}). In the following, we will always assume that $N \ge 2$ and that a $\CBB(0)$ space is endowed with the $N$-dimensional Hausdorff measure induced by its distance.

We recall that if $(X,\dist)$ is an $N$-dimensional $\CBB(0)$ metric space, then the triple $(X,\dist,\haus^N)$ is an $\RCD(0,N)$ metric measure space, see \cite{PetruninAlexandrovCD, GigliKuwadaOhta}. Moreover, $(X,\dist)$ is a $2$-dimensional $\CBB(0)$ metric space if and only if $(X,\dist,\haus^2)$ is $\RCD(0,2)$, see \cite{LytchakStadler}. Analogous results hold for any lower bound on the curvature.

\medskip

A fundamental property of $\RCD(K,N)$ spaces is the stability with respect to pmGH-convergence, meaning that a pmGH-limit of a sequence of (pointed) $\RCD(K_n,N_n)$ spaces for some $K_n\to K$ and $N_n\to N$ is an $\RCD(K,N)$ metric measure space

The classical Bishop--Gromov comparison theorem holds \cite{VillaniBook}, that is the functions
\[
r\mapsto \frac{\mathcal{H}^N(B_r(x))}{r^N},
\qquad
r\mapsto \frac{\Per(B_r(x))}{r^{N-1}},
\]
are nonincreasing on $(0,+\infty)$ for any $x \in X$. We shall denoted by $\omega_N$ for the volume of the unit ball in $\mathbb R^N$.
\smallskip

Let us now recall some results about sets of locally finite perimeter in $\RCD$ spaces. Given a Borel set \(E\subset X\) in an \(\RCD(K,N)\) space \((X,\dist,\haus^N)\) and any \(t\in[0,1]\), we denote by \(E^{(t)}\) the set
of \emph{points of density \(t\)} of \(E\), namely
\[
E^{(t)}\coloneqq\bigg\{x\in X\;\bigg|\;\lim_{r\to 0}\frac{\haus^N(E\cap B_r(x))}{\haus^N(B_r(x))}=t\bigg\}\, .
\]
The \emph{essential boundary} of \(E\) is defined as \(\partial^e E\coloneqq X\setminus(E^{(0)}\cup E^{(1)})\). 
The \emph{reduced boundary} \(\mathcal F E\subset\partial^e E\) of a set of locally finite perimeter \(E\subset X\) is defined as the set of the points of \(X\) where the unique tangent to \(E\), up to isomorphism, is the
half-space \(\{x=(x_1,\ldots,x_N)\in\setR^N\,:\,x_N>0\}\) in \(\setR^N\); see \cite[Definition 4.1]{AmbrosioBrueSemola19} for the notion of convergence.

{
It was proved in \cite{BruePasqualettoSemola} after \cite{Ambrosio02,AmbrosioBrueSemola19} that
\begin{equation}\label{eq:RepresentationPerimeter}
\Per(E,\cdot)=\haus^{N-1}|_{\mathcal F E}.
\end{equation}
Notice that the notion of perimeter that we are using does not charge the boundary of the space under
consideration, if any (see \cite{DePhilippisGigli18,MondinoKapovitch,BrueNaberSemola20} or the beginning of \cref{sec:IsoperimetricProblem} for background about boundaries of $\RCD(K,N)$ spaces $(X,\dist,\mathcal{H}^N)$). Indeed, by the very definition (see \cite[Definition 4.1]{AmbrosioBrueSemola19}), reduced boundary points for $E$ are regular points for the ambient space $X$, i.e., $\R^N$ is the unique tangent cone of $X$ at such points, while the boundary of the space is disjoint from the set of regular points.

Moreover, according to \cite[Proposition 4.2]{BPSGaussGreen}, 
\[
\mathcal F E=E^{(1/2)}=\bigg\{x\in X\;\bigg|\;\lim_{r\to 0}\frac{\haus^N(E\cap B_r(x))}{\haus^N(B_r(x))}=\frac{1}{2}\bigg\}\, ,
\quad\text{ up to }\mathcal H^{N-1}\text{-null sets}\, . 
\]
}
\medskip 

\subsubsection{Cylindrical splitting}
We will need a useful Kasue-type rigidity theorem \cite{Kasue83} on $\RCD(0,N)$ spaces following from the results in \cite{KettererKitabeppuLakzian, Ketterer21}.

\begin{theorem}[Cylindrical splitting]\label{thm:SplittingKetterer}
Let $(X,\dist,\haus^N)$ be an $\RCD(0,N)$ space with $N\ge 2$. Let $E\subset X$ be a bounded open set and denote $\Omega\eqdef X \setminus \overline{E}$. Suppose that $\Omega = u^{-1}(0,+\infty)$ for some locally Lipschitz proper harmonic\footnote{Harmonic means $\int\nabla u \cdot \nabla f=0$ for any Lipschitz function $f$ supported in $\Omega$.} function $u:\Omega\to \R$ such that $|\nabla u|=1$ almost everywhere.\\
Let $(\widetilde\Omega,\widetilde\dist)$ be the metric space given by the completion of $\Omega$ endowed with intrinsic distance induced by $\dist$, let $\widetilde\meas= \haus^N|_\Omega$, and let $(Y,\dist_{\rm int})$ be the metric space given by $\partial E$ endowed with intrinsic distance induced by $\dist$\footnote{In this setting, the distance between distinct connected components of $\widetilde\Omega$ or of $Y$ is infinite.}.\\
Then any connected component $\Omega_\alpha$ of $(\widetilde\Omega,\widetilde\dist, \widetilde\meas)$ is $\RCD(0,N)$ and it is isomorphic to either $(Y_\alpha\times [0,+\infty), \dist_{\rm int}\otimes \dist_{\rm eu}, \haus^N)$ or to $(Y_\alpha\times [0,D_\alpha], \dist_{\rm int}\otimes \dist_{\rm eu}, \haus^N)$ for some $D_\alpha>0$, where $Y_\alpha$ is a connected component of $Y$, and in particular $(Y_\alpha, \dist_{\rm int}, \haus^{N-1}_{\rm int})$\footnote{Here $\haus^{N-1}_{\rm int}$ is the $(N-1)$-dimensional measure induced by $\dist_{\rm int}$.} is $\RCD(0,N-1)$. Moreover, in the first case there is $a_\alpha>0$ such that $\widetilde\dist|_{\{u>a_\alpha\}}= \dist|_{\{u>a_\alpha\}}$.
\end{theorem}

\begin{proof}
The statement is a convenient reformulation of \cite[Theorem 4.11]{Ketterer21}, whose proof is based on \cite{KettererKitabeppuLakzian}, in the case of non-collapsed $\RCD(0,N)$ spaces. The very last assertion follows from the fact that, since $u$ is proper and $Y$ is bounded, one shows that for $a_\alpha$ large enough any geodesic in $X$ with endpoints in $\Omega_\alpha \cap \{u>a_\alpha\}$, for an unbounded connected component $\Omega_\alpha$, is contained in $\Omega_\alpha$, hence $\widetilde\dist|_{\{u>a_\alpha\}}= \dist|_{\{u>a_\alpha\}}$.
\end{proof}

We mention that a version of \cref{thm:SplittingKetterer} holds more generally for collapsed $\RCD(0,N)$ spaces $(X,\dist,\meas)$.

As pointed out to the authors by C. Ketterer, the statement of \cref{thm:SplittingKetterer} is optimal in the sense that, in the second part of the statement, it is not possible in general to say that $\widetilde\dist|_{\Omega_\alpha}= \dist|_{\Omega_\alpha}$, as the example of a one-ended round cylinder closed with a disk shows (see also the introduction of \cite{Ketterer21}).
\medskip

\subsubsection{Isoperimetric sets and profile on $\RCD$ spaces}
Let us record here a topological regularity result for isoperimetric sets which is a direct consequence of the arguments in \cite{MoS21}, and the results obtained in \cite{APPSa} and \cite{AntonelliPasqualettoPozzetta21}.  

The set $\mathcal{R}^E$ in the forthcoming statement is the set of points in $\partial E$ such that at least one blow-up is a half-space in $\mathbb R^N$, see \cite[Definition 6.19]{MoS21}. We also recall that with $E^{(t)}$ we denote the set of points of $E$ that have density $t$ with respect to $\mathcal{H}^N$, and $\partial^eE:=X\setminus (E^{(0)}\cup E^{(1)})$ is the {\em essential boundary} of $E$.

\begin{theorem}\label{thm:Regularity}
Let $(X, \dist , \mathcal{H}^N )$ be an $\RCD(K,N)$ space such that $\mathcal{H}^N(B_1(x))\geq v_0>0$ for every $x\in X$. Let $E$ be an isoperimetric set. Then $E^{(1)}$ is open and bounded, and the essential boundary $\partial^e E$ is equal to the topological boundary $\partial E^{(1)}$. We denote $\partial E:=\partial E^{(1)}=\partial ^eE$.

Let us assume $B_2(x_0)\cap \partial X=\emptyset$ for some $x_0\in X$. Then, for every $y\in \partial E \cap B_1(x_0)$, there exists $0<s<1$ such that the following holds. For every $\alpha\in (0,1)$, there exists a relatively open set $O_\alpha \subset \partial E\cap B_s(y)$  with $\mathcal{R}^E\cap B_s(y)\subset O_\alpha$ such that $O_\alpha$ is $\alpha$-biH\"older homeomorphic to an open, smooth $(N-1)$-dimensional manifold. 

In addition, we have that 
\[
\mathrm{dim} ((\partial E\setminus\mathcal{R}^E)\cap B_2(x_0))\leq N-3,
\]
where ${\rm dim}$ denotes Hausdorff dimension.

In particular, if $N=2$ and $\partial X=\emptyset$, the boundary $\partial E$ of every isoperimetric set $E$ is a $1$-dimensional (possibly disconnected) topological manifold without boundary.
\end{theorem}

The first part of \cref{thm:Regularity} is stated and proved in \cite[Theorem 1.3, Theorem 1.4]{AntonelliPasqualettoPozzetta21}.  Moreover, the second part of \cref{thm:Regularity} is the analog of \cite[Proposition 6.21]{MoS21}, which is a direct consequence of \cite[Theorem 6.8]{MoS21}. The statement of \cite[Theorem 6.8]{MoS21} for isoperimetric sets in $\RCD(K,N)$ spaces $(X,\dist,\mathcal{H}^N)$ follows verbatim from the proof of \cite[Theorem 6.8]{MoS21}, and the uniform $\Lambda$-minimality estimates in \cref{cor:UniformRegularityIsop}. Finally, the last part of \cref{thm:Regularity} is the analog of \cite[Theorem 6.29]{MoS21}, and can be obtained reasoning verbatim as in \cite{MoS21}, and using the density bounds and the uniform $\Lambda$-minimality estimates in \cref{cor:UniformRegularityIsop}.
\medskip

In the following, we will always assume without loss of generality that any isoperimetric set $E$ as in \cref{thm:Regularity} satisfies $E=E^{(1)}$, hence it is open, bounded and $\Per(E,\cdot)=\mathcal{H}^{N-1}\res \partial E$.
\medskip

We recall the forthcoming result, which has been proved in \cite[Proposition 3.1]{APPSb}, see also \cite[Proposition 3.11]{APPSa} for the last part of the statement.

\begin{theorem}[{\cite{APPSa, APPSb}}]\label{prop:C>0}
Let $(X,\dist,\mathcal{H}^N)$ be a non compact $\RCD(0,N)$ metric measure space for some $N\ge 2$, and assume that $\mathcal{H}^N(B_1(x))\ge v_0>0$ for every $x \in X$. Let $E\subset X$ be an isoperimetric set. Then there exists $c\in [0,\infty)$ such that, denoting by $f$ the signed distance function from $\overline{E}$, in such a way it is positive outside $E$ and negative inside $E$, it holds
\begin{equation}\label{eq:sharpLap0}
\mathbf\Delta f\le \frac{c}{1+\frac{c}{N-1}f}\, ,\quad\text{on $X\setminus \overline{E}$}\, ,\quad \mathbf\Delta f\ge  \frac{c}{1+\frac{c}{N-1}f}\, ,\quad\text{on $E$}\, ,
\end{equation}
where the bounds are understood in the sense of distributions, see \cref{def:MeasureLaplacian}.
Moreover, if $c=0$, we have that the rigidity stated in \cref{thm:SplittingKetterer} holds with $\Omega:=X\setminus\overline{E}$.\\ 
Finally, for every $t\ge 0$ it holds
\begin{equation}\label{eq:extareabd}
\Per(\{x\in X\, :\, \dist(x,\overline{E})\le t\})\le \left(1+\frac{c}{N-1}t\right)^{N-1} \Per(E).
\end{equation}\label{eqn:ExtPerEsterno}
and, for any $t\ge 0$,
\begin{equation}\label{eq:intareabound}
\Per(\{x\in X\, :\, \dist(x,X\setminus E)\le t\})\le \left(1-\frac{c}{N-1}t\right)^{N-1}\Per(E).
\end{equation}
\end{theorem}

\begin{remark}[Mean curvature barrier]\label{def:MeanCurvatureBarrier}
    We recall that any $c$ for which \eqref{eq:sharpLap0} holds is called {\em mean curvature barrier} for the isoperimetric set $E$, as it plays the role of a constant value for the mean curvature of the boundary of $E$. This terminology was introduced in \cite{APPSa}.
\end{remark}

Let us now recall a couple of results. The mass decomposition result in \cref{thm:MassDecompositionINTRO} is proved in \cite{AntonelliNardulliPozzetta} building on top of \cite{Nar14, AFP21, AntonelliPasqualettoPozzetta21}. The result in \cref{cor:UniformRegularityIsop} can be found in \cite[Corollary 4.17]{APPSa}.

\begin{theorem}[{Asymptotic mass decomposition \cite{AntonelliNardulliPozzetta}}]\label{thm:MassDecompositionINTRO}
Let $K\in\mathbb R$, and $N\geq 1$. Let $(X,\dist,\mathcal{H}^N)$ be a non compact $\RCD(K,N)$ space. Assume there exists $v_0>0$ such that $\mathcal{H}^N(B_1(x))\geq v_0$ for every $x\in X$. Let $V>0$. For every minimizing (for the perimeter) sequence $\Omega_i\subset X$ of volume $V$, with $\Omega_i$ bounded for any $i$, up to extract a subsequence, there exist an increasing and bounded sequence $\{N_i\}_{i\in\mathbb N}\subset \mathbb N$, disjoint finite perimeter sets $\Omega_i^c, \Omega_{i,j}^d \subset \Omega_i$, and points $p_{i,j}$, with $1\leq j\leq N_i$ for any $i$, such that
\begin{itemize}
    \item $\lim_{i} \dist(p_{i,j},p_{i,\ell}) = \lim_{i} \dist(p_{i,j},o)=\infty$, for any $j\neq \ell\leq \overline N$ and any $o\in X$, where $\overline N:=\lim_i N_i <\infty$;
    \item $\Omega_i^c$ converges to $\Omega\subset X$ in the sense of finite perimeter sets, and we have $\mathcal{H}^N(\Omega_i^c)\to_i \mathcal{H}^N(\Omega)$, and $ \Per( \Omega_i^c) \to_i \Per(\Omega)$. Moreover $\Omega$ is a bounded isoperimetric region for its own volume in $X$;
    \item for every $j\leq\overline N$, $(X,\dist,\mathcal{H}^N,p_{i,j})$ converges in the pmGH sense  to a pointed $\RCD(K,N)$ space $(X_j,\dist_j,\mathcal{H}^N,p_j)$. Moreover there are isoperimetric regions $Z_j \subset X_j$ such that $\Omega^d_{i,j}\to_i Z_j$ in $L^1$-strong and $\Per(\Omega^d_{i,j}) \to_i \Per (Z_j)$;
    \item it holds that
    \begin{equation}\label{eq:UguaglianzeIntro}
    I_{(X,\dist,\mathcal{H}^N)}(V) = \Per(\Omega) + \sum_{j=1}^{\overline{N}} \Per (Z_j),
    \qquad\qquad
    V=\mathcal{H}^N(\Omega) +  \sum_{j=1}^{\overline{N}} \mathcal{H}^N(Z_j).
    \end{equation}
\end{itemize}
\end{theorem}

\begin{corollary}[{\cite[Corollary 4.17]{APPSa}}]\label{cor:UniformRegularityIsop}
    Let $0<V_1<V_2<V_3$, and let $K \in \R, N\geq 2, v_0>0$. Then there exist $\Lambda, R>0$ depending on $K,N,v_0,V_1,V_2, V_3$ such that the following holds.

If $(X,\dist,\haus^N)$ is an $\RCD(K,N)$ space with $\inf_{x\in X}\haus^N(B_1(x))\geq v_0>0$, $\haus^N(X)\ge V_3$, and $E\subset X$ is an isoperimetric region with $\haus^N(E) \in [V_1,V_2]$, then $E$ is a $(\Lambda,R)$-minimizer, i.e., for any $F\subset X$ such that $F \Delta E \subset B_R(x)$ for some $x \in X$, then
\[
\Per(E) \le \Per(F) + \Lambda\haus^N(F\Delta E)\, .
\]
The constant $\Lambda$ can be taken to be equal to the Lipschitz constant of the isoperimetric profile $I$ on $[V_1/2,(V_2+V_3)/2]$, and $R>0$ can be taken to be a radius such that $\haus^N(B_R(x)) \le \min\{V_1/2,(V_3-V_2)/2 \}$ for every $x \in X$.

Moreover there are $C_1\in (0,1),C_2>0,R'$ depending on $K,N,v_0,V_1,V_2, V_3$ such that
\begin{equation}\label{eqn:DensityEstimates}
C_1 \le \frac{\haus^N(B_r(x) \cap E)}{\haus^N(B_r(x))}\le 1-C_1\, ,
\qquad
C_2^{-1} \le \frac{\Per(E,B_r(x))}{r^{N-1}} \le C_2\, ,
\end{equation}
for any $x \in \partial E$ and any $r\in(0,R']$.
\end{corollary}

We finally recall some basic properties of the isoperimetric profile and the asymptotic mass decomposition on $N$-dimensional $\RCD(0,N)$ spaces. The following \cref{cor:IsoperimetricProfileRCD0N} is obtained in \cite[Theorem 3.8]{APPSb}, specializing the study on the differential properties of the isoperimetric profile in the case of nonnegative curvature. The result in \cref{lem:IsoperimetricAtFiniteOrInfinite} is proved in \cite[Lemma 4.19 and Corollary 4.20]{APPSa}.

\begin{theorem}[{\cite[Theorem 3.8]{APPSb}}]\label{cor:IsoperimetricProfileRCD0N}
Let $N\geq 1$, and let $(X,\dist,\mathcal{H}^N)$ be an $\RCD(0,N)$ space with isoperimetric profile function $I$. Then the function $I^{\frac{N}{N-1}}$ is concave on $(0,\mathcal{H}^N(X))$. Then, a fortiori, $I$ is concave and strictly subadittive on $(0,\mathcal{H}^N(X))$. Finally, if $\mathcal{H}^N(X)=+\infty$, $I$ is nondecreasing on $(0,+\infty)$.
\end{theorem}

\begin{proposition}[{\cite[Lemma 4.19 and Corollary 4.20]{APPSa}}]\label{lem:IsoperimetricAtFiniteOrInfinite}
Let $(X,\dist,\mathcal{H}^N)$ be a noncompact $\RCD(0,N)$ space with $N\geq 2$. Let us assume that $\mathcal{H}^N(B_1(x))\geq  v_0>0$ for every $x\in X$.  Let $\{\Omega_i\}_{i\in\mathbb N}$ be a minimizing (for the perimeter) sequence of bounded finite perimeter sets of volume $V$ in $X$. Then, if one applies \cref{thm:MassDecompositionINTRO}, either $\overline N=0$, or $\overline N=1$ and $\mathcal{H}^N(\Omega)=0$.

In particular, for any $V>0$ there is an $\RCD(0,N)$ space $(Y,\dist,\mathcal{H}^N)$ which is either $X$ or a pmGH limit of $X$ along a diverging sequence $\{p_i\}_i \subset X$, and a set $E \subset Y$ such that $\mathcal{H}^N(E)=V$ and $I_X(V)=\Per(E)$.

Finally, let $c$ be any barrier given by \cref{prop:C>0} applied to either $\Omega$ or $Z_1$. Then
\begin{equation}\label{eqn:MeanCurvatureInequality}
I_+'(v)\leq c \leq I_-'(v),
\end{equation}
where $I_+'$ and $I_-'$ denote the right and left derivative, respectively, which exist due to \cref{cor:IsoperimetricProfileRCD0N}.
\end{proposition}

\medskip

We conclude this part by stating a useful lemma on the connected components of complements of isoperimetric sets, consequence of \cref{cor:IsoperimetricProfileRCD0N}. First, we recall the definition of {\em (number of) ends} of a space.

\begin{definition}
We say that a metric space $(X,\dist)$ has $k\in\N$ \emph{ends} with $k\ge 1$ if
\begin{itemize}
\item for every compact set $K\subset X$ the set $X\setminus K$ has at most $k$ unbounded connected components;
\item there exists a compact set $K'\subset X$ such that the set $X\setminus K'$ has exactly $k$ unbounded connected components.
\end{itemize}
\end{definition}

\begin{lemma}\label{lem:ConnCompRCD}
Let $(X,\dist,\haus^N)$ be a noncompact $\RCD(0,N)$ space such that $\haus^N(B_1(x))\ge v_0>0$ for any $x\in X$. Let $E\subset X$ be an isoperimetric set. Then the connected components of $X\setminus \overline{E}$ are unbounded. In particular, if $X$ has one end, $X\setminus \overline{E}$ is connected.
\end{lemma}

\begin{proof}
If by contradiction $\Omega$ is a bounded connected component of $X\setminus \overline{E}$, then the set $F\eqdef E\cup \Omega$ satisfies $\haus^N(F)>\haus^N(E)$ and $\Per(F)<\Per(E)$, which contradicts the fact that the isoperimetric profile is nondecreasing on noncompact {non-collapsed} $\RCD(0,N)$ spaces ({i.e., the reference measure is $\haus^N$}), see \cref{cor:IsoperimetricProfileRCD0N}.
\end{proof}

\medskip
\subsection{On the asymptotic geometry of $\mathrm{CBB}(0)$ spaces}
\label{sec:AsympGeometry}

We collect a series of lemmas about the geometry at infinity along $\CBB(0)$ spaces, focusing on the case of spaces having a pointed Gromov--Hausdorff limit along a diverging sequence equal to a cylinder over a compact space.

\begin{lemma}\label{Lemma:KeyEstimate}
	Let $(X,\dist)$ be a noncompact $N$-dimensional $\CBB(0)$ metric space with $N\geq 2$. Let $\{p_i\}$ be a sequence diverging at infinity on $X$. Let $\gamma:[0,+\infty)\to X$ be a ray from $o\in X$ such that
	\[
	\lim_{i\to +\infty}\frac{\dist(p_i,\gamma)}{|op_i|}=0,
	\qquad
	\liminf_{i\to+\infty} \dist(p_i,\gamma) >0.
	\]
	Let $t_i\geq 0$ be such that $\dist(p_i,\gamma)=|p_i\gamma(t_i)|$\footnote{Notice that, from the hypothesis, $t_i>0$ for $i$ large enough.}, let $\alpha_i$ be a unit-speed geodesic from $p_i$ to $o$, and let $\beta_i$ be a unit-speed geodesic from $p_i$ to $\gamma(\tau_i)$, where $\tau_i\geq t_i$ is chosen such that $|p_io|=|p_i\gamma(\tau_i)|$. 
	
	Then, for every $\varepsilon>0$ there exists $i_0\in\mathbb N$ such that for every $i\geq i_0$ the following hold:
	\begin{equation}\label{eq:StimaEpsilonLinea}
	r^2+s^2+2rs(1-\varepsilon)\leq |\alpha_i(r)\beta_i(s)|^2, \qquad \text{for all $r,s<|op_i|$},
	\end{equation}
	\begin{equation}\label{eq:StimaEpsilonAngoliLinea}
	\left|\angle (p_i)^o_{\gamma(t_i)} -\frac{\pi}{2}\right|+\left|\angle (p_i)^{\gamma(t_i)}_{\gamma(\tau_i)} -\frac{\pi}{2}\right|+\left|\angle (p_i)^o_{\gamma(\tau_i)} -\pi\right|\leq \varepsilon .
	\end{equation}
\end{lemma}

\begin{proof}
	By triangle inequality we estimate
	\begin{equation*}
		\begin{split}
		\tau_i &= \tau_i - t_i + |\gamma(t_i)o| \ge |p_i\gamma(\tau_i)| - |p_i \gamma(t_i)| + |o p_i |  - |p_i \gamma(t_i)| = 2 |op_i| - 2 |p_i \gamma(t_i)|,
		\end{split}
	\end{equation*}
	and since $ |op_i| - |p_i \gamma(t_i)| = \dist(\gamma(0),p_i) - \dist(p_i,\gamma) \ge 0$, then
	\begin{equation}\label{eq:zzTaui}
		\tau_i^2 \ge 4( |op_i| -  |p_i \gamma(t_i)|)^2.
	\end{equation}
	By monotonicity we get
	\[
	\begin{split}
	\frac{r^2 + s^2- |\alpha_i(r)\beta_i(s)|^2}{2rs} 
	&\le \frac{2|op_i|^2 - |\gamma(0)\gamma(\tau_i)|^2}{2|op_i|^2} 
	= \frac{2|op_i|^2 - \tau_i^2}{2|op_i|^2}
	\le \frac{2|op_i|^2 - 4(|op_i|-|p_i\gamma(t_i)|)^2}{2|op_i|^2} \\
	&= -1 -2 \frac{\dist(p_i,\gamma)^2}{|op_i|^2} + 4 \frac{\dist(p_i,\gamma)}{|op_i|},
	\end{split}
	\]
	for any $r,s \in(0,|op_i|)$. As the right hand side tends to $-1$ as $i\to+\infty$, estimate \eqref{eq:StimaEpsilonLinea} is proved.
	
	Now denote $\sigma_i=[p_i \gamma(t_i)]$. Up to subsequence, let $(X,\dist,p_i)$ converge in pGH sense to a limit space $(Y,\dist_\infty,p_\infty)$. 
	By monotonicity, hinge condition, and \eqref{eq:zzTaui}, exploiting the hypotheses we find
	\[
	\begin{split}
	    \liminf_i \angle (p_i)^o_{\gamma(t_i)} & \ge \liminf_i \arccos\left( \frac{|p_i\gamma(t_i)|^2 +|op_i|^2 - t_i^2}{2|p_i\gamma(t_i)|\,|op_i|}\right) \\&\ge \liminf_i \arccos\left( \frac{|p_i\gamma(t_i)|^2 + |p_i\gamma(t_i)|^2 + t_i^2  - t_i^2}{2|p_i\gamma(t_i)|\,|op_i|}\right) = \frac{\pi}{2},
	\end{split}
	\]
	\[
	\begin{split}
	    \liminf_i \angle (p_i)^{\gamma(t_i)}_{\gamma(\tau_i)}
	    &\ge \liminf_i  \arccos\left( \frac{|p_i\gamma(t_i)|^2 +|p_i\gamma(\tau_i)|^2 - (\tau_i-t_i)^2}{2|p_i\gamma(t_i)|\,|op_i|} \right) \\
	    &\ge \liminf_i \arccos\left( \frac{|p_i\gamma(t_i)|^2 +|p_i\gamma(t_i)|^2 + (\tau_i-t_i)^2 - (\tau_i-t_i)^2}{2|p_i\gamma(t_i)|\,|op_i|} \right) = \frac{\pi}{2},
	\end{split}
	\]
	\[
	\begin{split}
	    \liminf_i \angle (p_i)^o_{\gamma(\tau_i)} &\ge \liminf_i \arccos\left( \frac{|op_i|^2 + |p_i\gamma(\tau_i)|^2 - \tau_i^2}{2|op_i|\, |p_i\gamma(\tau_i)|} \right) \\
	    &\overset{\eqref{eq:zzTaui}}{\ge} \liminf_i \arccos\left( \frac{2|op_i|^2 - 4(|op_i|-|p_i\gamma(t_i)|)^2}{2|op_i|^2} \right) = \pi.
	\end{split}
	\]
	On the other hand, by four point comparison we have that
	\[
	\angle^0 \bar \alpha_i(t) \bar p_i \bar \sigma_i(t) + 
	\angle^0 \bar \sigma_i(t) \bar p_i \bar \beta_i(t) + 
	\angle^0 \bar \alpha_i(t) \bar p_i \bar \beta_i(t)  \le 2\pi,
	\]
	for any $i$ and $t>0$ sufficiently small. Letting $t\to0$ this implies
	\[
	\angle (p_i)^{\gamma(t_i)}_{\gamma(\tau_i)}
	+ \angle (p_i)^{\gamma(t_i)}_{\gamma(\tau_i)}
	+ \angle (p_i)^o_{\gamma(\tau_i)} \le 2\pi,
	\]
	for any $i$. Taking into account the previous lower bounds on the angles $\angle (p_i)^{\gamma(t_i)}_{\gamma(\tau_i)}$, $\angle (p_i)^{\gamma(t_i)}_{\gamma(\tau_i)}$, $\angle (p_i)^o_{\gamma(\tau_i)} $, estimate \eqref{eq:StimaEpsilonAngoliLinea} follows.
\end{proof}

We mention that the following \cref{lem:EveryLimitAtInfinitySplits} generalizes to arbitrary dimensions and arbitrary diverging sequences the result in \cite[Lemma 2.12]{PositiveScalarRicci}, that appeared on arXiv while the present paper was being completed. Later, we found that this property had already observed in a work of Kleiner--Lott, \cite[Page 2852, Appendix G]{KleinerLott}.

\begin{lemma}\label{lem:EveryLimitAtInfinitySplits}
Let $(X,\dist)$ be a noncompact $N$-dimensional $\CBB(0)$ metric space with $N\geq 2$. Let $\{p_i\}$ be a sequence diverging at infinity on $X$. Then every pGH limit of $(X,\dist,p_i)$ splits a line. 
\end{lemma}

\begin{proof}
Let us fix a base point $o\in X$. Let us first prove the following claim: there exists a ray $\gamma:[0,+\infty)\to X$ such that $\gamma(0)=o$ such that
\begin{equation}\label{eqn:Claim1}
\lim_{i\to +\infty}\frac{d(p_i,\gamma)}{|op_i|}=0.
\end{equation}
Indeed, let $\gamma_i:[0,|op_i|]\to X$ be a segment from $o$ to $p_i$, and let $\gamma:[0,+\infty)\to X$ be the local uniform limit of $\gamma_i$ as $i\to+\infty$. Thus, $\gamma$ is a unit-speed ray from $o$. Let us now prove \eqref{eqn:Claim1}. 

Fix $\varepsilon>0$ and take $i_0(\varepsilon)>1$ such that for all $i> i_0$ we have $|op_i|>1$, and $|\gamma_i(1)\gamma(1)|<\varepsilon$. Let us fix $i>i_0$. We can now estimate 
\begin{equation}\label{eqn:UnaStima1}
    \frac{\dist(p_i,\gamma)}{|op_i|}\leq \frac{|p_i\gamma(|op_i|)|}{|op_i|} = \frac{|\gamma_i(|op_i|)\gamma(|op_i|)|}{|op_i|}\leq |\gamma_i(1)\gamma(1)| < \varepsilon,
\end{equation}
where in the second-last inequality we used $|op_i|>1$ and the fact that for any two geodesics $\alpha,\beta:[0,\ell]\to X$ with $\alpha(0)=\beta(0)$ we have that
\[
[0,\ell]\ni t \mapsto \frac{|\alpha(t)\beta(t)|}{t}
\]
is nonincreasing in $t$ as a result of item (2) of \cref{def:CBB}. Thus \eqref{eqn:UnaStima1} gives \eqref{eqn:Claim1}.
\medskip

Let us now conclude the proof of the Lemma. Let $(Y,\dist_Y,y)$ denote the pGH limit of $(X,\dist,p_i)$, up to subsequence. Let us consider the ray $\gamma$ constructed as in the first part of the proof. For every $i$ sufficiently big, let us take $t_i$ such that $\dist(p_i,\gamma)=|p_i\gamma(t_i)|$, and $\tau_i>t_i$ such that $|p_i\gamma(\tau_i)|=|op_i|$. Let us now distinguish two cases:
\begin{itemize}
\item If $\liminf_{i\to + \infty}\dist(p_i,\gamma)=0$, thus the pGH limit of $(X,\dist,\gamma(t_i))$ is isometric to $Y$. Notice that $t_i\to +\infty$. Since $\gamma$ is a ray, if we consider $\tilde\gamma_i(t):[-t_i,+\infty)\to X$ defined by $\tilde\gamma_i(t):=\gamma_i(t+t_i)$, we get that it converges locally uniformly to a line in $Y$. Thus, applying the splitting theorem, we get that $Y$ splits a line.
\item If $\liminf_{i\to + \infty}\dist(p_i,\gamma)>0$, we can apply \cref{Lemma:KeyEstimate}. Let $\alpha_i,\beta_i$ be defined as in \cref{Lemma:KeyEstimate}. Let $\alpha_\infty,\beta_\infty$ the local uniform limits of $\alpha_i,\beta_i$ in $Y$. Notice that $\alpha_\infty,\beta_\infty$ are two rays starting from $y\in Y$. Thus, for every $t>0$ we have $|\alpha_\infty(t)\beta_\infty(t)|\leq 2t$. On the other hand, we can apply \eqref{eq:StimaEpsilonLinea} to get that, for every $t>0$, 
\[
|\alpha_\infty(t)\beta_\infty(t)|\geq 2t.
\]
Thus the curve $\delta:\mathbb R\to Y$ such that $\delta(t)=\beta_\infty(t)$ for $t\geq 0$ and $\delta(t)=\alpha_\infty(-t)$ for $t\leq 0$ is a line, and $Y$ splits thanks to the splitting theorem.
\end{itemize}
\end{proof}

\begin{remark}\label{rem:RCD0NNonsplitta}
    The previous \cref{lem:EveryLimitAtInfinitySplits} fails in the class of $\RCD(0,N)$ metric measure spaces. The example in \cite[pp. 913--914]{KasueWashio}, see also the discussion in \cite[Section 1.4]{AntBruFogPoz}, is a complete Riemannian manifold with $\mathrm{Ric}\geq 0$ such that it does not split a line, and there exists a pointed Gromov Hausdorff limit at infinity that is isometric to the same space. Another example is given by the Grushin-halfplane with a properly chosen measure that turns it into an $\RCD(0,N)$ space, see \cite[Remark 3.9]{GrushinHalf}.
\end{remark}

\begin{lemma}\label{lem:UppBoundDistDaGamma}
Let $(X,\dist)$ be a noncompact $N$-dimensional $\CBB(0)$ metric space with $N\geq 2$. Let $\{p_i\}$ be a sequence diverging at infinity on $X$. Let $\gamma:[0,+\infty)\to X$ be a ray emanating from $o\in X$ such that
\[
\lim_{i\to +\infty}\frac{\dist(p_i,\gamma)}{|op_i|}=0.
\]
Let us assume that $(X,p_i)$ pGH-converges to $(\mathbb R\times K,(0,k_0))$ for some compact $K$, with $k_0\in K$.
Then 
\[
\limsup_{i\to +\infty}\dist(p_i,\gamma)\leq\mathrm{diam}(K).
\]
\end{lemma}

\begin{proof}
We can assume that $\ell:=\limsup_i \dist(p_i,\gamma)= \lim_i \dist(p_i,\gamma)>0$. Let $t_i\geq 0$ be such that $\dist(p_i,\gamma)=|p_i\gamma(t_i)|$, and $\tau_i\ge t_i$ such that $|p_io|=|p_i\gamma(\tau_i)|$. Let $\alpha_i=[p_io]$, $\beta_i=[p_i\gamma(\tau_i)]$, and $\sigma_i=[p_i\gamma(t_i)]$. Passing to the pointed limit, $\alpha_i,\beta_i$ converge to rays from $(0,k_0)$ whose union is a line. Exploiting \cref{rem:AngleContinuity} and \cref{Lemma:KeyEstimate}, one sees that the limit $\sigma_\infty$ of $\sigma_i$ is orthogonal to such line. Hence applying \cref{prop:GeodeticheProdotto} one easily gets that the length of $\sigma_\infty$ is no more than $\diam\,K$, implying the claim. 
\end{proof}

\begin{corollary}\label{cor:RettaGenerataDaSuccessione}
Let $(X,\dist)$ be a noncompact $N$-dimensional $\CBB(0)$ metric space with $N\geq 2$. Let $\{p_i\}$ be a sequence diverging at infinity on $X$ such that $(X,p_i)$ pGH-converges to $(\mathbb R\times K,(0,k_0))$ for some compact $K$, with $k_0\in K$. Fix $o \in X$.

Then the ray $\gamma:[0,+\infty)\to X$ emanating from $o$ given by the limit, up to subsequence, of the geodesics $[op_i]$ satisfies $\lim_i \dist(p_i,\gamma)/|op_i|=0$ and there exists a sequence $t_i\to+\infty$ such that $(X,\gamma(t_i))$ pGH-converges to $\R\times K$.
\end{corollary}

\begin{proof}
The fact that $\lim_i \dist(p_i,\gamma)/|op_i|=0$ has been already observed in \eqref{eqn:Claim1}. Hence the full statement follows by applying \cref{lem:UppBoundDistDaGamma}.
\end{proof}

\begin{lemma}\label{lem:NoLinesOneEnd}
Let $N\geq 2$, and let $(X,\dist)$ be an $N$-dimensional noncompact $\CBB(0)$ metric space. 
Let $p_i$ be a diverging sequence of points on $X$ such that $(X,p_i)$ converges in the pGH topology to $\mathbb R\times K$, where $K$ is compact. 

Then $X$ splits no lines if and only if $X$ has one end.
\end{lemma}

\begin{proof}
If $X$ splits no lines, then obviously $X$ has one end. So assume $X$ has one end. Without loss of generality we can rewrite the space as $X=\R^\ell \times Y$, where $Y$ does not split lines, for some integer $\ell \ge 0$. We aim at showing that $\ell=0$.\\
So assume by contradiction that $\ell\ge1$. 
The points $p_i$ have the form $p_i=(v_i,y_i)$. Since $(v,y)\mapsto (v-v_i,y)$ is an isometry, we have that $(X,(0,y_i))$ converges to $\R\times K$, so let us assume that $v_i=0$ for any $i$. If $\{y_i\}$ remains in a compact subset of $Y$, then $X$ is isometric to $\R\times K$, and thus it has two ends, that is impossible. Hence $Y$ is noncompact and $y_i$ diverges along $Y$, up to subsequence.\\
Fix a point $o=(0,y_0)$ and let $\gamma,t_i$ be given by \cref{cor:RettaGenerataDaSuccessione}. Hence, from \cref{prop:GeodeticheProdotto}, $\gamma=(0,\sigma)$ and $\sigma$ is a unit-speed geodesic in $Y$ from $y_0$. Pick $w\in\R^\ell$ with unit norm and consider $\alpha(t)\eqdef(tw,y_0)$. Hence $|\alpha(t_i)\gamma(t_i)|^2 = 2t_i^2$. Define the curve $\beta_i:[0,\sqrt{2}\,t_i]\to X$ by
\[
\beta_i(s) \eqdef \left(\frac{s}{\sqrt{2}\,t_i} t_i w , \sigma \left(\left(  1- \frac{s}{\sqrt{2}\,t_i} \right) t_i \right) \right).
\]
Hence $\beta_i(0)=\gamma(t_i)$, $\beta_i(\sqrt{2}\,t_i)=\alpha(t_i)$, and $|\beta_i'|=1$ almost everywhere, that is, $\beta_i$ is a geodesic from $\gamma(t_i)$ to $\alpha(t_i)$. Up to subsequence, the lines $t\mapsto \gamma(t-t_i)$ converge to a line $\eta$ through some $(0,k_0) \in \R\times K$ along $(X,\gamma(t_i))$, and the geodesics $\beta_i$ converge to a ray $\beta$ from $(0,k_0)$.\\
An immediate computation shows that $\angle \gamma(t_i)^o_{\beta_i(t_i)}= \tfrac{\pi}{4}$ for any $i$. Moreover such angles pass to the limit by \cref{rem:AngleContinuity}, hence $\angle (0,k_0)^{\eta(-1)}_{\beta(1)}= \tfrac\pi4$. But \cref{prop:GeodeticheProdotto} implies that $\beta(t)=(t,k_0)$ or $\beta(t)=(-t,k_0)$ for any $t\ge0$, which implies $\angle (0,k_0)^{\eta(-1)}_{\beta(1)} \in \{0,\pi\}$, giving a contradiction.
\end{proof}

\begin{lemma}\label{lem:ConoHalfline}
Let $(X,\dist)$ be a noncompact $N$-dimensional $\CBB(0)$ metric space with $N\geq 2$. Let $\gamma:[0,+\infty)\to X$ be a ray emanating from $o\in X$ such that, for a sequence $t_i\to+\infty$, $(X,\gamma(t_i))$ pGH-converges to $(\mathbb R\times K,(0,k_0))$ for some compact $K$, with $k_0\in K$.

If $X$ has one end, then the asymptotic cone to $X$ is isometric to a halfline $([0,+\infty), \dist_{\rm eu})$ endowed with Euclidean distance $\dist_{\rm eu}$.
\end{lemma}

\begin{proof}
We want to show that the ideal boundary $X_\infty$ of $X$ contains only a point. Assume by contradiction that there is a ray $\sigma$ from $o$ such that $\lim_{t\to+\infty} |\gamma(t)\sigma(t)|/t = \delta >0$.

We claim that there exists $\eps>0$ and $i_0\in \N$ such that
\begin{equation}\label{eq:zw}
     |\gamma(t_i)\sigma(t_i)|<(1-\eps)2t_i \qquad \forall\,i\ge i_0.
\end{equation}
Indeed, suppose instead that there is $\varepsilon_j\searrow 0$ and a subsequence $i_j$ such that $|\gamma(t_{i_j})\sigma(t_{i_j})|\ge(1-\eps_j)2t_{i_j}$ for any $j$. The sequence of blow-downs $(X,\dist_j=t_{i_j}^{-1}\dist,o)$ pGH-converges to the cone $(C(X_\infty),\dist_C,o_C)$ over $X_\infty$. The curves $\gamma_j(t)\eqdef \gamma(t_{i_j}\, t)$, $\sigma_j(t)\eqdef \sigma(t_{i_j}\, t)$
are rays from $o$ in $(X,\dist_j)$ and, up to subsequence, converge to limit rays $\gamma_\infty$, $\sigma_\infty$ in the asymptotic cone. Moreover $(1-\eps_j)2 \le \dist_j(\gamma_j(1),\sigma_j(1)) \le 2$, which implies that $\dist_C(\gamma_\infty(1),\sigma_\infty(1))=2$. Thus $\delta=2$, and we get that for every $t>0$, $\dist_C(\gamma_\infty(t),\sigma_\infty(t))=2t$. Thus $\eta:\R\to C(X_\infty)$ given by $\eta(t)=\gamma_\infty(t)$ for $t\ge 0$ and $\eta(t)=\sigma_\infty(-t)$ for $t<0$ is a line in $C(X_\infty)$. Then \cref{thm:AsymptoticCones} implies that $X$ splits a line, but, since $X$ is assumed to have one end, this contradicts \cref{lem:NoLinesOneEnd}.

Consider now a sequence of geodesics $\alpha_i$ from $\gamma(t_i)$ to $\sigma(t_i)$. Observe that the length $L(\alpha_i)=|\gamma(t_i)\sigma(t_i)|\sim \delta t_i \to +\infty$ by absurd hypothesis, as $i\to\infty$.

Let $\triangle \bar o \bar \gamma(t_i) \bar \sigma(t_i)$ be a comparison triangle for $\triangle  o  \gamma(t_i) \sigma(t_i)$. By \eqref{eq:zw} there is $\eps_1>0$ such that
\begin{equation}\label{eq:zzz}
\angle \gamma(t_i)^o_{\sigma(t_i)} \ge \angle^0 \bar \gamma(t_i)^{\bar o}_{\bar\sigma(t_i)} \overset{\eqref{eq:zw}}{\ge} \eps_1 >0
 \qquad \forall\,i\ge i_0,
\end{equation}
where the first inequality follows by angle condition.

On the other hand, we claim also that there is $\eps_2<\pi$ such that 
\begin{equation}\label{eq:zzz1}
    \angle \gamma(t_i)^o_{\sigma(t_i)} \le \eps_2
\end{equation}
for any sufficiently large $i$.

Indeed, suppose that, up to subsequence, $\lim_i \angle \gamma(t_i)^o_{\sigma(t_i)} =\pi$. Consider the sequence of blow-downs $(X,\dist_i\eqdef t_i^{-1} \dist, o)$, which converges to the asymptotic cone $C(X_\infty)$ of $X$. Similarly as above, the curves $\gamma_i(t)\eqdef \gamma(t_i\, t)$, $\sigma_i(t)\eqdef \sigma(t_i\, t)$ converge in the asymptotic cone to limit rays $\gamma_\infty$, $\sigma_\infty$  from a tip $o_C$. Also the geodesics $\tilde\alpha_i(t)\eqdef\alpha_i(t_i\,t)$ converge to a geodesic $\tilde\alpha_\infty$ from $\gamma_\infty(1)$ to $\sigma_\infty(1)$. Observe that $L(\tilde\alpha_\infty)=\delta>0$ and then $\tilde\alpha_\infty$ is non-constant. By \cref{rem:AngleContinuity} the angles $\angle \gamma(t_i)^o_{\sigma(t_i)}$, which are invariant under the rescaling of the distance, pass to the limit, therefore $\angle \gamma_\infty(1)^{o_C}_{\sigma_\infty(1)} =\pi$. But this contradicts \cref{lem:AngoliRaggiCrossSection}.

Now we want to show that \eqref{eq:zzz}, \eqref{eq:zzz1} lead to a contradiction. Up to subsequence, along the sequence $(X,\gamma(t_i))$ the curve $\alpha_i$ (resp. $\gamma_i(t)\eqdef \gamma(t-t_i)$) converge locally uniformly to a ray $\alpha_\infty$ from $(0,k_0)$ (resp. to a line $\gamma_\infty$ through $(0,k_0)$) in $\R\times K$. Up to subsequence, by \cref{rem:AngleContinuity}, angles pass to the limit and then $0<\eps_1 \le \angle (0,k_0)^{\gamma_\infty(-1)}_{\alpha_\infty(1)} \le \eps_2<\pi$, which contradicts the third bullet of \cref{prop:GeodeticheProdotto}.

Therefore we conclude that $\delta=0$, that is, every ray is equivalent to $\gamma$ and $X_\infty$ is a point.
\end{proof}

\begin{remark}\label{rem:SormaniResults}
    In the setting of manifolds with nonnegative Ricci curvature and linear volume growth, it is known that the diameters of the level sets of any Busemann function grow sublinearly \cite[Theorem 1]{SormaniCAG}, and the asymptotic cone can be either a half-line or a line, see also \cite{SormaniJDG2000}. Moreover, in the same setting, an almost splitting theorem for cylinders on level sets of a Busemann function is proved in \cite[Theorem 34]{SormaniCAG}, where the GH-closeness is in term of the diameters of the level sets of the Busemann function. For statements analogous to the previous ones on $\RCD$ spaces, see \cite{Tao18, Tao20}.
    {After the recent \cite{XingyuZhu}, it is known that the conclusion of \cref{lem:ConoHalfline} holds on one-ended smooth manifolds with nonnegative Ricci curvature and linear volume growth.} 
\end{remark}

\section{Structure of asymptotically cylindrical $\mathrm{CBB}(0)$ spaces}\label{sec:Structure}

{

\begin{definition}\label{def:MeasureLaplacian}
    Let $(X,\dist)$ be a noncompact $N$-dimensional $\CBB(0)$ metric space. Let $F:X\to \R$ be a locally Lipschitz function
    and let $\nu$ be a Radon measure on $X$. We say that $F$ has \emph{measure Laplacian} if there exists a Radon measure $\mu$ on $X$ such that for any Lipschitz function $f$ on $X$ with bounded support, there holds
    \[
    \int \nabla F \cdot \nabla f \de \haus^N = - \int f \de \mu.
    \]
    In such a case we write $\boldsymbol{\Delta}F = \mu$.\\
    We say that $\boldsymbol{\Delta}F \ge \nu$ if $F$ has measure Laplacian and
    \[
    -\int \nabla F \cdot \nabla f \de \haus^N \ge \int f \de \nu,
    \]
    for any nonnegative Lipschitz function $f$ on $X$ with bounded support.
\end{definition}

\begin{definition}[Busemann function]\label{def:Busemann}
    Let $(X,\dist)$ be a noncompact $N$-dimensional $\CBB(0)$ metric space. Let $\gamma:[0,+\infty)\to X$ be a ray emanating from $o\in X$. For any $t\ge0$ we define $F_t:X\to \R$ by $F_t(x)\eqdef t- \dist(x,\gamma(t))$.\\
    By triangle inequality, $F_t$ is $1$-Lipschitz and $F_t\le F_s$ for any $t\le s$. The locally uniform limit $F(x)=\lim_{t\to+\infty} F_t(x)$ is called the \emph{Busemann function associated to $\gamma$}.
\end{definition}

\begin{remark}\label{rem:BasicPropertiesBusemann}
Let $(X,\dist)$ be a noncompact $N$-dimensional $\CBB(0)$ and let $F$ be the Busemann function associated to some ray $\gamma$. Then
\begin{itemize}
    \item $\lip  F (x)=1$ for any point $x\in X$,
    
    \item $F$ has measure Laplacian $\boldsymbol{\Delta}F\ge0$.
    
    \item {For any $r \le R$ there exists a $1$-Lipschitz surjective map $\Psi^R_r: \{F=R\} \to \{F=r\}$.}
\end{itemize}
Indeed $\lip  F \le 1$ being $F$ the uniform limit of $1$-Lipschitz functions. Also let $z_{t,\eps}\eqdef [x\gamma(t)](\eps)$, for $t,\eps>0$, and $x\in X$. On a sequence $t_n\to+\infty$ we have that $z_{t_n,\eps} \to z_\eps \in \partial B_\eps (x)$ and
\[
\frac{|F(x)-F(z_\eps)|}{|xz_\eps|}
= \lim_n \frac{|F_{t_n}(x)-F_{t_n}(z_{t_n,\eps})|}{\varepsilon} =1.
\]
Hence $\lip F(x) \geq  \limsup_{\eps\to0} \frac{|F(x)-F(z_\eps)|}{|xz_\eps|} = 1$. The second item follows from \cite[Proposition 5.19]{Gigli12}. {Finally, the third item follows by the fact that $F$ is a convex function, see \cite{CheegerGromollSoul}, and then there exist well-defined gradient curves of $F$ starting from any point, see \cite[Proposition 2.12]{PetruninSurvey}.}
\end{remark}

}

{
The next lemma follows from \cref{lem:ConoHalfline} together with a non-smooth version of \cite[Theorem 4.3]{Kasue88}. An alternative proof consists in arguing by contradiction, exploiting \cref{lem:ConoHalfline} and the monotonicity condition; details are left to the reader.
}

\begin{lemma}\label{lem:CompactLevels}
Let $(X,\dist)$ be a noncompact $N$-dimensional $\CBB(0)$ metric space with $N\geq 2$. Let $\gamma:[0,+\infty)\to X$ be a ray emanating from $o\in X$ such that $(X,\gamma(t_i))$ pGH-converges to $(\mathbb R\times K,(0,k_0))$ for some compact $K$, with $k_0\in K$. Suppose that $X$ has one end. Let $F$ be the Busemann function associated to $\gamma$.

Then $\{F\le T\}$ is compact for any $T \in \R$.
\end{lemma}

{
\begin{remark}\label{rem:CompletenessBusemann}
Using our convention of signs, it is known that $-\sup_{\gamma}F_\gamma$, where $\gamma$ are rays emanating from a fixed point, is proper on Riemannian manifolds under the nonnegative sectional curvature assumption \cite{CheegerGromollSoul}, and similarly for nonnegatively curved Alexandrov spaces. On manifolds with nonnegative Ricci curvature and linear volume growth every Busemann function is proper, and this comes from \cite[Theorem 19]{SormaniJDG98}, while in general it might be false due to the very recent example by J. Pan and G. Wei \cite{PanWei}.
\end{remark}
}

\begin{lemma}[Monotonicity of diameters]\label{lem:MonotonicityDiam}
Let $(X,\dist)$ be a noncompact $N$-dimensional $\CBB(0)$ metric space with $N\geq 2$. Let $\gamma:[0,+\infty)\to X$ be a ray emanating from $o\in X$. Let $F$ be the Busemann function associated to $\gamma$.

If there exists a sequence $t_j\to+\infty$ such that $\{F\le t_j\}$ is compact for any $j$, then
\[
\diam\{ F= r\} \le \diam\{F=R\},
\]
for any $r\le R$.
\end{lemma}

\begin{proof}
{Apply the third item in \cref{rem:BasicPropertiesBusemann}.}
\end{proof}

\begin{proposition}\label{prop:StimaFq}
    Let $N\ge 2$, and let $(X,\dist)$ be a noncompact $N$-dimensional $\CBB(0)$ space. Let $\gamma$ be a ray emanating from $o\in X$. Let $F$ be the Busemann function associated to the ray $\gamma$. Let $t_q\ge0$ be a nonnegative real number such that $\dist(q,\gamma)=|q\gamma(t_q)|$, for any $q \in X$. Then, for every $q\in X$ \begin{equation}\label{eqn:EstimateBusemannLower}
        F(q)\geq t_q.
    \end{equation}
    
    Let us further assume that there exists a constant $C>0$ for which
    \[
    \dist(p,\gamma)\leq C,
    \]
    for all $p\in X$. 
    Then, for every $\varepsilon>0$ there exists $R>0$ such that for every $q\in X$ with $|qo|\geq R$ we have 
    \begin{equation}\label{eqn:EstimateBusemann}
        F(q)\leq t_q+\varepsilon.
    \end{equation}
\end{proposition}

\begin{proof}
Fix $q\in X$, and let $t_q$ be as in the assumptions. Let $t>t_q$ be an arbitrary real number, and let us denote $p_q:=\gamma(t_q)$.\\
Let us first prove the lower bound. By the hinge condition on $(p_q)_q^{\gamma(t)}$ and the fact that $\angle (p_q)_q^{\gamma(t)}=\pi/2$ by \cref{lem:GeodCadePerpendicolare}, we get that 
\begin{equation}\label{eqn:StimaClaim1Zero}
    \begin{split}
    |q\gamma(t)|^2\leq |qp_q|^2+|p_q\gamma(t)|^2 &\leq (t-t_q)^2\left(\frac{|qp_q|^2}{(t-t_q)^2}+1\right), \quad \text{and then}\\
    |q\gamma(t)|&\leq (t-t_q)\left(1+\frac{|qp_q|^2}{2(t-t_q)^2}\right).
    \end{split}
\end{equation}
From the last inequality, together with the fact that $\dist(p,\gamma)\leq C$, we get that 
\[
F_t(q)=t-|q\gamma(t)| \geq t_q-\frac{C^2}{2(t-t_q)},
\]
and taking $t\to +\infty$ in the previous inequality we get $F(q)\geq t_q$, that is the sought lower bound. 
\smallskip

Let us now prove the upper bound. \\
\textbf{Claim 1}. For every $\delta>0$ there exists $R=R(\delta)>0$ such that the following holds. For every $q\in X$ with $|oq|\geq R$ there exists $T=T(\delta,q)>0$ such that for every $t\geq T$ and for every unit speed geodesics $\sigma,\beta$ connecting $q$ to $o$, and $q$ to $\gamma(t)$, respectively, we have 
\[
s^2+r^2+2sr(1-\delta)\leq |\sigma(s)\beta(r)|^2, \quad \forall t\geq T, \quad \forall r\in (0,|q\gamma(t)|), \quad \forall s\in (0,|oq|).
\]
\smallskip 
Let us prove the above Claim with some estimates that are similar to those in the proof of \cref{Lemma:KeyEstimate}. By the monotonicity condition, and the hinge condition, together with \cref{lem:GeodCadePerpendicolare}, we have 
\begin{equation}\label{eqn:StimaClaim1a}
\begin{split}
    \frac{s^2+r^2-|\sigma(s)\beta(r)|^2}{2sr}&\leq \frac{|oq|^2+|q\gamma(t)|^2-t^2}{2|oq||q\gamma(t)|}\leq \frac{t_q^2+2|qp_q|^2+(t-t_q)^2-t^2}{2|oq||q\gamma(t)|}\\
    &=-\frac{(t-t_q)t_q}{|oq||q\gamma(t)|}+\frac{|qp_q|^2}{|oq||q\gamma(t)|}.
\end{split}
\end{equation}
From the triangle inequality we have $t_q\geq |oq|-|qp_q|$, and thus 
\begin{equation}\label{eqn:StimaClaim1b}
\frac{t_q}{|oq|}\geq 1-\frac{|qp_q|}{|oq|}.
\end{equation}
Moreover, from the triangle inequality we also have 
\begin{equation}\label{eqn:StimaClaim1b1}
|q\gamma(t)|\geq (t-t_q)-|qp_q|= (t-t_q)\left(1-\frac{|qp_q|}{t-t_q}\right).
\end{equation}
Putting the second of \eqref{eqn:StimaClaim1Zero}, \eqref{eqn:StimaClaim1b}, and \eqref{eqn:StimaClaim1b1} in \eqref{eqn:StimaClaim1a} we finally get
\begin{equation}\label{eqn:StimaClaim1c}
\begin{split}
    \frac{s^2+r^2-|\sigma(s)\beta(r)|^2}{2sr} \leq -\frac{1}{1+\frac{|qp_q|^2}{2(t-t_q)^2}}\left(1-\frac{|qp_q|}{|oq|}\right)+\frac{|qp_q|^2}{|oq|(t-t_q)}\frac{1}{1-\frac{|qp_q|}{t-t_q}}.
\end{split}
\end{equation}
Now, inspecting the previous inequality, one immediately realizes that, since $|qp_q|=\dist(q,\gamma)\leq C$ by assumption, for every fixed $\delta>0$ one can find $R$ such that the statement of the Claim holds.
\smallskip 

\textbf{Claim 2}. For every $\varepsilon>0$ there exists $R=R(\varepsilon)>0$ such that the following holds. For every $q\in X$ with $|oq|\geq R$, and that does not belong to the ray $\gamma$, there exists $T=T(q,\delta)$ such that for every $t\geq T$ we have $|q\gamma(t)|\geq |oq|$, and
\[
\angle q_{p_q}^{\gamma(t)}\leq \pi/2+\varepsilon.
\]
Let us prove the above Claim with some estimates that are similar to those in the proof of \cref{Lemma:KeyEstimate}. First, by taking $T$ large enough with respect to $q$ we can always choose $T$ in such a way that for every $t\geq T$ we have $|q\gamma(t)|\geq |oq|$.
By exploiting the hinge condition and $|q\gamma(t)|\geq |oq|$ we get
\begin{equation}\label{eqn:StimaClaim2a}
    \frac{|oq|^2+|qp_q|^2-t_q^2}{2|oq||qp_q|}\leq \frac{|qp_q|^2+t_q^2+|qp_q|^2-t_q^2}{2|oq||qp_q|} = \frac{|qp_q|}{|oq|},
\end{equation}
\begin{equation}\label{eqn:StimaClaim2b}
    \frac{|qp_q|^2+|q\gamma(t)|^2-(t-t_q)^2}{2|qp_q||q\gamma(t)|}\leq \frac{|qp_q|^2+|qp_q|^2+(t-t_q)^2-(t-t_q)^2}{2|qp_q||q\gamma(t)|}=\frac{|qp_q|}{|q\gamma(t)|}\leq\frac{|qp_q|}{|oq|}.
\end{equation}
By using \eqref{eqn:StimaClaim2a} and \eqref{eqn:StimaClaim2b}, and noticing that $|qp_q|=\dist(q,\gamma)\leq C$, we can find $R$ big enough such that 
\begin{equation}\label{eqn:AngoliGrandi}
\begin{split}
\min\{\angle q_{p_q}^o, \angle q_{p_q}^{\gamma(t)}\} &\geq \min\left\{\arccos\left(\frac{|oq|^2+|qp_q|^2-t_q^2}{2|oq||qp_q|}\right),\arccos\left(\frac{|qp_q|^2+|q\gamma(t)|^2-(t-t_q)^2}{2|qp_q||q\gamma(t)|}\right) \right\} \\ 
&\geq \frac{\pi}{2}-\frac{\varepsilon}{2},
\end{split}
\end{equation}
for every $q\in X$ that is not on $\gamma$, with $|oq|\geq R$, and for every $t$ sufficiently big depending on $q$. 
Moreover, by using the Claim 1, see in particular the first inequality in \eqref{eqn:StimaClaim1a}, we have that, by taking a possibly larger $R$, we can also ensure 
\begin{equation}\label{eqn:AngoliGrandi2}
\angle q_o^{\gamma(t)}\geq \pi-\frac{\varepsilon}{2},
\end{equation}
for every $q\in X$ that is not on $\gamma$, with $|oq|\geq R$, and for every $t$ sufficiently big depending on $q$. Finally, by the four-point comparison, as in the end of \cref{Lemma:KeyEstimate}, we get that 
\[
\angle q_{p_q}^o + \angle q_{p_q}^{\gamma(t)} + \angle q_o^{\gamma(t)}\leq 2\pi.
\]
By comparing the previous inequality with \eqref{eqn:AngoliGrandi}, and \eqref{eqn:AngoliGrandi2} we finally get 
\[
\angle q_{p_q}^{\gamma(t)}\leq \frac{\pi}{2}+\varepsilon,
\]
which is the sought claim. 
\smallskip

Let us now finish the proof of the upper bound. We can assume without loss of generality that $q$ is not a point of the ray $\gamma$. Let us fix $\varepsilon>0$. Let us consider $\varepsilon'>0$ small enough such that the following implication holds 
\begin{equation}\label{eqn:CondAngolo}
0\leq \alpha\leq \frac{\pi}{2}+\varepsilon' \Rightarrow -C\cos\alpha\leq \varepsilon.
\end{equation}
Now, by using {Claim 2}, we can find $R=R(\varepsilon')>0$ such that for every $q\in X$ with $|oq|\geq R$ that does not belong to the ray $\gamma$, there exists $T=T(q,\delta)$ such that for every $t\geq T$ we have
\[
\angle q_{p_q}^{\gamma(t)}\leq \pi/2+\varepsilon'.
\]
Now, by the hinge condition, by \eqref{eqn:CondAngolo}, and the fact that $|qp_q|=\dist(q,\gamma)\leq C$ we have that for every $q\in X$ with $|oq|\geq R$ and $t\geq T$, 
\begin{equation}
    \begin{split}
    (t-t_q)^2&\leq |qp_q|^2+|q\gamma(t)|^2-2|qp_q||q\gamma(t)|\cos \angle q_{p_q}^{\gamma(t)} \\
    &\leq |q\gamma(t)|^2\left(1+\frac{|qp_q|^2}{|q\gamma(t)|^2}+\frac{2\varepsilon}{|q\gamma(t)|}\right), \\
    \end{split}
\end{equation}
and then, for $t$ sufficiently big,
\begin{equation}
    \begin{split}
    (t-t_q) &\leq |q\gamma(t)|\left(1+\frac{|qp_q|^2}{2|q\gamma(t)|^2}+\frac{\varepsilon}{|q\gamma(t)|}\right) \Rightarrow \\
    t-|q\gamma(t)|&\leq t_q + \frac{|qp_q|^2}{2|q\gamma(t)|}+\varepsilon.
    \end{split}
\end{equation}
Finally, taking $t\to +\infty$ in the previous inequality we get, by exploiting that $|qp_q|=\dist(q,\gamma)\leq C$, and $|q\gamma(t)|\to +\infty$, that $F(q)\leq t_q+\varepsilon$ for every $q\in X$ with $|oq|\geq R$; which is the sought upper bound.
\end{proof}

Let us now prove a Lemma that provides a bound on the distance of points from a ray in a $\CBB$ space whose asymptotic cone is a half-line.

{
\begin{lemma}\label{lemma:LimitiDistanzeDaRaggi}
Let $(X,\dist)$ be a noncompact $N$-dimensional $\CBB(0)$ metric space with $N\geq 2$. Let $\gamma:[0,+\infty)\to X$ be a ray emanating from $o\in X$ such that $(X,\gamma(t_i))$ pGH-converges to $(Y,y):=(\mathbb R\times K,(0,k_0))$ for some compact $K$, with $k_0\in K$. 

Let us assume $X$ has one end. Then 
\begin{equation}\label{eqn:LimitiPalle}
\limsup_{i\to+\infty} \max_{\partial B_{t_i}(o)} \dist(p,\gamma)<+\infty.
\end{equation}
Moreover, if $F$ is the Busemann function associated to $\gamma$, see \cref{def:Busemann}, we also have that 
\begin{equation}\label{eqn:LimitiBusemann}
\limsup_{i\to+\infty} \max_{\{F=t_i\}} \dist(p,\gamma)<+\infty.
\end{equation}
\end{lemma}

\begin{proof}
Since $X$ has one end, $t_i\to+\infty$, and without loss of generality we can assume that $t_i$ is increasing.\\
\textbf{Proof of \eqref{eqn:LimitiPalle}}.
By using \cref{lem:ConoHalfline} we get that the asymptotic cone of $X$ is a half-line. As a consequence we get that 
\begin{equation}\label{eqn:Massimioti}
    \lim_{i\to+\infty}\frac{\max_{\partial B_{t_i}(o)}\dist(p,\gamma)}{t_i} = 0.
\end{equation}
Indeed, if not, we can find, for $i$ sufficiently big, points $q_i\in \partial B_{t_i}(o)$ such that $\dist(q_i,\gamma(t_i))\geq \dist(q_i,\gamma) > t_i\delta$, for some $\delta>0$. Then, in the realization of the convergence of $(X,t_i^{-1}\dist,o)$ to the asymptotic cone, $q_i$ and $\gamma(t_i)$ are converging to two different points in the boundary of the unit ball centered at the tip of the asymptotic cone. Anyway, this is not possible since the asymptotic cone is a half-line.

Let us assume, by contradiction, that up to subsequences we have 
\begin{equation}\label{eqn:Assurdo1}
\lim_{i\to +\infty} \max_{\partial B_{t_i}(o)}\dist(p,\gamma)=+\infty.
\end{equation}
Let us denote by $M_i:=\max_{\partial B_{t_i}(o)}\dist(p,\gamma)$ and by $\dist_i\eqdef\dist/M_i$ the rescaled distance.
\smallskip

\textbf{Claim 1}. We aim now at showing that for every $\varepsilon>0$, and $r\geq 1$, we have that, for every $i\geq i_\varepsilon$,
\begin{equation}\label{eqn:LoVoglioLAssurdo}
B_r^{\dist_i}(\gamma(t_i)) \subset \{p:\dist_i(p,\gamma)<\varepsilon\}.
\end{equation}

In order to prove \eqref{eqn:LoVoglioLAssurdo} let us prove that for every $\varepsilon>0$, and for every $r\geq 1$, we have $i_\varepsilon$ such that 
\begin{equation}\label{eqn:ClaimBordi2}  
\dist(\sigma(t),\gamma)\leq \varepsilon M_i,
\end{equation}
for every $i\geq i_\varepsilon$, and for every segment $\sigma:[0,T]\to X$ such that $\sigma(0)=\gamma(t_i)$, and $T\leq rM_i$.

If not, up to passing to subsequences, we have $\sigma_i:[0,T_i]\to X$ such that $\sigma_i(0)=\gamma(t_i)$, $T_i\leq rM_i$, and for $s_i\in [0,T_i]$ we have $\dist(\sigma_i(s_i),\gamma)\geq \varepsilon M_i$. {Notice that, since $M_i\to +\infty$ both $s_i,T_i\to +\infty$}. Then, up to subsequence, $\sigma_i$ locally uniformly converge to a ray $\sigma:[0,+\infty)\to X$ in the limit $Y$ emanating from $y$. Moreover $\gamma_i:=\gamma(t+t_i)$ converge to a line $\ell$ such that $\ell(0)=y$. By \cref{prop:GeodeticheProdotto} we either have that $\sigma=\ell|_{[0,+\infty)}$ or $\sigma=\ell|_{(-\infty,0]}$.

Thus we have, where in the following we take the sign $+$ when $\sigma=\ell|_{[0,+\infty)}$ and the sign $-$ when $\sigma=\ell|_{(-\infty,0]}$,
\begin{equation}
\begin{split}
    \lim_{i\to +\infty} \lim_{t\to 0^+}\left(1-\frac{|\sigma_i(t)\gamma(t_i\pm t)|^2}{2t^2}\right)&= \lim_{i\to +\infty} \lim_{t\to 0^+}\frac{2t^2-|\sigma_i(t)\gamma(t_i\pm t)|^2}{2t^2} \\
    &=\lim_{i\to +\infty}\cos\angle \sigma_i'(0)(\gamma^{\pm})'(t_i) = 1,
\end{split}
\end{equation}
where we are using the convergence of angles, see \cref{rem:AngleContinuity}. Thus 
\[
\lim_{i\to +\infty}\lim_{t\to 0^+} \frac{|\sigma_i(t)\gamma(t_i\pm t)|}{t}=0.
\]
Notice that, for $i$ sufficiently big, $t_i-s_i\geq t_i-rM_i \geq t_i/2$, since $M_i/t_i\to 0$. Then, for $i$ sufficiently big we can estimate, by using the monotonicity condition in \cref{def:CBB} and the fact that $s_i\leq T_i\leq rM_i$ 
\begin{equation}
\varepsilon M_i\leq \dist(\sigma_i(s_i),\gamma)\leq |\sigma_i(s_i)\gamma(t_i\pm s_i)|\leq s_i\lim_{t\to 0^+}\frac{|\sigma_i(t)\gamma(t_i\pm t)|}{t}\leq s_i\frac{\varepsilon}{2r}\leq M_i\frac{\varepsilon}{2},
\end{equation}
which is a contradiction. Thus \eqref{eqn:ClaimBordi2} is established. Now, from \eqref{eqn:ClaimBordi2} we immediately get the sought claim \eqref{eqn:LoVoglioLAssurdo}.
\smallskip

Let us now exploit the {Claim 1}, and in particular \eqref{eqn:LoVoglioLAssurdo}, to find the sought contradiction. By the very definition of $\dist_i$ and $M_i$ we get that, for every $i$, there exists $p_i\in\partial B^{\dist_i}_{t_i/M_i}(o)$ such that $\dist_i(p_i,\gamma)=1$. Let us denote $\gamma_i(\cdot):=\gamma(M_i\cdot)$ for every $i$. Thus there exists $s_i\in [0,+\infty)$ such that $\dist_i(p_i,\gamma)=\dist_i(p_i,\gamma_i(s_i))=1$. We can estimate, by using the triangle inequality,
\begin{equation}\label{eqn:UnaStimaAT}
1=\dist_i(p_i,\gamma_i(s_i))\geq |\dist_i(p_i,o)-\dist_i(o,\gamma_i(s_i))|=\left|\frac{t_i}{M_i}-s_i\right|.
\end{equation}
Thus, 
\begin{equation}\label{eqn:UnaStima2}
\begin{split}
\dist_i(p_i,\gamma(t_i))&=\dist_i(p_i,\gamma_i(t_i/M_i))\leq \dist_i(p_i,\gamma_i(s_i)) + \dist_i(\gamma_i(s_i),\gamma_i(t_i/M_i))\\
&= 1 + M_i^{-1}\dist(\gamma(M_is_i),\gamma(t_i)) \leq 2,
\end{split}
\end{equation}
where in the last inequality we are using \eqref{eqn:UnaStimaAT}. Notice that, for $i$ large enough, the previous inequality gives the sought contradiction. Indeed from \eqref{eqn:LoVoglioLAssurdo} with $r=2$ and $\varepsilon=1/2$, and from \eqref{eqn:UnaStima2}, we deduce that for $i$ large enough we have $p_i\in B^{\dist_i}_2(\gamma(t_i))\subset \{p:\dist_i(p,\gamma)<1/2\}$, and thus $\dist_i(p_i,\gamma)<1/2$, while on the contrary we have $\dist_i(p_i,\gamma)=1$.
\medskip

\textbf{Proof of \eqref{eqn:LimitiBusemann}}. By \cref{lem:NoLinesOneEnd} we get that the asymptotic cone of $X$ is a half-line. We now claim that 
\begin{equation}\label{eqn:ClaimRaggioInfinito}
\lim_{R\to +\infty}\sup_{p\in X\setminus B_R(o)}\frac{\dist(p,\gamma)}{|op|}=0.    
\end{equation}
Indeed, if not, we have a sequence of diverging points $p_n$ such that $\dist(p_n,\gamma)\geq \delta|op_n|$ for all $n\in\mathbb N$, and for some $\delta>0$. Thus we have $|p_n\gamma(|op_n|)|\geq \delta|op_n|$, and rescaling the metric by $|op_n|^{-1}$ we get that the boundary of the ball of radius one in the asymptotic cone has at least two points, which is a contradiction because it is a half-line.

Now, since the level sets of $F$ are compact, see \cref{lem:CompactLevels}, we get that there exist $R_i\to+\infty$ such that $\{F=t_i\}\subset X\setminus \overline B_{R_i}$. Now by the triangle inequality, see in particular \eqref{eqn:StimaClaim1b}, and by \eqref{eqn:EstimateBusemannLower} we have that for every $q\in\{F=t_i\}$ the following inequality holds
\[
t_i=F(q)\geq t_q \geq |oq|\left(1-\frac{\dist(q,\gamma)}{|oq|}\right).
\]
Thus, by exploiting the previous inequality we get that for every $i\in\mathbb N$ and $q\in\{F=t_i\}$ we have 
\[
\frac{\dist(q,\gamma)}{t_i}\leq \frac{\dist(q,\gamma)}{|oq|}\frac{1}{1-\frac{\dist(q,\gamma)}{|oq|}}.
\]
Hence, by also using \eqref{eqn:ClaimRaggioInfinito}, and the fact that $\{F=t_i\}\subset X\setminus\overline{B}_{R_i}$ for a sequence of $R_i\to \infty$, from the previous inequality we conclude that 
\begin{equation}\label{eqn:Massimioti2}
    \lim_{i\to+\infty}\max_{\{F=t_i\}}\frac{\dist(q,\gamma)}{t_i}=0,
\end{equation}
which is the analog of \eqref{eqn:Massimioti}.

Let us assume, by contradiction, that up to subsequences we have 
\begin{equation}\label{eqn:Assurdo2}
\lim_{i\to +\infty} \max_{\{F=t_i\}}\dist(p,\gamma)=+\infty.
\end{equation}
Let us call $M_i':=\max_{\{F=t_i\}}\dist(p,\gamma)$. Observe that $M_i'<+\infty$ by \cref{lem:CompactLevels}. Denote $\dist_i'\eqdef \dist/M_i'$.

\smallskip

\textbf{Claim 2}. We aim at showing that
for every $\varepsilon>0$, and $r\geq 1$, we have that, for every $i\geq i_\varepsilon$,
\begin{equation}\label{eqn:LoVoglioLAssurdo2}
B_r^{\dist_i'}(\gamma(t_i)) \subset \{p:\dist_i'(p,\gamma)<\varepsilon\}.
\end{equation}

Notice that the proof of the previous claim can be done arguing verbatim as in the proof of \eqref{eqn:LoVoglioLAssurdo}, where we only exploited \eqref{eqn:Massimioti}, whose analog we proved in \eqref{eqn:Massimioti2}, and $(X,\gamma(t_i))\to \mathbb (R\times K,(0,k_0))$.
\smallskip

\textbf{Claim 3}. For every $i$ there exists $x_i\in\{F=t_i\}$ and $\tau_i\in[0,+\infty)$ such that $1=\dist_i'(x_i,\gamma)=\dist_i'(x_i,\gamma(\tau_i))$. We claim that 
\begin{equation}\label{eqn:BastaStiClaim}
    \lim_{i\to +\infty}\lim_{t\to+\infty} \frac{t-\tau_i-\dist(x_i,\gamma(t))}{M_i'}=0.
\end{equation}
In order to prove \eqref{eqn:BastaStiClaim}, let us introduce two symbols that will simplify the notation. We denote by $o_R(1)$ an arbitrary function of $p\in X$ such that for all $\varepsilon>0$ there exists $R_\varepsilon>0$ such that for every $q\in X\setminus B_{R_\varepsilon}(o)$ we have $|o_R(1)|\leq \varepsilon$. Moreover, we denote by $o_R^t(1)$ an arbitrary function  of $p\in X$, and $t\geq 0$ such that for all $\varepsilon>0$ there is an $R_\varepsilon>0$ such that for every $q\in X\setminus B_{R_{\varepsilon}}(o)$ we have $|o_R^t(1)|\leq \varepsilon$ for every $t\geq T_q$, where $T_q$ might depend on $q$. 

By the computations in \eqref{eqn:StimaClaim1a} and \eqref{eqn:StimaClaim1c}, together with the information in \eqref{eqn:ClaimRaggioInfinito}, we get 
\[
\frac{|oq|^2+|q\gamma(t)|^2-t^2}{2|oq||q\gamma(t)|} \leq -1 + o_R^t(1).
\]
Furthermore, by the computations in \eqref{eqn:StimaClaim2a} and \eqref{eqn:StimaClaim2b}, together with \eqref{eqn:ClaimRaggioInfinito}, we deduce that 
\[
\frac{|oq|^2+\dist(q,\gamma)^2-t_q^2}{2|oq|\dist(q,\gamma)}\leq o_R(1),
\]
and
\[
\frac{\dist(q,\gamma)^2+|q\gamma(t)|^2-(t-t_q)^2}{2\dist(q,\gamma)|q\gamma(t)|}\leq o_R^t(1),
\]
where $t_q\geq 0$ is any real number such that $\dist(q,\gamma)=\dist(q,\gamma(t_q))$. Using the previous three inequalities, together with the four-point comparison, as done, for example, in the end of the proof of Claim 2 of \cref{prop:StimaFq}, we conclude that 
\begin{equation}\label{eqn:BastaStiClaim2}
\frac{\dist(q,\gamma)^2+|q\gamma(t)|^2-(t-t_q)^2}{2\dist(q,\gamma)|q\gamma(t)|}= o_R^t(1).
\end{equation}
Thus, since $\dist(q,\gamma)^2/(\dist(q,\gamma)|q\gamma(t)|)=o_R^t(1)$ by \eqref{eqn:ClaimRaggioInfinito}, from \eqref{eqn:BastaStiClaim2} we get
\begin{equation}\label{eqn:BastaStiClaim3}
\frac{\left(|q\gamma(t)|-(t-t_q)\right)\left(|q\gamma(t)|+(t-t_q)\right)}{\dist(q,\gamma)|q\gamma(t)|}=\frac{|q\gamma(t)|^2-(t-t_q)^2}{\dist(q,\gamma)|q\gamma(t)|}= o_R^t(1).
\end{equation}
From the triangle inequality, see \eqref{eqn:StimaClaim1b1}, and the inequality in the second line of \eqref{eqn:StimaClaim1Zero}, we get
\[
1-\frac{\dist(q,\gamma)}{|t-t_q|}\leq \frac{|q\gamma(t)|}{t-t_q}\leq 1+\frac{\dist(q,\gamma)^2}{2(t-t_q)^2},
\]
and thus
\begin{equation}\label{eqn:BastaStiClaim4}
\frac{|q\gamma(t)|+(t-t_q)}{|q\gamma(t)|} = 1+o_R^t(1).
\end{equation}
From \eqref{eqn:BastaStiClaim4} and \eqref{eqn:BastaStiClaim3}, we finally get 
\begin{equation}\label{eqn:BastaStiClaim5}
\frac{|q\gamma(t)|-(t-t_q)}{\dist(q,\gamma)}=o_R^t(1).
\end{equation}
Now, \eqref{eqn:BastaStiClaim} is a direct consequence of \eqref{eqn:BastaStiClaim5} and how $\tau_i$ is defined. 
\smallskip 

Let us now exploit the previous claims to find the sought contradiction, and then conclude \eqref{eqn:LimitiBusemann}. By using the notation in the Claim 3 above we have
\[
\frac{t_i}{M_i'}=\frac{F(x_i)}{M_i'}=\lim_{t\to+\infty}\frac{t-\dist(x_i,\gamma(t))}{M_i'} = \lim_{t\to+\infty}\left(\frac{t-\tau_i-\dist(x_i,\gamma(t))}{M_i'}\right)+\frac{\tau_i}{M_i'},
\]
and thus, by exploiting \eqref{eqn:BastaStiClaim}, we conclude 
\[
\lim_{i\to+\infty}\left|\frac{t_i-\tau_i}{M_i'}\right|=0.
\]
Thus, for $i$ large enough, we have 
\[
\dist_i'(x_i,\gamma(t_i))\leq \dist_i'(x_i,\gamma(\tau_i))+\dist_i'(\gamma(\tau_i),\gamma(t_i))\leq 2.
\]
In particular, choosing $r=2,\varepsilon=1/2$ in \eqref{eqn:LoVoglioLAssurdo2}, for $i$ large enough we have $x_i\in B_2^{\dist_i'}(\gamma(t_i))\subset \{p\in X:\dist_i'(p,\gamma)<1/2\}$, which is a contradiction since $\dist'_i(x_i,\gamma)=1$.
\end{proof}
}

Putting together the previous results, we obtain the following theorem, which characterizes $\CBB(0)$ spaces having a pGH limit at infinity equal to a cylinder over a compact space as those spaces which are tubular neighborhoods of some line or ray.

\begin{theorem}\label{thm:raggiosselimiteainfinito}
Let $N\geq 2$, and let $(X,\dist)$ be a noncompact $N$-dimensional $\CBB(0)$ metric space. The following are equivalent.
\begin{enumerate}
    \item There exists a diverging sequence of points $p_i$ on $X$ such that $(X,p_i)$ converges in the pGH topology to $\mathbb R\times K$, where $K$ is compact. 
    
    \item Either $X$ is isometric to $\R\times K$, or $X$ has one end and there exists a ray $\gamma:[0,+\infty)\to X$ and a constant $C$ such that 
    \[
    \dist(p,\gamma)\leq C\qquad \text{for all $p\in X$}.
    \]
\end{enumerate}
\end{theorem}

{
\begin{proof}
Implication (2)$\Rightarrow$(1) is trivial. Implication (1)$\Rightarrow$(2) easily follows combining \cref{cor:RettaGenerataDaSuccessione}, \cref{lemma:LimitiDistanzeDaRaggi}, and \cref{lem:MonotonicityDiam}. In fact, if $X$ has one end, the previous results yield a ray $\gamma:[0,+\infty)\to X$ emanating from a point $o$ and a uniform upper bound on the diameter of any level set of the Busemann function associated to $\gamma$. Since any level set $\{F=t\}$ contains $\gamma(t)$, for $t\ge0$, the claim follows.
\end{proof}
}

\begin{theorem}\label{thm:ConvergenzaBusemann}
Let $N\geq 2$, and let $(X,\dist)$ be a noncompact $N$-dimensional $\CBB(0)$ metric space such that $\mathcal{H}^N(B_1(x))\geq v_0>0$ for some $v_0>0$ for every $x\in X$.
Let $\gamma:[0,+\infty)\to X$ be a ray such that, for some constant $C>0$,
\[
\dist(p,\gamma)\leq C, \qquad \text{for all $p\in X$}.
\]
Let $t_j\to +\infty$ and assume that $(X,\gamma(t_j))$ pGH-converges to $(\mathbb R\times K,(0,k_0))$ for some $(N-1)$-dimensional compact $\CBB(0)$ space $K$.

Let us denote by $f:\mathbb R\times K\to \R$ the coordinate function $f(t,k)=t$. Let $F_j:=F-t_j$. Then the following hold.
\begin{itemize}
    \item For every $j\in\mathbb N$ and every $T\geq 0$, the set $\{F_j\leq T\}$ is compact.
    
    \item $F_j$ converges locally uniformly to $f$ along the pointed sequence of converging metric spaces. In particular, for any $s, s_1,s_2 \in \R$ with $s_1<s_2$, we have that
    \begin{equation}\label{eq:ConvergenceLevelSetBusemann}
    \begin{split}
    \{F_j<s\}\to \{f<s\},\qquad &\text{in $L^1_{\rm loc}$-strong}, \\
    \{s_1<F_j<s_2\}\to \{s_1<f<s_2\},\qquad &\text{in $L^1$-strong}.
    \end{split}
    \end{equation}
    
    \item For any $s \in \R$, the set $\{F<  s\}$ has finite perimeter and
    \begin{equation}\label{eq:MonotoniaPerimetriBusemann}
    \mathrm{Per}(\{F< s_1\})\leq \mathrm{Per}(\{F< s_2\}),\qquad \text{$\forall\, s_1, s_2 \in \R$ with $s_1<s_2$,}
    \end{equation}
    \begin{equation}\label{eq:LimitePerimetriBusemann}
        \lim_j \Per(\{F_j<t\}) =\lim_{s\to +\infty}\mathrm{Per}(\{F< s\})=\mathcal{H}^{N-1}(K) 
        \qquad \forall\, t \in \R.
    \end{equation}
    Moreover, if $\Per (\{F< \bar s\})=\mathcal{H}^{N-1}(K)$ for some $\bar s \in \R$, then there exists $s_0\ge \bar s$ and a connected component $\{F=\bar s\}_\alpha$ of $(\{F=\bar s\},\dist_{\rm int})$ such that $(\{F=s_0\}, \dist|_{\{F=s_0\}})$ and $(\{F=\bar s\}_\alpha, \dist_{\rm int})$ are isometric to $K$, and $(X\setminus \{F<s_0\}, \dist|_{X\setminus \{F<s_0\}})$ is isometric to $(\{F=s_0\}\times [0,+\infty),\dist|_{\{F=s_0\}} \otimes \dist_{\rm eu})$, which is isometric to $K\times[0,+\infty)$.\\
    In particular $X\setminus \{F<s_0\}$ splits isometrically as a connected cylindrical end\footnote{We observe that, arguing as in the proof of \cref{prop:ProfiloCostiffSPlitCilindro}, one can prove that $\{F<s_0\}$ is an isoperimetric set, up to choosing a greater $s_0$. In particular the end of $X$ is a cylinder over the boundary of an isoperimetric set.} and
    \[
    \haus^{N-1}(K) = \haus^{N-1} (\{F=s_0\}) =\mathcal{H}^{N-1}_{\rm int}(\{F=\bar s\}_\alpha),
    \]
    where $\mathcal{H}^{N-1}_{\rm int}$ is the $(N-1)$-dimensional Hausdorff measure on $(\{F=\bar s\}_\alpha, \dist_{\rm int})$ and $\dist_{\rm int}$ is the intrinsic distance on $\{F= \bar s\}$.
    
    \item {If $q_j$ is an arbitrary sequence of diverging points on $X$, then up to subsequences $(X,q_j)$ pGH-converges to $\mathbb R\times K$.}
\end{itemize}
\end{theorem}

\begin{proof}
Observe that the assumption implies that $X$ has one end.
The first item follows from \cref{lem:CompactLevels}. Since $\mathcal{H}^N(B_1(x))\geq v_0>0$ for some $v_0>0$ for every $x\in X$, the convergence of $(X,\gamma(t_j))$ to $\R\times K$ holds in pmGH-sense by \cite{DePhilippisGigli18}.
The fact that $F_j$ locally uniformly converges to $f$ immediately follows from \cref{prop:StimaFq}. The functions $F_j$ are uniformly bounded in balls $B_R(\gamma(t_j))$ and are $1$-Lipschitz by \cref{rem:BasicPropertiesBusemann}. It is then possible to derive local convergence in $L^1$-strong of the sublevels $\{F_j<s\}$ to $\{f<s\}$ for any $s \in \R$, see for example the proof of \cite[Proposition 3.2]{AmbrosioHondaTewodrose}, and \eqref{eq:ConvergenceLevelSetBusemann} follows.

Recalling \cref{rem:BasicPropertiesBusemann}, by \cite[Proposition 6.1]{BPSGaussGreen} we get that there exists a set of full measure $S\subset \R$ such that $\{F<s\}$ is a set of finite perimeter for any $s \in S$. Moreover, in the notation of \cite{BPSGaussGreen}, we have $(\nabla F \cdot \nu_{\{F<s\}})_{\rm int}=-1$ at $\Per(\{F<s\},\cdot)$-a.e. point, for any $s \in S$. Hence we can apply that Gauss--Green formula in \cite[Theorem 5.2]{BPSGaussGreen}, invoking \cite[Theorem 4.11]{BPSGaussGreen}, to get
\begin{equation*}
    0\le \int_{\{s_1<F<s_2\}} \de \boldsymbol{\Delta}F = \Per(\{F<s_2\}) - \Per(\{F<s_1\}),
\end{equation*}
for any $s_1,s_2 \in S$ with $s_1<s_2$. By lower semicontinuity of perimeters, choosing a sequence $(s_1')_n\searrow s_1<s_2$ with $(s_1')_n \in S$, we deduce that 
\begin{equation}\label{eq:zzGG1}
    \Per(\{F<s_1\})\le \Per(\{F<s_2\}),
\end{equation}
for any $s_1\in \R, s_2 \in S$ with $s_1<s_2$, and that $\{F<  s\}$ has finite perimeter for any $s \in \R$. Therefore we can apply the Gauss--Green formula in \cite[Theorem 5.2]{BPSGaussGreen} again integrating by parts $\boldsymbol{\Delta}F$ on $\{s_1<F<s_2\}$ for any $s_1<s_2$ with $s_1 \in S, s_2 \in \R$ to get
\begin{equation}\label{eq:zzGG2}
\begin{split}
    0&\le \int_{\{s_1<F<s_2\}} \de \boldsymbol{\Delta}F = \int - (\nabla F \cdot \nu_{\{F<s_2\}})_{\rm int} \de \Per(\{F<s_2\},\cdot)  - \Per(\{F<s_1\}) \\
    &\le \|(\nabla F \cdot \nu_{\{F<s_2\}})_{\rm int}\|_{L^\infty(\Per(\{F<s_2\},\cdot)} \Per(\{F<s_2\}) - \Per(\{F<s_1\}) \\
    &
    \le \Per(\{F<s_2\}) - \Per(\{F<s_1\}).
\end{split}
\end{equation}
Therefore we finally get that for any $s_1,s_2 \in \R$ with $s_1<s_2$, taking any $s'\in(s_1,s_2)\cap S$, there holds
\[
\Per(\{F<s_1\}) \overset{\eqref{eq:zzGG1}}{\le} \Per(\{F<s'\}) \overset{\eqref{eq:zzGG2}}{\le} \Per(\{F<s_2\}),
\]
proving \eqref{eq:MonotoniaPerimetriBusemann}.

Now by \eqref{eq:ConvergenceLevelSetBusemann}, and exploiting \cref{rem:SemicontPerimeterConverging}, we have that $\liminf_j \Per(\{F_j<s\}) \ge \Per(\{f<s\}) = \haus^{N-1}(K)$, for any $s \in \R$. We want to show that equality actually holds. Indeed, if by contradiction, up to subsequence, there holds $\lim_j \Per(\{F_j<s\}) > \haus^{N-1}(K)$ for some $s \in \R$, then for $s_2>s$ we estimate
\[
\begin{split}
    \haus^N&(\{s<f<s_2\}) = (s_2-s)\haus^{N-1}(K) < \lim_j (s_2-s) \Per(\{F_j<s\})\\& \overset{\eqref{eq:MonotoniaPerimetriBusemann}}{\le} \limsup_j \int_s^{s_2} \Per(\{ F_j <t\}) \de t 
    = \limsup_j \haus^N(\{ \{s<F_j<s_2\}) \overset{\eqref{eq:ConvergenceLevelSetBusemann}}{=} \haus^N(\{s<f<s_2\}).
\end{split}
\]
This proves that $\lim_j \Per(\{F_j<s\}) = \haus^{N-1}(K)$ for any $s \in \R$. Moreover we find
\[
\Per(\{F<s\}) \overset{\eqref{eq:MonotoniaPerimetriBusemann}}{\le} \Per(\{F<s+t_j\}) = \Per(\{F_j<s\}) \qquad\forall\,j,
\]
hence $\limsup_{s\to\infty} \Per(\{F<s\}) \le \limsup_{s\to\infty} \lim_j \Per(\{F_j<s\}) = \haus^{N-1}(K)$. On the other hand
\[
\liminf_{s\to+\infty} \Per(\{F<s\}) \overset{\eqref{eq:MonotoniaPerimetriBusemann}}{\ge} \Per(\{F< t_j\}) = \Per(\{F_j<0\}) \qquad\forall\,j,
\]
hence $\liminf_{s\to+\infty} \Per(\{F<s\}) \ge \lim_j \Per(\{F_j<0\}) = \haus^{N-1}(K)$, completing the proof of \eqref{eq:LimitePerimetriBusemann}.

Suppose now that $\Per (\{F< \bar s\})=\mathcal{H}^{N-1}(K)$ for some $\bar s \in \R$, hence the same holds for any $s\ge \bar s$ by \eqref{eq:MonotoniaPerimetriBusemann}. Recalling \eqref{eq:zzGG2}, this implies that $\boldsymbol{\Delta} F =0$ on $\Omega\eqdef X\setminus \{F\le \bar s\}$. Hence we can apply \cref{thm:SplittingKetterer} (notice that $X\setminus \{F \le \bar s\}$ has only one unbounded connected component as $X$ has one end). In the notation of \cref{thm:SplittingKetterer} and letting $\Omega\eqdef X \setminus \{F\le \bar s\}$, we get the existence of a connected component $\{F= \bar s\}_\alpha$ of $(\{F= \bar s\}, \dist_{\rm int})$ such that the unbounded connected component $\widetilde\Omega_\alpha$ of the completion of $X\setminus \{F \le \bar s\}$ with respect to the induced intrinsic distance is isometric to $(\{F= \bar s\}_\alpha \times [0,+\infty), \dist_{\rm int} \otimes \dist_{\rm eu})$. {For $s_0\ge \bar s$ sufficiently large,} we have that $\dist|_{\{F\ge s_0-1\}}=\widetilde \dist$ by \cref{thm:SplittingKetterer}. Hence $(\{F= \bar s\}_\alpha, \dist_{\rm int})$ is isometric to $K$ and $\haus^{N-1}(K) = \mathcal{H}^{N-1}_{\rm int}(\{F=\bar s\}_\alpha)$.

Let $\Phi:\widetilde\Omega_\alpha \to \{F=\bar s\}_\alpha \times[0,+\infty)$ be the isometry given by \cref{thm:SplittingKetterer}. Computing the Busemann function at points of $\{F\ge s_0\} \cap \widetilde\Omega_\alpha$ with respect to $\dist$ or $\widetilde\dist$ is equivalent, thus $\Phi$ sends level sets $\{F=\rho\}\cap \widetilde\Omega_\alpha$ into level sets of the Busemann function of the ray $\Phi \circ \gamma|_{[\bar s,\infty)}$ in $(\{F=\bar s\}_\alpha \times[0,+\infty), \dist_{\rm int} \otimes \dist_{\rm eu})$, for $\rho\ge s_0$. Hence the product structure of $\widetilde\Omega_\alpha$ implies that $\Phi(\{F=\rho\}\cap \widetilde\Omega_\alpha)=\{F=\bar s\}_\alpha \times \{\phi(\rho)\}$ for some $\phi:[s_0,+\infty) \to [0,+\infty)$, for $\rho \ge s_0$.
Since $\Phi$ is an isometry, we get that $(\{F=s_0\}\cap \widetilde\Omega_\alpha, \dist|_{\{F=s_0\}\cap \widetilde\Omega_\alpha})=(\{F=s_0\}\cap \widetilde\Omega_\alpha, \widetilde\dist|_{\{F=s_0\}\cap \widetilde\Omega_\alpha})$ is isometric to $(\{F=\bar s\}_\alpha, \dist_{\rm int})$.\\
Since $\dist=\widetilde\dist$ in a neighborhood of $\{F=s_0\}\cap \widetilde\Omega_\alpha$, we deduce that $\haus^{N-1}(\{F=s_0\}\cap \widetilde\Omega_\alpha)<+\infty$ and
\begin{equation}\label{eq:zuw}
\haus^{N-1}(\{F=s_0\}\cap \widetilde\Omega_\alpha) = \mathcal{H}^{N-1}_{\rm int}(\{F=\bar s\}_\alpha) = \haus^{N-1}(K). 
\end{equation}
Since $\Omega$ has only one unbounded connected component, up to choosing a greater $s_0$, we have that $\{F=s_0\} \subset \widetilde\Omega_\alpha$. Thus $\{F=s_0\}\cap \widetilde\Omega_\alpha= \{F=s_0\}$, which completes the proof of the rigidity in the third item.

{
It remains to prove the last item. Up to passing to subsequence, by assumption there are $\tau_j\nearrow+\infty$ such that $\dist(q_j,\gamma)=|q_j\gamma(\tau_j)| \le C$, hence $(X,\gamma(\tau_j))$ pGH-converges to some product $\R\times K'$, where $K'$ is compact. Up to further subsequences, we can assume that $\tau_{j-1}<t_j < \tau_j$ for any $j$. By \cref{rem:BasicPropertiesBusemann} we find a sequence of $1$-Lipschitz surjective functions $\Psi_j: \{ F = \tau_j\} \to \{F=t_j\}$. By uniform convergence of sequences of functions $F-\tau_j$ and $F- t_j$ proved above, an Arzelà--Ascoli argument implies the existence of a $1$-Lipschitz surjective map $\Psi:K'\to K$. Analogously, there exists a $1$-Lipschitz surjective map $\Phi:K\to K'$. Hence $K$ and $K'$ are isometric, see, e.g., \cite[Lemma 4.5]{AntBruFogPoz}.}
\end{proof}
{
In the case of manifolds with nonnegative Ricci curvature, or for arbitrary $\RCD(0,N)$ spaces, the picture depicted in \cref{cor:Equivalenze} is more complicated. In fact, there exist Riemannian manifolds with nonnegative Ricci curvature and linear volume growth such that the diameters of the level sets of a Busemann function are unbounded, \cite[Example 26]{SormaniJDG98}. On the other hand, as mentioned in the Introduction, it was recently proved in \cite{XingyuZhu} that on smooth manifolds $(M,g)$ with $\Ric\ge0$ and $\inf_{x\in M} \haus^N(B_1(x))>0$, conditions (1) and (2) in the next \cref{cor:Equivalenze} are equivalent.
}

\begin{corollary}\label{cor:Equivalenze}
Let $N\geq 2$, and let $(X,\dist)$ be a noncompact $N$-dimensional $\CBB(0)$ metric space such that $\mathcal{H}^N(B_1(x))\geq v_0>0$ for some $v_0>0$ for every $x\in X$. Let $o \in X$. Then the following are equivalent. 
\begin{enumerate}
    \item $\limsup_{r\to+\infty} \haus^N(B_r(o)) /r <+\infty$;
    
    \item for any diverging sequence of points $p_i\in X$ such that $(X,p_i)$ pGH-converges to a limit $(Y,y)$, it holds that $Y$ is isometric to a product $\R\times K$, where $K$ is compact;
    
    \item there exists a diverging sequence of points $p_i\in X$ such that $(X,p_i)$ pGH-converges to a limit $\R\times K$, where $K$ is compact.

    \item There exists a constant $C>0$ such that $I(V)\leq C$ for every $V>0$.
\end{enumerate}
{If any of the previous items holds, then there exists a unique compact $\CBB(0)$ space $K$ such that any pGH limit of $X$ along a sequence of diverging points is isometric to $\R\times K$.}
\end{corollary}

\begin{proof}
We observe first that on any noncompact $\RCD(0,N)$ space $(Z,\dist,\meas)$, for any $z \in Z$ the estimate $\meas(B_r(z))\ge C(N,\meas(B_1(z))) r$ holds for any $r\ge 1$. Indeed, this follows by performing the proof of \cite[Theorem 4.1]{SchoenYauLectures} verbatim exploiting the Laplacian comparison from \cite[Corollary 5.15]{Gigli12}. We skip this proof, whose strategy is however analogous to the first implication below. In particular, an $\RCD(0,N)$ space with finite measure is compact.

(1)$\Rightarrow$(2). We know from \cite[Corollary 5.15]{Gigli12} that the squared distance function $\dist_o^2$ from $o$ has measure Laplacian $\boldsymbol{\Delta} \dist_o^2 \le 2N \haus^N$. Let $x \in X$, $r\in(0,|ox|)$, and define $r_1\eqdef |ox|-r$ and $r_2\eqdef |ox|+r$. Integrating the function
\[
f(y) \eqdef \begin{cases}
1 & |oy|\le r_1 , \\
\frac{1}{r_2-r_1}(r_2-\dist_o(y)) & r_1<|oy|<r_2,\\
0 & r_2\le |oy|,
\end{cases}
\]
with respect to $\boldsymbol{\Delta} \dist_o^2$, by chain rule (see \cite[Eq. (4.16)]{Gigli12}), we get
\[
\begin{split}
    2N\haus^N(B_{r_2}(o)) &\ge \int f \de \boldsymbol{\Delta} \dist_o^2 = - \int \nabla f \cdot (2 \dist_o \nabla\dist_o) \de \haus^N
    =2 \int_{B_{r_2}(o)\setminus B_{r_1}(o)}  \frac{\dist_o}{r_2-r_1} |\nabla \dist_o|^2 \de \haus^N \\
    &\ge \frac{2}{r_2-r_1} r_1 \haus^N(B_{r_2}(o)\setminus B_{r_1}(o)) \ge \frac{2}{r_2-r_1} r_1 \haus^N(B_r(x)).
\end{split}
\]
Letting $c>\limsup_{r\to+\infty} \haus^N(B_r(o)) /r$, we find that
\begin{equation}\label{eq:zxw}
    \haus^N(B_r(x)) \le 4cN \frac{|ox|+r}{|ox|-r} \, r,
\end{equation}
for any $x\in X$ such that $|ox|\ge R_0$ and any $r \in (0,|ox|)$, for some $R_0$ large enough.\\
Let now $p_i$ and $Y$ as in (2). By \cref{lem:EveryLimitAtInfinitySplits} we know that $Y$ is isometric to some product $\R\times K$. Plugging $x=p_i$ in \eqref{eq:zxw}, by volume convergence we get that $\haus^N(B_r(y)) \le 4cN r$ for any $r>0$, hence $K$ must have finite total measure. Since $K$ is $\CBB(0)$, it must be compact.

(2)$\Rightarrow$(3). Trivial.

(3)$\Rightarrow$(1). By \cref{thm:raggiosselimiteainfinito}, either $X$ is isometric to $\R\times K$, or $X$ has one end and there exists a ray $\gamma:[0,+\infty)\to X$ from $o$ and a constant $C$ such that $\dist(p,\gamma)\leq C$ for all $p\in X$. Without loss of generality we can assume the second alternative. Hence there exist $t_j\to+\infty$ such that $(X,\gamma(t_j))$ pGH-converges to $\R\times K$, and \cref{thm:ConvergenzaBusemann} applies.\\
Let $r>0$, $x \in B_r(o)$ and $t_x\ge0$ such that $\dist(x,\gamma(t_x))=\dist(x,\gamma)$. Hence $t-\dist(x,\gamma(t)) \le t-(t-t_x-\dist(x,\gamma)) \le t_x + C$, for any $t>0$ sufficiently big. Hence the Busemann function $F$ of the ray $\gamma$ satisfies $F(x) \le t_x+C$. Moreover, $|x\gamma(r)|< 2r$ and $|x\gamma(t)|\ge t-r - |x\gamma(r)|> t-3r$ for any $t\ge r$, hence $t_x< \max\{r,C+3r\}= 3r+C$. Therefore $F(x)< 3r+2C$, and thus we proved the containment
\[
B_r(o) \subset \{F < 3r+2C\},
\]
for any $r>0$. Since sublevel sets of $F$ are compact, $F$ has a global minimum $m\in\R$. Hence, by \cref{thm:ConvergenzaBusemann}, coarea formula, \eqref{eq:MonotoniaPerimetriBusemann}, \eqref{eq:LimitePerimetriBusemann}, we conclude that
\[
\begin{split}
    \haus^N(B_r(o)) &\le \haus^N(\{F < 3r+2C\}) = \int_m^{3r+2C} \Per( \{F<s\}) \de s \\
    &\le \haus^{N-1}(K) (3r+2C-m).
\end{split}
\]

(3)$\Rightarrow$(4). If $X$ has one end, it directly follows from \cref{thm:raggiosselimiteainfinito} and the third item of \cref{thm:ConvergenzaBusemann} (see also the third item of \cref{thm:AsintoticaProfilo} for comparison). If $X$ has at least two ends, then $X$ is isometric to $\mathbb R\times K$ and the statement follows from \cref{lem:IsoperimetricLargeVolumesProduct}.

(4)$\Rightarrow$(1) or (3). Let $V(r):=\inf_{x\in X}\mathcal{H}^N(B_r(x))$. For $v>0$, let us define 
\[
\phi(v):=\inf\{r>0:V(r)\geq v \}.
\]
Arguing verbatim as in \cite[Theorem 7.4]{LeonardiRitore}, since the proof there only uses the Bishop--Gromov inequality and the relative isoperimetric inequality, both available in our setting, we have that there exists a constant $C:=C(N)$ such that the following holds. For every $v>0$ we have
\begin{equation}\label{eqn:IsoperimetricInequalityBella}
I(v)\geq C\frac{v}{\phi(2v)}.
\end{equation}
We now claim that $\limsup_{r\to+\infty} V(r)/r < +\infty$. Indeed, if it is not the case, for every $M$ we find $r_0>1$ such that for every $r\geq r_0$ we have $V(r)\geq Mr$. Then, for every $v\geq Mr_0$, we have that $V(v/M)\geq v$, and thus $\varphi(v)\leq v/M$. Hence, for $v\geq Mr_0$, we have 
\[
I(v)\geq C\frac{v}{\varphi(2v)}\geq C\frac{v}{2v/M}= \frac{CM}{2}.
\]
This implies that $\lim_{v\to +\infty} I(v)=+\infty$, which is against to our assumption on $I$. Therefore
$$
\limsup_{r\to+\infty} V(r)/r =:C_0 < +\infty.
$$

We now prove that if $\limsup_{r\to+\infty}\mathcal{H}^N(B_r(o))/r=+\infty$, then there is a diverging sequence of points $q_i\in X$ such that $(X,q_i)$ pGH-converges to a limit $\mathbb R\times K$, where $K$ is compact. This would conclude the proof of the implication (4)$\Rightarrow$(1) or (3).\\
So, let us now assume $\limsup_{r\to+\infty}\mathcal{H}^N(B_r(o))/r=+\infty$. In particular we have $r_i\to +\infty$ such that $\mathcal{H}^N(B_{r_i}(o))/r_i\to +\infty$ as $i\to +\infty$. Let us take $p_i\in X$ such that 
\[
V(r_i)\geq \mathcal{H}^N(B_{r_i}(p_i))-i^{-1}.
\]
We now claim that
\begin{equation}\label{eq:zqw}
    \frac{|op_i|}{r_i}\to +\infty
\end{equation}
If not, up to subsequences, there is a constant $\bar C$ such that $|op_i|\leq \bar Cr_i$. Thus, by using Bishop--Gromov inequality, we have, up to subsequences 
\[
\mathcal{H}^N(B_{r_i}(o))\leq \mathcal{H}^N(B_{(\bar C+1)r_i}(p_i))\leq (\bar C+1)^N\mathcal{H}^N(B_{r_i}(p_i)) \leq (\bar C+1)^N\left(2C_0r_i+i^{-1}\right),
\]
which is in contradiction with the fact that $\mathcal{H}^N(B_{r_i}(o))/r_i\to +\infty$. 

Let us now prove that there exists a constant $\tilde C$ for which the following holds. For every $i$ large enough, we have $\rho_i\in (r_i/4,3r_i/4)$, and $q_i\in\partial B_{\rho_i}(p_i)$ such that $\mathcal{H}^N(B_h(q_i))\leq \tilde Ch$ for every $h\in (0,r_i/4)$.\\
Let us notice that up to subsequences we have 
\[
\int_0^{r_i}\mathrm{Per}(B_\rho(p_i))\de\rho = \mathcal{H}^N(B_{r_i}(p_i))\leq 2C_0r_i,
\]
and then, by using the Hardy--Littlewood maximal inequality, we have that there exists a constant $C>0$ such that for every $t>0$, and every $i$ we have
\[
t\, \mathscr{L}^1\left(\left\{\rho\in(0,r_i):\sup_{h>0}\fint_{\rho-h}^{\rho+h}\mathrm{Per}(B_\xi(p_i))\de \xi>t\right\}\right)\leq Cr_i.
\]
Thus choosing $t=4C$ in the previous inequality we conclude that for $i$ large enough there exists $\rho_i\in (r_i/4,3r_i/4)$ such that 
\[
\mathcal{H}^N\left(B_{\rho_i+h}(p_i)\setminus B_{\rho_i-h}(p_i)\right) = \int_{\rho_i-h}^{\rho_i+h}\mathrm{Per}(B_s(p_i))\de s\leq 8Ch, \quad \forall h\in (0,r_i/4).
\]
Then taking an arbitrary $q_i\in\partial B_{\rho_i}(p_i)$ we have, by the triangle inequality, that for every $i$ and every $h\in (0,r_i/4)$, we have
\[
\mathcal{H}^N\left(B_h(q_i)\right) \leq \mathcal{H}^N\left(B_{\rho_i+h}(p_i)\setminus B_{\rho_i-h}(p_i)\right)\leq 8Ch,
\]
and thus the claim is proved with $\tilde C:=8C$.\\
Let now $(\tilde X,\tilde\dist,\tilde q)$ be a pGH limit of $(X,\dist,q_i)$ up to subsequences. By the results of \cite{DePhilippisGigli18}, since $\mathcal{H}^N(B_1(x))\geq v_0$ for every $x\in X$, we get that $(\tilde X,\tilde\dist,\mathcal{H}^N,\tilde q)$ is the pmGH limit of $(X,\dist,\mathcal{H}^N,q_i)$, and then, by the volume convergence and the fact that $\rho_i\to +\infty$, we conclude that $\mathcal{H}^N(B_h(\tilde q))\leq \tilde Ch$ for every $h>0$.\\
Since
\[
|oq_i|\geq |op_i|-|p_iq_i|\geq |op_i|-r_i \xrightarrow[i]{\eqref{eq:zqw}} +\infty,
\]
we get that $\tilde X$ splits as $\mathbb R\times K$, see \cref{lem:EveryLimitAtInfinitySplits}. From the latter volume estimates on $B_h(\tilde q)$, we get that $(K,\dist_K)$ is a $\CBB(0)$ space with finite volume. Hence, it is compact, and we have proved (3).

{Finally, the claimed uniqueness on the pGH limits at infinity follows from the last item of \cref{thm:ConvergenzaBusemann}.}
\end{proof}

\section{Volume non-collapsedness of $2$-dimensional $\mathrm{CBB}(0)$ spaces}\label{sec:Dimension2}

{As pointed to us by A. Lytchak, the forthcoming \cref{lem:LowBoundDiametri} is likely to hold in arbitrary dimensions by using the gradient flow of the Busemann function of a ray. Since we will not need it at this level of generality, we provide a direct and elementary proof of it in dimension $2$.}
\begin{lemma}\label{lem:LowBoundDiametri}
Let $(X,\dist)$ be a noncompact $2$-dimensional $\CBB(0)$ metric space. Then, for every $o\in X$, there exists $C>0$ such that 
\[
\mathrm{diam}(\partial B_r(o))\geq C, \qquad \text{for all $r\geq 1$}.
\]
\end{lemma}

\begin{proof}
Assume by contradiction that there is an increasing sequence $r_j\nearrow+\infty$ such that $\mathrm{diam}(\partial B_{r_j}(o)) \to 0$.
Hence, rescaling distance by $1/r_j$ and passing to the limit we get that the asymptotic cone to $X$ is (isometric to) a halfline $([0,+\infty),\dist_{\rm eu})$.\\
Fix a ray $\gamma$ from $o$. We observe that
\begin{equation}\label{eq:zz2}
    \lim_{\rho\to+\infty} \frac1\rho \max_{x \in \partial B_\rho(o)} \dist(x,\gamma)  =0.
\end{equation}
Indeed, if $x_\rho$ maximizes in \eqref{eq:zz2}, as $(X,\dist/\rho)\to ([0,+\infty),\dist_{\rm eu})$, we see that $\dist(x_\rho,\gamma)/\rho\le \dist(x_\rho,\gamma(\rho))/\rho \to 0$ as $\rho\to\infty$, as $x_\rho$ and $\gamma(\rho)$ converge to the same point in the asymptotic cone.

We want to show that $\partial B_\rho(o)=\{\gamma(\rho)\}$ for some $\rho$ large enough. Suppose instead that for any $\rho>0$ there is a point $x_\rho \in \partial B_\rho(o)\setminus\{\gamma(\rho)\}$. Let $\tau>\rho$ and let $\alpha=[\gamma(\tau) x_\rho]$. Choosing $\tau\ge 3\rho$ there is $T \in (0, \tau]$ such that $\dist(x_\rho,\gamma)=|x_\rho \gamma(\tau-T)|$. Observe that $\rho = |ox_\rho| \le \dist(x_\rho,\gamma) + (\tau-T)$, and thus
\begin{equation}\label{eq:zz3}
   |x_\rho\gamma(\tau)| \le |x_\rho \gamma(\tau-T)| + T
   \le 2\dist(x_\rho,\gamma) + \tau -\rho 
   \overset{\eqref{eq:zz2}}{=} o(\rho) -\rho + \tau < \tau,
\end{equation}
for any $\rho$ large enough.

As $t\mapsto f(t)\eqdef \dist(\alpha(t),\gamma(\tau-t))/t$ is nonincreasing, $L(\alpha)=|x_\rho \gamma(\tau)|<\tau$ by \eqref{eq:zz3}, then $f(|x_\rho \gamma(\tau)|)>0$ is well-defined, and thus $\alpha \cap \gamma = \{\gamma(\tau)\}$, for any $\rho$ large enough.

For large $\rho$, we also have that $\alpha \subset X\setminus B_{\frac\rho2}(o)$. Indeed, if $\alpha(\bar t) \in B_{\frac\rho2}(o)$, letting $t_1\in(0,\bar t)$ such that $\alpha(t_1)\in \partial B_\rho(o)$, we would get
\[
|x_\rho \gamma(\tau)| = L(\alpha) = |\gamma(\tau)\alpha(t_1)| + L (\alpha|_{[t_1,|x_\rho \gamma(\tau)|]})\ge \tau - \rho + |\alpha(t_1)\alpha(\bar t)| +|\alpha(\bar t) x_\rho| \ge \tau- \rho + 2\frac\rho2  = \tau, 
\]
contradicting \eqref{eq:zz3}. Since for $\rho$ large enough we have $\alpha \subset X\setminus B_{\frac\rho2}(o)$, and \eqref{eq:zz2} implies that $\dist(x,\gamma) < |ox|/10$ for any $x \in X \setminus B_{\frac{\hat \rho}{2}}(o)$, then $\dist(\alpha(t),\gamma)<|o\alpha(t)|/10$ for any $\rho>\hat\rho$ large enough. Thus the distance $\dist(\alpha(t),\gamma)$ is never achieved at $o$, for $t \in [0,|x_\rho\gamma(\tau)|]$.

Finally we want to show that, for given $\rho$, choosing $\tau$ large enough guarantees that $\dist(\alpha(t),\gamma)$ is never achieved at $\gamma(\tau)$, for $t \in (0,|x_\rho\gamma(\tau)|]$. Indeed, let $r_j>\rho$ for suitable $j$ and define $\tau_j=r_j+1$, let $w_j\in \partial B_{r_j}(o) \cap \alpha$. For large $j$, we have that $\angle \gamma(\tau_j)^{w_j}_o < \tfrac\pi2$. Indeed, any pGH limit of $(X,\gamma(\tau_j))$ is of the form $(\R\times X', (0,x'))$ and the absurd assumption implies that $w_j\to (-1,x')$ along $(X,\gamma(\tau_j)) \to (\R\times X', (0,x'))$. Convergence of angles from \cref{rem:AngleContinuity} implies that $\angle \gamma(\tau_j)^{w_j}_o \to 0$ as $j\to+\infty$. So, for $j\ge j_\rho$ large enough, choosing any $\tau=\tau_j$ implies that $\angle \gamma(\tau)^{w_j}_o < \tfrac\pi2$ and thus \cref{lem:GeodCadePerpendicolare} implies that $\dist(\alpha(t),\gamma)$ cannot be achieved at $\gamma(\tau)$, for $t \in (0,|x_\rho\gamma(\tau)|]$.

All in all, for $\rho$ large and $\tau>\rho$ of the form $\tau=r_j+1$ for $j\ge j_\rho$, we have proved that $\alpha \cap \gamma = \{\gamma(\tau)\}$, $\dist(\alpha(t),\gamma)$ is not achieved at $o$ or $\gamma(\tau)$, for $t \in (0,|x_\rho\gamma(\tau)|]$. Hence we are in position to apply \cref{prop:ConcavityDistance} with $\sigma=\alpha$, getting that
\[
(0,|x_\rho \gamma(\tau)|) \ni t \mapsto \dist(\alpha(t),\gamma)
\]
is concave.

Now select $j_\rho \le j_0 < j_1 < j$ such that $r_{j_0}>\rho$ and take $\tau=r_j+1$. Let $t_0>t_1>s>0$ be the first time that $\alpha$ intersects $\partial B_{r_{j_0}}(o)$, $\partial B_{r_{j_1}}(o)$, $\partial B_{r_{j}}(o)$, respectively. By concavity we have that
\begin{equation*}
    \frac{\dist(x_\rho,\gamma) - \dist(\alpha(t_0),\gamma)}{|x_\rho\gamma(\tau)|- t_0} \le \frac{\dist(\alpha(t_1),\gamma) - \dist(\alpha(s),\gamma) }{t_1-s}.
\end{equation*}
By the absurd assumption, $\diam (\partial B_{r_k}(o) )\le D$ for any $k$. Moreover $|x_\rho\gamma(\tau)|- t_0 = |x_\rho \alpha(t_0)| \le |ox_\rho| + |o \alpha(t_0)| = \rho + r_{j_0}$. Hence
\begin{equation*}
\begin{split}
    \dist(x_\rho,\gamma) &\le
    \dist(\alpha(t_0),\gamma)
    +  (\rho + r_{j_0}) \, \frac{\dist(\alpha(t_1),\gamma) + \dist(\alpha(s),\gamma) }{t_1-s}
    \\
    &\le 
    \diam ( \partial B_{r_{j_0}}(o) ) +  (\rho + r_{j_0}) \, \frac{\diam ( \partial B_{r_{j_1}}(o) ) + \diam ( \partial B_{r_{j}}(o) ) }{t_1-s}
    \\
    &\le 
    \diam ( \partial B_{r_{j_0}}(o) ) +  (\rho + r_{j_0}) \, \frac{2D}{t_1-s}.
\end{split}
\end{equation*}
We first pass to the limit $j\to+\infty$, keeping other indices fixed. Hence $t_1\xrightarrow[]{j} +\infty$ and $s\xrightarrow[]{j} 1$. Then
\[
\dist(x_\rho,\gamma) \le \diam ( \partial B_{r_{j_0}}(o) ).
\]
The previous inequality holds for any $j_0\ge j_\rho$ with $r_{j_0}>\rho$. Letting now $j_0\to+\infty$, we deduce that $\dist(x_\rho,\gamma)=0$, i.e., $\partial B_{\rho}(o)$ is a point and it coincides with $\{\gamma(\rho)\}$.

Therefore we proved that $\partial B_{\rho}(o)=\{\gamma(\rho)\}$ for any $\rho$ large enough. Hence $X\setminus B_{\bar \rho}(o) = \gamma \setminus B_{\bar \rho}(o)$ for some $\bar\rho>0$. By monotonicity condition, any two geodesics $[\gamma(\bar\rho+1)x]$, $[\gamma(\bar\rho+1)y]$, for any $x,y \in  B_{\bar \rho}(o)$, must coincide on the common interval of definition, as they coincide on the time interval $[0,1]$. Hence $X$ coincides with the ray $\gamma$, implying that $X$ has Hausdorff dimension equal to $1$, that contradicts the hypothesis.
\end{proof}

We mention that an alternative proof of the forthcoming \cref{thm:Noncollapsed} is sketched in \cite[Remark 2.13]{PositiveScalarRicci}, that appeared on arXiv while the present paper was being completed.

\begin{theorem}[Volume non-collapsedness]\label{thm:Noncollapsed}
Let $(X,\dist)$ be a $2$-dimensional $\CBB(0)$ metric space. Then there exists $v_0>0$ such that 
\[
\mathcal{H}^2(B_1(x))\geq v_0,
\]
for all $x\in X$.
\end{theorem}

\begin{proof}
Without loss of generality, $X$ is noncompact, it has no lines and ${\rm AVR}(X)=0$. Suppose by contradiction that there is a sequence $x_i \in X$ such that $\haus^2(B_1(x_i))\searrow 0$. Hence $x_i$ diverges along $X$. Up to subsequence, we can assume that $(X,x_i)$ pGH-converges to a limit $(Y,y_0)$. By \cite[Theorem 1.2]{DePhilippisGigli18}, ${\rm dim}(Y) \le 1$. By \cref{lem:EveryLimitAtInfinitySplits} we deduce that $(Y,y_0)=(\R,y_0)$ is (isometric to) a line. By \cref{cor:RettaGenerataDaSuccessione} we find $o \in X$, a ray $\gamma$ from $o$, and times $t_i\nearrow+\infty$ such that $(X,\gamma(t_i))\to (\R,y_1)$ in pGH sense. By \cref{lem:LowBoundDiametri}, there is $C>0$ such that $\diam(\partial B_r(o))\ge C$ for any $r>1$. Hence
\begin{equation}\label{eq:zaz1}
    \exists\,\eps>0 \st \forall\,i \,\exists\, p_i \in \partial B_{t_i}(o) \st
    \dist(p_i,\gamma)\ge 2\eps .
\end{equation}
Indeed, if $|p\gamma(s_p)|\eqdef \dist(p,\gamma)\le \delta_i\to 0$ for any $p \in \partial B_{t_i}(o) $ along some (non-relabeled) subsequence, then $|p\gamma(t_i)|\le |p\gamma(s_p)| + |s_p-t_i|\le 2|p\gamma(s_p)| \le 2\delta_i$ for any $p \in \partial B_{t_i}(o) $. Hence $\diam (\partial B_{t_i}(o)) \le 4\delta_i\to0$, that is impossible.

Denote $|p_i\gamma(s_i)|\eqdef \dist(p_i,\gamma)$. By \cref{lem:NoLinesOneEnd}, $X$ has one end, and then $|s_i-t_i|\le |p_i\gamma(s_i)|=\dist(p_i,\gamma)\le C_1$ by \eqref{eqn:LimitiPalle}. Hence we can write that $(X, \gamma(s_i)) \to (\R,0)$ in pGH sense. By \eqref{eq:zaz1}, the geodesic $\alpha_i\eqdef [\gamma(s_i)p_i]$ is defined on $[0,2\eps]$ for any $i$. Let $w_i\eqdef \alpha_i(\eps)$. Hence $w_i$ converges, up to subsequences, to $w \in \{\eps,-\eps\}\subset \R$, and \cref{lem:GeodCadePerpendicolare} and \cref{rem:AngleContinuity} imply that
\[
\frac\pi2 = \angle \gamma(s_i)_{\gamma(s_i+\eps)}^{p_i} \xrightarrow[i]{} \angle 0_{\eps}^w \in \{0,\pi\},
\]
reaching a contradiction.
\end{proof}

\section{Isoperimetric problem}\label{sec:IsoperimetricProblem}

First, let us recall some notions about the boundary of non-smooth spaces. In the setting of the {non-collapsed} $\RCD(K,N)$ spaces ({i.e., the reference measure is $\haus^N$}), there are two notions of boundary. The first introduced by De Philippis--Gigli in \cite{DePhilippisGigli18}, according to which 
\[
\partial^* X=\overline{\mathcal{S}^{N-1}\setminus \mathcal{S}^{N-2}},
\]
where we refer to \cite{DePhilippisGigli18} for the definition of the singular strata $\mathcal{S}^k$. Notice that for every $x\in\mathcal{S}^{N-1}\setminus \mathcal{S}^{N-2}$ the tangent cone is unique and isomorphic to $\mathbb R_+^N$, see \cite[item (iii), Theorem 1.4]{BrueNaberSemola20}. In particular, we also have, by the volume convergence \cite{DePhilippisGigli18}, that the density of $\mathcal{H}^N$ at every point in $\mathcal{S}^{N-1}\setminus \mathcal{S}^{N-2}$ is $1/2$, and $\leq 1/2$ for every $x\in\partial^* X$ by the lower semicontinuity of the density. 

In \cite{MondinoKapovitch} a new notion of boundary $\partial X$ for the {non-collapsed} $\RCD(K,N)$ spaces ({i.e., the reference measure is $\haus^N$}) has been introduced \cite[Definition 4.2]{MondinoKapovitch}. As noticed in that paper, this notion agrees with the notion of boundary in the Alexandrov setting, which is studied, e.g., in \cite{Perelman91}. 

In the setting of arbitrary {non-collapsed} $\RCD(K,N)$ spaces ({i.e., the reference measure is $\haus^N$}) it is known that the notion of having nonempty boundary is independent of the definition of boundary we choose, see \cite[Theorem 6.6]{BrueNaberSemola20}. In the setting of Alexandrov spaces with curvature bounded from below it has been proved in \cite[Item (iii) of Theorem 7.4]{BrueNaberSemola20} that, whenever $\partial^* X\neq 0$, we have that $\partial^* X = \partial X$, so that the two notions of boundary fully coincide on Alexandrov spaces with curvature bounded from below with $\partial^* X\neq \emptyset$.

Finally, by \cite[Theorem 1.6]{BrueNaberSemola20}, we know that the pmGH limit $X$ of {non-collapsed} $\RCD(K,N)$ spaces $X_i$ ({i.e., the reference measure is $\haus^N$}) {with uniform bounds on the volume of unit balls}, and without boundary $\partial^* X_i$, has no boundary $\partial^* X$. We remark that some of the results obtained in the $\RCD$ setting in \cite{BrueNaberSemola20} were already known in the Alexandrov setting. E.g., the non-collapsing of the boundary along pGH limits, and the volume convergence of the boundary in \cite[Theorem 9.2, Remark 9.13]{KapovitchPerelmanStability}, and \cite[Theorem 1.3]{Fujioka}, respectively.
\medskip

Before giving the proof of \cref{thm:Existence2D}, let us record here the following fact, which implicitly plays a fundamental role in the proof of \cref{thm:Existence2D}. 
{The next corollary immediately follows from \cref{lem:EveryLimitAtInfinitySplits}, the last item of \cref{thm:ConvergenzaBusemann}, the boundary stability result from \cite[Theorem 1.6]{BrueNaberSemola20} and the 
non-collapsing of boundaries from \cite{KapovitchPerelmanStability}.}

\begin{corollary}
    Let $(X,\dist)$ be a $2$-dimensional $\CBB(0)$ metric space. Then the following holds. 
    \begin{enumerate}
        \item If $\partial X$ is empty or bounded, then either all the pGH limits at infinity are isometric to the Euclidean plane $\mathbb R^2$, or all the pGH limits at infinity are isometric to $\mathbb R\times \mathbb S^1(\rho)$ for a fixed $\rho>0$;
        \item If $\partial X$ is unbounded, then either all the pGH limits at infinity are isometric to $\mathbb R\times [0,\ell]$ for a fixed $\ell$, or all the pGH limits at infinity are either isometric to $\mathbb R^2$, or $\mathbb R\times [0,+\infty)$.
    \end{enumerate}
\end{corollary}

{We now state a Lemma analogous to the classical Morgan--Johnson result \cite{MorganJohnson00}, and that will be employed to possibly pull-back the mass lost at infinity when exploiting the direct method of the Calculus of Variations to solve the isoperimetric problem. The proof easily follows as in \cite[Theorem 3.5]{MorganJohnson00} by Bishop--Gromov monotonicity, taking into account the rigidity from \cite[Theorem 8.2]{BrueNaberSemola20}.}

\begin{lemma}\label{lem:MJHalfspaces}
Let $N\geq 2$, and let $(X,\dist,\mathcal{H}^N)$ be a noncompact {non-collapsed} $\RCD(0,N)$ space ({i.e., the reference measure is $\haus^N$}) metric space. Let us denote $\mathbb R^N_+\eqdef [0,+\infty)\times\mathbb R^{N-1}$ the $N$-dimensional halfspace.

Let us assume $\partial^* X\neq\emptyset$. Then, for every $x\in\partial^* X$, and for any $V>0$, there exists a ball $B^X\subset X$ centered at $x$ such that $\haus^N(B^X)=V$ and $\Per (B^X) \le \Per (B_{\mathbb R^N_+}(V))$, where $B_{\mathbb R^N_+}(V)$ is the ball of center $0$ and volume $V$ in $\mathbb R^N_+$. Moreover if equality holds in the previous inequality we have that every ball $\tilde B$ centered at $x$ compactly contained in $B^X$ is isometric to a ball of the same radius centered at a boundary point of $\mathbb R_+^N$.
\end{lemma}

\begin{remark}
    Notice that \cref{lem:MJHalfspaces} gives the following result. For every noncompact {non-collapsed} $\RCD(0,N)$ space ({i.e., the reference measure is $\haus^N$}) $X$ with nonempty boundary $\partial^*X\neq\emptyset$, we have 
    \[
    I_X(V)\leq I_{\mathbb R^N_+}(V), \quad \forall V>0.
    \]
    Moreover, if equality holds for a volume $V_0$, we have that there exists a radius $\rho_0$ such that for every $x\in\partial^*X$ the ball $\overline{B}_{\rho_0}^X(x)$ is isometric to $\overline{B}_{\rho_0}^{\mathbb R^N_+}(0)$, and it is an isoperimetric set. This recovers in a more general setting \cite[Corollary 6.17]{LeonardiRitore}.
\end{remark}

\begin{lemma}\label{lem:MJCylindric}
Let $(X,\dist)$ be a noncompact $N$-dimensional $\CBB(0)$ metric space with $N\geq 2$ such that $\haus^N(B_1(x))\ge v_0$ for some $v_0>0$ for any $x$. Let $\gamma:[0,+\infty)\to X$ be a ray such that, for some constant $C>0$,
\[
\dist(p,\gamma)\leq C, \qquad \text{for all $p\in X$},
\]
and let $F$ be the Busemann function associated to $\gamma$. {Let $\mathfrak{V}\eqdef \mathcal{H}^{N-1}(K)$ where $K$ is the compact fiber of any pGH limit at infinity of $X$\footnote{The latter is well-defined by \cref{thm:ConvergenzaBusemann}.}.}

Then, for every $V>0$, there exists $s\in\mathbb R$ such that 
\[
\mathcal{H}^N(\{F<s\})=V, \qquad \mathrm{Per}(\{F<s\})\leq \mathfrak{V}.
\]
Moreover if for some $\bar s\in\mathbb R$ we have
\[
\mathrm{Per}(\{F< \bar s\})= \mathfrak{V},
\]
then there exists $s_0\ge \bar s$ such that $(\{F=s_0\}, \dist|_{\{F=s_0\}})$ is $\CBB(0)$ and $(X\setminus \{F< s_0\}, \dist|_{X\setminus \{F< s_0\}})$ is isometric to $(\{F= s_0\} \times [0,+\infty), \dist|_{\{F=s_0\}} \otimes \dist_{\rm eu})$.
\end{lemma}

\begin{proof}
The assumption implies that $X$ has one end.
From \cref{lem:CompactLevels} we have that for every $s\in\mathbb R$, $\{F\leq s\}$ is compact. Since $X$ has infinite volume, and since the function $\rho\to \mathcal{H}^N(\{F<\rho\})$ is continuous, we get that for every $V>0$ there exists $s$ such that $\mathcal{H}^N(\{F<s\})=V$. 
Now the fact that 
\[
\mathrm{Per}(\{F<s\})\leq \mathfrak{V},
\]
is a direct consequence of the third item of \cref{thm:ConvergenzaBusemann}, as well as the rigidity part.
\end{proof}

We are now ready to exploit our previous results to prove that on $2$-dimensional $\CBB(0)$ metric spaces the isoperimetric problem has always solution.

\begin{theorem}\label{thm:Existence2D}
Let $(X,\dist)$ be a noncompact $2$-dimensional $\CBB(0)$ metric space. Then for every $V>0$ there exists an isoperimetric region of volume $V$.
\end{theorem}

\begin{proof}
First, notice that as a consequence of \cref{thm:Noncollapsed} we get that for some $v_0>0$ we have $\mathcal{H}^2(B_1(x))\geq v_0$ for every $x\in X$. Let $V>0$. If we take a minimizing (for the perimeter) sequence of bounded sets of finite perimeter of volume $V$ and we apply \cref{lem:IsoperimetricAtFiniteOrInfinite} we get that either all the mass stays in $X$ or it escapes to precisely one pGH limit at infinity. 

If all the mass stays in $X$, then we end up with an isoperimetric set of volume $V$, and we are done. Otherwise, there exists a pGH limit at infinity $X_\infty$ of $X$ and a set of finite perimeter $E\subset X_\infty$, which is isoperimetric for its own volume in $X_\infty$, and such that $\mathcal{H}^N(E)=V$, and $I_X(V)=\mathrm{Per}(E)$. By \cref{lem:EveryLimitAtInfinitySplits} we get that $X_\infty$ splits a line, and thus $X_\infty$ is isometric to $\mathbb R\times X'$, where $X'$ is a one-dimensional $\CBB(0)$ metric space. Thus, due to the characterization of $1$-dimensional $\CBB$ metric spaces, see \cite[Theorem 15.18]{AKP}, $X_\infty$ is isometric to either $\mathbb R^2$, $\mathbb R\times[0,+\infty)$, $\mathbb R\times\mathbb S^1(\rho)$ for some $\rho>0$, or $\mathbb R\times [0,\ell]$ for some $\ell>0$. Thus we have four cases:
\begin{enumerate}
    \item $X_\infty\cong\mathbb R^2$. Since $E$ is an isoperimetric set, we deduce it is a ball of volume $V$. By applying the non-smooth version of Morgan--Johnson's result \cite[Theorem 3.5]{MorganJohnson00}, see \cite[Proposition 3.3]{AntonelliPasqualettoPozzetta21}, we get that for any ball $B\subset X$ of volume $V$ we have 
    \[
    I_X(V)\leq \mathrm{Per}(B)\leq \mathrm{Per}(E) = I_X(V).
    \]
    Thus, every ball of volume $V$ is an isoperimetric set and we are done\footnote{Actually, if this is the case, by the classical rigidity of Morgan--Johnson's result \cite[Theorem 3.5]{MorganJohnson00}, see also \cite[Corollary 1.7]{DePhilippisGigli18}, we can even deduce that all the balls of volume $V'<V$ in $X$ are isometric to flat Euclidean balls.}.
    \item $X_\infty\cong \mathbb R\times[0,+\infty)$. Since $E$ is an isoperimetric set, we deduce it is, up to isometry of $X_\infty$, a ball in $\mathbb R\times[0,+\infty)$ centered at 0, see \cite[Pages 5-6]{RosIsoperimetric}. Then by \cite[Theorem 1.6]{BrueNaberSemola20} we have $\partial X\neq\emptyset$, and we can argue as in (1) exploiting \cref{lem:MJHalfspaces}.
    \item $X_\infty\cong \mathbb R\times\mathbb S^1(\rho)$. Since $E$ is an isoperimetric set, we deduce it is, up to isometry of $X_\infty$, a geodesic ball of $\mathbb R\times \mathbb S^1$ that is isometric to a Euclidean flat ball, or a cyilinder $[a,b]\times \mathbb S^1(\rho)$, see \cite[Pages 5-6]{RosIsoperimetric} and \cite{PedrosaRitore}. Then, in the first case we argue verbatim as in item (1) above and we are done. In the second case, either $X\cong \mathbb R\times\mathbb S^1(\rho)$, and we are clearly done, or $X$ has one end. In the latter case, by using \cref{thm:raggiosselimiteainfinito} and \cref{lem:MJCylindric}, arguing as in item (1), we get a contradiction.
    \item $X_\infty\cong \mathbb R\times [0,\ell]$. Since $E$ is an isoperimetric set, we deduce it is either a geodesic ball centered on the boundary, or a slab $[a,b]\times [0,\ell]$, see \cite[Pages 5-6]{RosIsoperimetric}. Then in the first case we argue verbatim as in item (2) above, while in the second case we argue verbatim as in item (3) above, substituting $\mathbb S^1(\rho)$ with $[0,\ell]$.
\end{enumerate}
\end{proof}

A direct outcome of the \cref{thm:Existence2D}, which to the authors' knowledge was unknown, is the following.
\begin{corollary}
    Let $X$ be either an unbounded convex body in $\mathbb R^2$, or the boundary of an unbounded convex body in $\mathbb R^3$ (endowed with intrinsic distance). Then $X$ admits isoperimetric regions for every volume $V>0$.
\end{corollary}

From the \cref{thm:Existence2D} and the classification of low-dimensional $\RCD$ spaces \cite{KitabeppuLakzian}, we get existence of isoperimetric sets on \emph{any} $\RCD(0,2)$ space, possibly collapsed. {The proof of the next corollary easily follows by \cref{thm:Existence2D} and by the fact that, in case $(X,\dist,\mathfrak{m})$ is a collapsed $\RCD(0,2)$ space, then its essential dimension is either $0$ or $1$ \cite{DePhilippisGigli18, BSConstancy, BrenaGigliHondaZhu}, and the case of essential dimension equal to $1$ can be treated directly exploiting \cite[Theorem A.2]{CavallettiMilmanCD}.}

\begin{corollary}\label{cor:RCDK2}
    Let $(X,\dist,\mathfrak{m})$ be an $\RCD(0,2)$ metric measure space. Then for every $V>0$ there exists an isoperimetric set of volume $V$ in $X$.
\end{corollary}

\begin{remark}
    Notice that $K=0$ in \cref{cor:RCDK2} cannot be relaxed. Indeed, there are examples of $2$-dimensional surfaces with uniformly bounded but not non-negative sectional curvature such that no isoperimetric regions exist for every volume $V$, see, e.g., \cite[Example 3.5]{AFP21}, \cite{Rit01NonExistence}.
\end{remark}

On $2$-dimensional $\CBB(0)$ spaces without boundary we can deduce that isoperimetric regions either have connected boundary or are isometric to truncated cylinders $(0,b)\times \S^1(\rho)$.

\begin{theorem}\label{thm:IsopBoundaryConnected}
Let $(X,\dist)$ be a noncompact $2$-dimensional $\CBB(0)$ metric space with $\partial X=\emptyset$. Let $E$ be an isoperimetric set. Then one of the following alternatives holds.
\begin{enumerate}
\item $\partial E$ is connected.
\item There exist $b,\rho>0$ such that $X$ is isometric to a right cylinder $\mathbb R\times \mathbb S^1(\rho)$ and $E=(0,b)\times\mathbb S^1(\rho)$, up to translation along the factor $\mathbb R$.
\end{enumerate}
In particular, in case $X$ is not isometric to a Euclidean $2$-dimensional right cylinder, the boundary of every isoperimetric set is connected. In the case (1) above, $\partial E$ is parametrized by a Lipschitz curve homeomorphic to $\mathbb S^1$.
\end{theorem}

\begin{proof}
Notice that there exists $v_0>0$ such that $\mathcal{H}^2(B_1(x))\geq v_0>0$ for every $x\in X$, see \cref{thm:Noncollapsed}. Let $E$ be the isoperimetric set in the assumption. Let $V:=\mathcal{H}^N(E)$, and $P:=\mathrm{Per}(E)$. Recall $E$ is bounded and open, up to choice of representative. Let $c$ be a mean curvature barrier for $E$, \cref{def:MeanCurvatureBarrier}.
Assume $\partial E$ is not connected in $X$. By \cref{thm:Regularity} we know that $\partial E$ is a topological manifold, hence $\partial E$ has finitely many connected components, which are compact by boundedness. Take $\Omega_1,\Omega_2$ two connected components in $\partial E$ and let $P_i:=\mathcal{H}^{N-1}(\Omega_i)$ for $i=1,2$.
It holds that $P_i>0$ for $i=1,2$, as any component of $\partial E$ is parametrized by a nonconstant Lipschitz closed curve by \cref{thm:Regularity}. Notice that since $\partial E$ has finitely many compact connected components $\Omega_i$, we have that $\dist(\Omega_i,\Omega_j)\ge d_0>0 $ for any $i\neq j$.\\
Hence there exists $\varepsilon>0$ small enough such that if we define
\[
\Omega_1^t:=\{x\in X\setminus E:\dist(x,\Omega_1)\leq t\},
\]
\[
\Omega_2^t:=\{x\in \overline E:\dist(x,\Omega_2)\leq t\},
\]
for every $t\in (0,\varepsilon)$, then 
$\dist(x,\Omega_1)=\dist(x,\overline E)$ on $\Omega_1^t$, and $\dist(x,\Omega_2)=\dist(x,X\setminus E)$ on $\Omega_2^t$. Let us define
\[
f(t):=\mathrm{Per}(\Omega_1^t,X\setminus\overline E), \quad g(t):=\mathrm{Per}(\Omega_2^t,E).
\]
By integrating \eqref{eq:sharpLap0} we get the bounds
\[
f(t)\leq P_1(1+c\cdot t), \quad g(t)\leq P_2(1-c\cdot t),
\]
for every $t\in (0,\varepsilon)$, compare with \eqref{eq:extareabd}, and \eqref{eq:intareabound}. Recall that by \cref{prop:C>0} we have $c\geq 0$. Our aim is now to show that $c=0$. So, let us derive a contradiction if $c>0$ under our assumptions. First, notice that there exists $\eta<\varepsilon$ small enough such that, by continuity, for every $\alpha\in (0,\eta)$ we can find $t_1,t_2$ such that $\mathcal{H}^N(\Omega_1^{t_1})=\mathcal{H}^N(\Omega_2^{t_2})=\alpha$. 

\textbf{Claim}. If $c>0$, for every $\alpha\in (0,\eta)$, and $\eta<\varepsilon$ small enough, we have $P_1t_1<P_2t_2$.

In order to prove the claim let us first prove that for $\xi>0$ small enough we have 
\[
g\left(\frac{P_1}{P_2}\xi\right)<\frac{P_2}{P_1}f(\xi).
\]
To show this, let us introduce some notation. Let $E_1^t:=\Omega_1^t\cup \overline{E}$, $E_2^t:=(X\setminus E)\cup \Omega_2^t$. Let $\beta_1(t):=\mathcal{H}^N(E_1^t)$, and $f_i(t):=\mathrm{Per}(E_i^t)$, for $i=1,2$. Notice that by coarea formula $\beta_1'(0)=P_1$. Moreover $\beta_1(0)=V$. 
From the very definition of the isoperimetric profile we have, for $t$ small enough,
\[
I(\beta_1(t))\leq f_1(t)=f(t)+P-P_1.
\]
Let us now prove the sought inequality on $g,f$ above. If it is not true we then have $\xi_i\to 0^+$ such that 
\[
g\left(\frac{P_1}{P_2}\xi_i\right)\geq\frac{P_2}{P_1}f(\xi_i).
\]
Thus, by exploiting the inequality on $g$, and the one on $f_1$ above we get
\begin{equation}
\begin{split}
P_2\left(1-\frac{P_1}{P_2}c\xi_i\right)&\geq g\left(\frac{P_1}{P_2}\xi_i\right)\geq 
\frac{P_2}{P_1}f(\xi_i)
\geq \frac{P_2}{P_1}\left(I(\beta_1(\xi_i))-P\right)+P_2 \\
&\Rightarrow 
\frac{P_2}{P_1}\frac{I(\beta_1(\xi_i))-I(\beta_1(0))}{\xi_i}+P_1c\leq 0.
\end{split}
\end{equation}
Taking $i\to +\infty$ in the previous inequality, and since $I$ has right-derivative, see \cref{cor:IsoperimetricProfileRCD0N}, we conclude that 
\[
\frac{P_2}{P_1}I'_+(V)\cdot P_1+P_1c \leq 0,
\]
which is a contradiction with $c>0$ and the fact that $I'_+(v)\geq 0$ as a consequence of the fact that $I$ is non-decreasing, see \cref{cor:IsoperimetricProfileRCD0N}. Thus the sought claim with $f,g$ is proved. Let us now show the \textbf{Claim}. By contradiction $P_1t_1\geq P_2t_2$. Then, by using coarea formula and the previous bounds we get that if $t_1,t_2$ are small enough, which we can always ensure by taking $\eta$ small enough,
\begin{equation}
\begin{split}
\alpha=\mathcal{H}^N(\Omega_2^{t_2})&=\int_0^{t_2}g(s)\de s = \int_0^{\frac{P_2}{P_1}t_2}g\left(\frac{P_1}{P_2}\vartheta\right)\frac{P_1}{P_2}\de \vartheta\leq  \\ 
&\leq \int_0^{t_1}g\left(\frac{P_1}{P_2}\vartheta\right)\frac{P_1}{P_2}\de \vartheta < \int_0^{t_1}f(\vartheta)\de \vartheta=\mathcal{H}^{N}(\Omega_1^{t_1})=\alpha,
\end{split}
\end{equation}
which is a contradiction.
\medskip

Now let us exploit the \textbf{Claim} to see that we cannot have $c>0$. If $c>0$, by the \textbf{Claim} we can find $t_1,t_2$ sufficiently small such that $\mathcal{H}^N(\Omega_1^{t_1})=\mathcal{H}^N(\Omega_2^{t_2})$, and $P_1t_1<P_2t_2$. Thus, for $t_1,t_2$ sufficiently small, the competitor $F:=E\cup \Omega_1^{t_1}\setminus \Omega_2^{t_2}$ has volume $V$, and
\[
\mathrm{Per}(F)-\mathrm{Per}(E)\leq f(t_1)-P_1+g(t_2)-P_2\leq c(P_1t_1-P_2t_2) <0,
\]
by using the claim and the fact that $c>0$. This is a contradiction with the fact that $E$ is isoperimetric.
\medskip

Then we have that $c=0$. Let us now distinguish two cases. 

If $X$ has at least two ends, then by the splitting theorem and the classification of $1$-dimensional $\CBB$ spaces \cite[Theorem 15.18]{AKP} we get that $X$ is as in (2). By the classification of isoperimetric sets in these spaces, see \cite{RosIsoperimetric}, we get that either $E$ is a geodesic ball with connected boundary, or, up to translation in the $\mathbb R$ factors, the set in (2). So this case is proved.

Let us now deal with the remaining case in which $X$ has one end. By \cref{lem:ConnCompRCD}, the set $\Omega\eqdef X\setminus \overline{E}$ is connected and unbounded. 
Applying the rigidity part in \cref{prop:C>0},  
and using the fact that $\Omega$ is connected, imply that the space $(Y,\dist_{\rm int})$ given by $\partial E$ endowed with intrinsic distance is connected, and thus $\partial E$ is connected, against the standing assumption of the proof.
Hence the two alternatives follow.

The rest of the statement follows from (1), (2) and \cref{thm:Regularity}.
\end{proof}

\begin{remark}
    The proof \cref{thm:IsopBoundaryConnected} works also in case $\partial X\neq\emptyset$, provided $\partial E$ has finitely many connected components.
\end{remark}

\begin{remark}\label{rem:Foliations}
Let $X$ be a $2$-dimensional $\CBB(0)$ metric space. Let $E$ be an isoperimetric set. As a consequence of \cref{thm:Regularity}, for every $x\in \partial E\setminus \partial X$ we have $(X,{r_i}^{-1}\dist,x,E)\to (\mathbb R^2,\dist_{\mathrm{eu}},0,\{x_2>0\})$ in the pointed Gromov--Hausdorff sense as $i\to +\infty$, see \cref{def:GHconvergence} and \cref{def:L1strong}. Hence, all the $x\in\partial E\setminus \partial X$ are regular points in $X$. As a consequence, if we take $X$ a $2$-dimensional $\CBB(0)$ metric space with singular points, then
\[
\bigcup_{\text{$E$ isoperimetric}} \partial E \neq X.
\]
We stress that there are examples \cite{OtsuShioya} of $2$-dimensional $\CBB(0)$ metric spaces with dense singular points. In particular it is in general not possible to foliate open subsets of a $\CBB(0)$ space by boundaries of isoperimetric sets.
\end{remark}

\medskip

The following \cref{thm:ExistenceLargeVolumesAsympCyl} and \cref{thm:AsintoticaProfilo}, which are a direct by-product of our study, give a full generalization of \cite[Theorem 6.20 and Theorem 6.21]{LeonardiRitore} to the case of Alexandrov spaces of arbitrary dimensions. Notice that this is a non-trivial generalization of the results in \cite{LeonardiRitore}, which dealt with the case of convex bodies in $\mathbb R^n$. Also, compare item (4) in \cref{thm:AsintoticaProfilo} with the analogous statements in \cite[Theorem 4.1 and Theorem 4.13]{RVcylindricallyBounded}.

\begin{theorem}\label{thm:ExistenceLargeVolumesAsympCyl}
Let $N\ge 2$ and let $(X,\dist,\mathcal{H}^N)$ be a noncompact {non-collapsed} $\RCD(0,N)$ metric measure space. Let us assume that all the pGH limits at infinity of $X$ are either $\mathbb R^{N}$, or $\mathbb R^{N-1}\times[0,+\infty)$.\\
Then $\mathcal{H}^N(B_1(x))\geq v_0$ for every $x\in X$, for some $v_0>0$, and isoperimetric regions exist for all volumes.
\end{theorem}

\begin{proof}
{The proof easily follows by direct method applying \cref{lem:IsoperimetricAtFiniteOrInfinite} and exploiting the classical Morgan--Johnson type estimate \cite[Proposition 3.3]{AntonelliPasqualettoPozzetta21} or its boundary version \cref{lem:MJHalfspaces}.}
\end{proof}

{As mentioned in the Introduction, after this work was posted on arXiv, X. Zhu \cite{XingyuZhu} partially addressed \cref{problemequivalencelineargrowth}. In particular, combining his results with the proof of \cref{thm:AsintoticaProfilo} we provide below, the author has obtained a version of \cref{thm:AsintoticaProfilo} on smooth manifolds with nonnegative Ricci curvature, linear volume growth, and volume of unit balls uniformly bounded from below away from zero \cite[Theorem 4.2]{XingyuZhu}.}

\begin{theorem}\label{thm:AsintoticaProfilo}
Let $N\ge 2$ and let $(X,\dist)$ be an $N$-dimensional noncompact $\CBB(0)$ metric space such that $\mathcal{H}^N(B_1(x))\geq v_0$ for every $x\in X$ for some $v_0>0$. Let us assume that there exists a pGH limit at infinity of $X$ that is isometric to $\mathbb R\times K$ for a compact $K$, and that $X$ has one end\footnote{If it has at least two ends, then it is isometric to $\mathbb R\times K$ by the splitting theorem.}. Then the following hold.
\begin{enumerate}
\item All the pGH limits at infinity of $X$ {are isometric to $\mathbb R\times K$. Denote $\mathfrak{V}\eqdef \mathcal{H}^{N-1}(K)$;}
\item There exists a constant {$V_0:=V_0(v_0,N,\mathfrak{V})$} such that isoperimetric regions exist for every volume $V\geq V_0$;
\item For every $V>0$ we have that 
    \[
    I(V)\leq \mathfrak{V}.
    \]
    If for some $V'$ the equality is attained, then it is attained for all $V\geq V'$ and there exists an isoperimetric set $\Omega$ of volume $\haus^N(\Omega)\ge V'$ such that $(\partial \Omega, \dist|_{\partial \Omega})$ is $\CBB(0)$ and $(X \setminus\Omega, \dist|_{X\setminus\Omega})$ is isometric to $(\partial\Omega\times [0,+\infty), \dist|_{\partial \Omega}\otimes\dist_{\rm eu})$.
\item {If $E_i$ is a sequence of isoperimetric sets with $\haus^N(E_i)\to+\infty$, $x_i \in \partial E_i$, and $(X,\dist,x_i)$ pGH-converges to $(\R\times K,\dist_{\rm eu} \otimes \dist_{K}, (0,k_0))$, then $E_i$ converges to $(-\infty,0)\times K$ in $L^1_{\rm loc}$ and $\partial E_i$ converges to $(-\infty,0)\times K$ in Hausdorff distance, in a realization of the pGH-convergence.} Moreover there exists a sequence of isoperimetric sets $\Omega_i\Subset \Omega_{i+1}$ such that $\cup_i \Omega_i = X$, and there holds
    \[
    \lim_{V\to+\infty} I(V) =  \mathfrak{V}.
    \]
\end{enumerate} 
\end{theorem}

\begin{proof}
First, thanks to \cref{thm:raggiosselimiteainfinito}, since $X$ has one end, there is a ray $\gamma$ such that $\dist(p,\gamma)\leq C$ for every $p\in X$. Thus, we can find $t_j\to +\infty$ such that $(X,\dist(\gamma(t_j))$ pGH-converges to $\mathbb R\times K$. Thus, item (1) directly comes from the fourth item in \cref{thm:ConvergenzaBusemann}. 

Let us now prove item (2). {
Let {$V_0:=V_0(v_0,N,\mathfrak{V})$} be the constant provided by \cref{lem:IsoperimetricLargeVolumesProduct}. 
Let $V\geq V_0$. If we take a minimizing (for the perimeter) sequence of bounded sets of finite perimeter of volume $V$ and we apply \cref{lem:IsoperimetricAtFiniteOrInfinite} we get that either all the mass stays in $X$ or it escapes to precisely one pGH limit at infinity. If all the mass stays in $X$ we are done, since we obtain an isoperimetric set of volume $V$. Otherwise, there exists a pGH limit at infinity $X_\infty\cong \mathbb R\times K$, due to item (1), and a set of finite perimeter $E\subset X_\infty$, which is isoperimetric for its own volume in $X_\infty$, and such that $\mathcal{H}^N(E)=V$, and $I_X(V)=\mathrm{Per}(E)$. Since $V\geq V_0$, from \cref{lem:IsoperimetricLargeVolumesProduct} we get that, up to translation in the $\mathbb R$ factor, $E=(0,a)\times K$, for some $a>0$. Now applying the last part of \cref{lem:MJCylindric} we have that there is a set $G\subset X$, which is a sublevel set of the Busemann function $F$, such that $\mathcal{H}^N(G)=V$, and $\mathrm{Per}(G)\leq \mathfrak{V}$. Thus 
\[
I_X(V)\leq \mathrm{Per}(G)<2\mathfrak{V} = \mathrm{Per}(E) = I_X(V),
\]
which is a contradiction. Thus, in the hypotheses above, in the minimization process there is no loss of mass at infinity, and thus an isoperimetric region of volume $V$ in $X$ always exists.}

Let us now prove item (3). Again by \cref{lem:MJCylindric} we get that for every $V>0$ there exists $G:=\{F<s\}\subset X$ such that $\mathcal{H}^N(G)=V$ and $\mathrm{Per}(G)\le \mathfrak{V}$. Thus $I(V)\leq \mathfrak{V}$ for every $V>0$. In addition, if equality holds for $V'>0$ then it holds for every $V\geq V'$ because $I$ is nondecreasing, see \cref{cor:IsoperimetricProfileRCD0N}.
Let us now prove the rigidity part in item (3). Assume there is $V'>0$ such that $I(V')=\mathfrak{V}$. Thus, applying \cref{lem:MJCylindric}, there exists $\bar s\in\mathbb R$ such that, being $F$ the Busemann function associated to $\gamma$,
\[
\mathfrak{V}= I(V')\leq \mathrm{Per}(\{F<\bar s\}) \leq \mathfrak{V}.
\]
By \eqref{eq:MonotoniaPerimetriBusemann} and \eqref{eq:LimitePerimetriBusemann} we have that $\{F<s\}$ is isoperimetric for any $s\ge \bar s$. Taking $\Omega\eqdef \{F<s_0\}$ for $s_0\ge \bar s$ given by \cref{lem:MJCylindric} implies the desired rigidity.

{Finally we prove item (4). Let $E_i$ be a sequence of isoperimetric sets with $\haus^N(E_i)\to+\infty$, and let $x_i \in \partial E_i$ such that $|ox_i|\ge |ox|$ for any $x \in \partial E_i$. Up to subsequence, assume that $(X,\dist,x_i)$ pGH-converges to $(\R\times K,\dist_{\rm eu} \otimes \dist_{K}, (0,k_0))$. Since $I$ is concave and nondecreasing, see \cref{cor:IsoperimetricProfileRCD0N}, and uniformly bounded by $\mathfrak{V}$, the Lipschitz constant of $I$ on $[V_0,+\infty)$ is bounded and tends to zero as $V_0\to+\infty$. Hence by \cref{cor:UniformRegularityIsop} we deduce that the sets $E_i$ are $(\Lambda_i,R_i)$-minimizers with $\Lambda_i \to 0$ and $R_i\to+\infty$.
 By the first item in \cref{rem:SemicontPerimeterConverging}, up to subsequences, $E_i$ converges in $L^1_{\rm loc}$ to a set $E \subset \R\times K$. Since $\Lambda_i \to 0$ and $R_i\to+\infty$, exploiting the approximation from the second item in \cref{rem:SemicontPerimeterConverging}, we have that $E$ is a perimeter minimizer in the sense of \cref{lem:PerimeterMinimizer}.\\
Moreover $\haus^N(E)=+\infty$ and $\Per(E)\le \liminf_i \Per(E_i) =\lim_{V\to+\infty} I(V)  \le \mathfrak V$ by lower semicontinuity, see \cref{rem:SemicontPerimeterConverging}. By \eqref{eqn:DensityEstimates}, the measure $\Per(E,\cdot)$ satisfies uniform $(N-1)$-dimensional density estimates, hence $\partial E$ must be bounded.\\
By the choice of $x_i$, \cref{lem:PerimeterMinimizer} implies that $E=(-\infty,0)\times K$. It now follows as in \cite[Theorem 2.43]{MoS21} that $\overline{E}_i \cap \overline{B}_R(x_i)$ (resp. $\partial E_i \cap B_R(x_i)$) converges in Hausdorff distance to $\overline{E} \cap \overline{B}_R((0,k_0))$ (resp. $\partial E \cap B_R((0,k_0)) = \{0\}\times K$) for any $R>2\,{\rm diam }\, K$ and that perimeter measures $\Per(E_i,\cdot) \res  B_R(x_i)$ converge in duality with continuous functions with bounded support in a realization of $(X,\dist,x_i) \to (\R\times K,\dist_{\rm eu} \otimes \dist_{K}, (0,k_0))$.

Fix a ray $\gamma$ from an origin $o\in X$. Let $R_0>2(2+\, {\rm diam}\, K)$ and let $t_i$ be such that $|x_i\gamma(t_i)|=\dist(x_i,\gamma)$. Let $F$ be the Busemann function associated to $\gamma$. By \cref{prop:StimaFq} we have that $t_i\le F(x_i)\le t_i+\eps_i$ for large $i$, with $\eps_i\to0$. By the second item in \cref{thm:ConvergenzaBusemann} and the choice of $R_0$ we have $\{t_i-2<F<t_i+2\} \subset B_{R_0}(x_i)$ for large $i$.\\
 Hence convergence in Hausdorff distance of the boundaries to the boundary $\{0\}\times K$ implies that
\begin{equation}\label{eq:zzwwzz}
\partial E_i \cap B_{R_0}(x_i) \subset \{t_i-\delta_i<F<t_i+\delta_i\} \subset  \{t_i-2<F<t_i+2\} \subset B_{R_0}(x_i),
\end{equation}
for large $i$, for some $\delta_i\searrow0$.  We claim that $\{F<t_i-2\}\subset E_i$. Indeed we have convergence in Hausdorff distance of $\overline{B}_{R_0}(x_i)\setminus E_i$ to $\overline{B}_{R_0}((0,k_0))\setminus E$. By \cref{lem:ConnCompRCD}, if we could join a point $z \in \{F<t_i-2\} \setminus \overline{E}_i$ to a point in $X\setminus \overline{B}_{|ox_i|+1}(o)$ with a curve contained in $X\setminus \overline{E}_i$ for any $i$ large, passing to the limit the sequence of curves we would get a limit curve contained in $X\setminus E$ that intersects $\{-2\}\times K$, that is not possible.\\
Hence we proved that $\partial E_i = \partial E_i \cap B_{R_0}(x_i)$ for large $i$, hence the convergence of perimeter measures in duality with continuous functions with bounded support implies $\lim_i \Per(E_i) = \Per (E) = \haus^{N-1}(K)$, and thus $\lim_{V\to+\infty} I(V) = \mathfrak V$. Finally, since also $\{F<t_i-2\} \subset E_i$, the desired exhaustion $\Omega_i$ is constructed by taking a suitable subsequence of $E_i$.}
\end{proof}

With techniques similar to the previous \cref{thm:AsintoticaProfilo}, we can also prove the following result.
\begin{proposition}\label{prop:ProfiloCostiffSPlitCilindro}
    Let $(X,\dist)$ be an $N$-dimensional noncompact $\CBB(0)$ metric space with $\mathcal{H}^N(B_1(x))\geq v_0>0$ for every $x\in X$. The following are equivalent.
    \begin{enumerate}
        \item The isoperimetric profile $I$ is eventually constant;
        \item Either $X$ is isometric to $\mathbb R\times K$ for some compact $K$, or there exists an isoperimetric set $A$ such that $(\partial A, \dist|_{\partial A})$ is $\CBB(0)$ and $(X \setminus A, \dist|_{X\setminus A})$ is isometric to $(\partial A\times [0,+\infty), \dist|_{\partial A}\otimes\dist_{\rm eu})$;
        \item Either $X$ is isometric to $\mathbb R\times K$ for some compact $K$, or there exists an open bounded set $A$ such that $(X \setminus A, \dist|_{X\setminus A})$ is $\CBB(0)$ and it is isometric to $(\partial A\times [0,+\infty), \dist|_{\partial A}\otimes\dist_{\rm eu})$.
    \end{enumerate}
In particular, if any of the above holds, either $X$ is isometric to $\mathbb R\times K$ for some compact $K$ and there exists $V_0>0$ such that for every $V\geq V_0$ isoperimetric regions of volume $V$ exist coincide with $(0,a)\times K$ for some $a>0$, up to translations; or there exists $V_0>0$ such that for every $V\geq V_0$ isoperimetric regions of volume $V$ exist, are unique, and coincide with $\{x \st \dist(x,A) < b_V\}$ for some $b_V\geq 0$.
\end{proposition}

\begin{proof}
(2)$\Rightarrow$(1). If $X$ is isometric to $\mathbb R\times K$ for a compact $K$, then the result follows from \cref{lem:IsoperimetricLargeVolumesProduct}. In the other case we can argue as follows. Since all the pGH limits at infinity of $X$ are isometric to $\partial A\times \mathbb R$, from item (2) of \cref{thm:AsintoticaProfilo} we have existence of isoperimetric sets for large volumes. Moreover, from item (3) of \cref{thm:AsintoticaProfilo}, we have that $I\leq \mathcal{H}^{N-1}(\partial A)$. Let us show that $I(V)=\mathcal{H}^{N-1}(\partial A)$ for every $V\geq V_0$, where $V_0$ is sufficiently big. We just sketch the argument, since it relies on a classical symmetrization technique already present in \cref{AppendixA}. By taking $V_0$ sufficiently big, we have that all the isoperimetric sets $E$ of volume $V\geq V_0$ intersect $X\setminus \overline{A}$ in a $\mathcal{H}^N$-positive measured set. By symmetrization, see \cref{lemma:Symmetrization}, since $E$ is isoperimetric, we can assume that, for some function $f$, 
\[
E\cap X\setminus \overline A=\{(x,t):x\in \pi(E\setminus \overline A), 0<t<f(x)\},
\]
where $\pi$ is the projection onto $\partial A$ in the product $\partial A\times [0,+\infty)$, where we identify $\partial A$ and $\partial A \times\{0\}$. By the same argument in the beginning of \cref{lem:IsoperimetricLargeVolumesProduct} we have that, for $V_0$ sufficiently big, $\mathcal{H}^{N-1}(\partial A\setminus\pi(E\setminus \overline A))=0$, and $f\in\mathrm{BV}_{\mathrm{loc}}(\partial A)$. Now, as a consequence of the inequality in \cite[Theorem 5.1, item (b)]{APS15}, we get that 
\[
\mathrm{Per}(E)\geq \mathrm{Per}(E,X\setminus\overline A)\geq \mathcal{H}^{N-1}(\partial A).
\]
Finally, for an isoperimetric set $E$ of volume $V\ge V_0$, since $\mathrm{Per}(E)=I(V)\leq\mathcal{H}^{N-1}(\partial A)$, we get the equality in the previous inequality. Thus, by inspecting the equality cases, taking into account also \cite[Theorem 5.1, item (b)]{APS15}, we see that for $V_0$
sufficiently big, we have that for every $V\geq V_0$ the equality $I(V)=\mathcal{H}^{N-1}(\partial A)$ holds, and the unique isoperimetric regions of volume $V$ is $\{x \st \dist(x,A) < b_V\}$ for some $b_V\geq 0$. The last assertion comes since the last argument implies that we have the equality in $\mathrm{Per}(E)=\mathrm{Per}(E,X\setminus\overline A)$ and in the inequality in \cite[Theorem 5.1, item (b)]{APS15}.

(1)$\Rightarrow$(2). Let us assume there exists $V_0$ such that for every $V\geq V_0$ we have $I(V)=C$ for some constant $C>0$. Let us consider the volume $V_0+1$. If there exists an isoperimetric set $E$ of volume $V_0+1$, then, since $I'(V_0+1)=0$, we get that its mean curvature barrier $c=0$, see the last part of \cref{lem:IsoperimetricAtFiniteOrInfinite}.
{From the fact that $I$ is non-decreasing, cf. \cref{cor:IsoperimetricProfileRCD0N}, and by using \eqref{eqn:ExtPerEsterno}, we have that $E_t\eqdef \{\dist(\cdot,E)<t\}$ is isoperimetric for any $t>0$. 
Moreover, since $c=0$, we can use the rigidity in \cref{prop:C>0}, see also \cref{thm:SplittingKetterer}, on $X\setminus\overline{E}$.
If $X$ has at least two ends, we thus get that $X$ is isometric to a product $\R\times K$ from the splitting theorem, and $K$ is isometric to a connected component $\partial E_\alpha$ of $(\partial E, \dist_{\rm int})$ from the previous rigidity. Thus $K$ is compact. If instead $X$ has one end, then $\Omega\eqdef X\setminus\overline{E}$ is connected by \cref{lem:ConnCompRCD}. In the notation of \cref{thm:SplittingKetterer}, we have that $\widetilde\dist(\cdot, E)= \dist(\cdot, E)$ on $\Omega$ by definition of intrinsic distance, hence $E_t=\{\widetilde\dist(\cdot, E)< t\}$. Recalling that $\widetilde\dist=\dist$ on $\{\dist_E>t_0\}$ for some $t_0$, the claim follows by taking $A\eqdef E_t$ for $t>t_0$.}

If there is no isoperimetric set of volume $V_0+1$, then a minimizing sequence (for the perimeter) of sets of volume $V_0+1$ must converge to an isoperimetric set with mean curvature barrier $c=0$ in a pGH limit at infinity of space, see \cref{lem:IsoperimetricAtFiniteOrInfinite}. From \cref{lem:EveryLimitAtInfinitySplits}, we get that such a limit splits as $\mathbb R\times K$. Hence, arguing as in the first part of the implication (1)$\Rightarrow$(2), we conclude that $K$ is compact as above.
Thus by \cref{thm:AsintoticaProfilo} we get existence of isoperimetric sets for sufficiently big volumes, and then we can argue again as in the first part of the proof of (1)$\Rightarrow$(2) to conclude.

(2)$\Longleftrightarrow$(3). The nontrivial implication (3)$\Rightarrow$(2) follows by arguing as in (2)$\Rightarrow$(1).
\medskip

The last part of the statement is a direct by-product of the proof above.
\end{proof}

\appendix

\section{Large isoperimetric sets in cylinders}\label{AppendixA}

{
Let us briefly recall the Steiner symmetrization process in our setting. Let us consider the metric measure space $(\mathbb R\times X,\dist_{\mathrm{eu}}\otimes \dist_X,\mathcal{H}^1\otimes\mathfrak{m})$, and let $E$ be a subset of $\mathbb R\times X$. Given $x\in X$, $t\in\mathbb R$, we denote
\[
E_x:=\{s\in\mathbb R:(s,x)\in E\}, \qquad E^t:=\{y\in X:(t,y)\in E\}.
\]
We denote 
\[
\mathrm{sym}(E):=\left\{(t,y)\in \mathbb R\times X:|t|\leq \frac{\mathcal{H}^1(E_y)}{2}\right\}.
\]
We say that a set $E\subset\mathbb R\times X$ is {\em normalized} when, for every $x\in X$, $E_x$ is an open segment whose midpoint is $0$. It is easy to notice that when $E$ is normalized, then 
\begin{equation}\label{eqn:MonotonicitySlices}
0\leq t\leq s \Rightarrow E^s\subset E^t.
\end{equation}

There holds the following result, compare with \cite[Theorem 14.4]{MaggiBook} and \cite[Section 5.1]{RitoreBook}. Since the proof is classical, we just sketch it.
\begin{lemma}\label{lemma:Symmetrization}
    Let $(\mathbb R\times X,\dist_{\mathrm{eu}}\otimes \dist_X,\mathcal{H}^1\otimes\mathfrak{m})$ be a metric measure space. Let $E\subset \mathbb R\times X$ be a measurable set. Then the following hold:
    \begin{itemize}
        \item $\mathrm{sym}(E)$ is measurable, and $(\mathcal{H}^1\otimes\mathfrak{m})(\mathrm{sym}(E))=(\mathcal{H}^1\otimes\mathfrak{m})(E)$;
        \item $\mathrm{Per}(\mathrm{sym}(E))\leq \mathrm{Per}(E)$.
    \end{itemize}
\end{lemma}
\begin{proof}

{The first item is a consequence of Fubini--Tonelli Theorem. The proof of the second item can be reached arguing verbatim as in \cite[Proposition 3]{MorganHoweHarman}, cf. with \cite[Theorem 5.6 and Corollary 5.7]{RitoreBook}. In particular, first one obtain that for every $r>0$ the volume of the $r$-enlargment $(\mathrm{sym}(E))_r$ is less or equal than the volume of the $r$-enlargment $E_r$. Thus, the lower Minkowski content of $\mathrm{sym}(E)$ is less or equal than the lower Minkowski content of $E$. By finally using the result in \cite[Theorem 3.6]{AmbrosioGigliDimarino} we get the sought conclusion.}
\end{proof}

\begin{remark}
    If $(K,\dist,\mathcal{H}^N)$ is a compact {non-collapsed} $\RCD(k,N)$ metric measure space such that $\inf_{x\in K} \haus^N(B_1(x)) \ge v_0' > 0$, then we can find a constant $C:=C(\mathcal{H}^N(K),k,N, v_0')$ such that, for every measurable set $E$ with $\mathcal{H}^N(E)\leq \mathcal{H}^N(K)/2$, we have
    \begin{equation}\label{eqn:Isop}
    \mathrm{Per}(E)\geq C\max\left\{\mathcal{H}^N(E),\left(\mathcal{H}^N(E)\right)^{(N-1)/N}\right\}.
    \end{equation}
This is an immediate consequence of the isoperimetric inequality for small volumes in this setting, see, e.g., \cite[Proposition 3.20 and Remark 3.21]{AntonelliPasqualettoPozzetta21}, and of the uniform positive lower bound satisfied by the isoperimetric profile under these assumptions, see, e.g., \cite[Corollary 4.14]{APPSa}.\\
Observe that if $X=(\R\times K, \dist_{\rm eu}\otimes \dist_K, \haus^{N+1})$ is $\RCD(k,N+1)$ and $\inf_{x \in X} \haus^{N+1}(B_1(x)) \ge v_0>0$, then $\inf_{z \in K}\haus^N(B_1(z))\ge v_0'>0$ for $v_0'=v_0/2$, indeed $B_1((0,z))\subset (-1,1)\times B^K_1(z)$ for any $z \in K$.
\end{remark}

\begin{remark}
    Let $(X\times\mathbb R,\dist_X\otimes\dist_{\mathrm{eu}},\mathfrak{m}\otimes\mathcal{H}^1)$ be a metric measure space. Let $E\subset X\times\mathbb R$ be a measurable set. It directly follows from the definition that, for every $a,b\in\mathbb R\cup\{-\infty,+\infty\}$, we have 
    \begin{equation}\label{eqn:EstimatePer}
    \mathrm{Per}(E,X\times (a,b))\geq \int_a^b \mathrm{Per}(E^s)\de s.
    \end{equation}
    Indeed, let us take $f_i\in \mathrm{LIP}_{\mathrm{loc}}(X\times(a,b))$ such that $f_i\to \chi_E$ in $L^1_{\mathrm{loc}}(X\times (a,b))$. Thus for $\mathcal{H}^1$-a.e. $s\in (a,b)$ we have that $f_i|_{E^s}\to \chi_{E^s}$ in $L^1_{\mathrm{loc}}(X)$. We further have, by using Fatou's lemma and the fact that for every $s\in (a,b)$ we have $\mathrm{lip}f_i\geq \mathrm{lip}(f_i|_{E_s})$ 
    \begin{equation}
    \begin{split}
    \liminf_{i\to +\infty}\int_{X\times (a,b)} \mathrm{lip}f_i\de (\mathfrak{m}\otimes\mathcal{H}^1) &= \liminf_{i\to +\infty}\int_a^b \int _X \mathrm{lip}f_i(x,s)\de \mathfrak{m}(x)\de \mathcal{H}^1(s)\\
    &\geq \int_a^b \int _X \liminf_{i\to +\infty}\mathrm{lip}f_i(x,s)\de \mathfrak{m}(x)\de \mathcal{H}^1(s)\\
    &\geq \int_a^b \int _X \liminf_{i\to +\infty}\mathrm{lip}(f_i|_{E_s})(x,s)\de \mathfrak{m}(x)\de \mathcal{H}^1(s)\\
    &\geq \int_a^b \mathrm{Per}(E_s)\de s.
    \end{split}
    \end{equation}
Thus taking the infimum in the left-hand-side of the previous inequality we get the sought claim.
\end{remark}

The proof of the next two results are inspired from \cite[Theorem 3.9 and Corollary 4.12]{RVcylindricallyBounded}, cf. \cite[Section 7.2]{RitoreBook}.

\begin{proposition}\label{prop:IsopRinorm}
    Let $X=(\mathbb R\times K,\dist_{\mathrm{eu}}\otimes \dist_K,\mathcal{H}^{N+1})$ be an $(N+1)$-dimensional {$\RCD(k,N+1)$ metric measure space such that $\inf_{x \in X}\haus^{N+1}(B_1(x)) \ge v_0>0$, with $K$ compact. Then there exists a constant $\ell:=\ell(\mathcal{H}^N(K), k,N,v_0)$} such that the following holds. 

    For every normalized isoperimetric set $E\subset \mathbb R\times X$ with $\mathcal{H}^{N}(K\setminus E^0)>0$, we have 
    \begin{equation}
        \mathrm{Per}(E)\geq \ell\cdot \mathcal{H}^{N+1}(E).
    \end{equation}
\end{proposition}
\begin{proof}
Notice that $(K,\dist_K,\mathcal{H}^N)$ is a compact {$\RCD(k,N)$} metric measure space, and $\mathcal{H}^{N+1}=\mathcal{H}^1\otimes \mathcal{H}^N$.
Moreover, let us notice that since $E$ is normalized, $E$ is the intersection of the subgraph of $f\in \mathrm{BV}_{\mathrm{loc}}(K)$ and the complement of the subgraph of $-f$, see \cite[Theorem 3]{AntonelliBrenaPasqualetto2}.

Let us consider $\tau\geq 0$ with the property that 
\[
\mathcal{H}^N(E^t)\geq \frac{\mathcal{H}^N(K)}{2} \qquad \forall t\in [0,\tau), \qquad\qquad \mathcal{H}^N(E^t)\leq \frac{\mathcal{H}^N(K)}{2} \qquad \forall t\in [\tau,+\infty).
\]
By using \eqref{eqn:EstimatePer}, \eqref{eqn:Isop}, and the co-area formula we get
\begin{equation}\label{eqn:Finally2}
\mathrm{Per}(E)\geq \mathrm{Per}(E,(\tau,+\infty)\times K)\geq C\mathcal{H}^{N+1}(E\cap (\tau,+\infty)\times K),
\end{equation}
{for $C=C(\haus^N(K),k,N,v_0)>0$.}
Thus, if $\tau=0$, the previous inequality, together with the fact that $E$ is normalized, gives the sought conclusion. 

Let us then assume, from now on, that $\tau>0$. Notice that, by comparing $E$, which is an isoperimetric set, with a slab with the same volume, we get that 
$$
\mathrm{Per}(E)\leq 2\mathcal{H}^{N}(K).
$$ 
By using the inequality in \cite[Theorem 5.1, item (b)]{APS15}, we also get that 
\begin{equation}\label{eqn:Finally}
\mathrm{Per}(E,(t,+\infty)\times K)\geq \mathcal{H}^N(E^t),
\end{equation}
for every $t\geq 0$. Using the previous two inequalities we get that for almost every $t\in [0,\tau)$ we have 
\begin{equation}\label{eqn:UnaStima}
\begin{split}
\mathcal{H}^N(K\setminus E^t)&=\mathcal{H}^N(K)-\mathcal{H}^N(E^t) \\
&\geq \frac{\mathrm{Per}(E)}{2} -\mathrm{Per}(E,(t,+\infty)\times K) \geq \mathrm{Per}(E,(0,t)\times K),
\end{split}
\end{equation}
where in the last equality we are using that $\mathrm{Per}(E)\geq 2\,\Per(E,(0,+\infty)\times K)$ from the fact that $E$ is normalized. By using \eqref{eqn:UnaStima}, \eqref{eqn:EstimatePer}, and \eqref{eqn:Isop} together with the fact that $\mathcal{H}^N(K\setminus E^t)\leq \mathcal{H}^N(K)/2$ for every $t\in [0,\tau)$, we get that, calling $f(t):=\mathcal{H}^N(K\setminus E^t)$, we have
\begin{equation}
f(t)\geq C\int_0^t f(s)^{\frac{N-1}{N}} \de s,
\end{equation}
{for some $C=C(\haus^N(K),k,N,v_0)>0$.}
By integrating the previous inequality and using that $f(t)>0$ for every $t\in [0,\tau)$ by assumption, we get 
\begin{equation}\label{eqn:lalala}
\mathcal{H}^N(K)\geq f(\tau) \geq C'\tau ^N \qquad\Rightarrow\qquad \tau \leq \bar{C}\left(\mathcal{H}^N(K)\right)^{1/N},
\end{equation}
where $\bar{C}=\bar{C}(\haus^N(K),k,N,v_0)>0$. Finally, since the slices $E^t$ are decreasing, we get, from the coarea formula, from the estimate in \eqref{eqn:lalala}, and from \eqref{eqn:Finally} with $t=0$, that
\begin{equation}\label{eqn:lalala2}
\mathcal{H}^{N+1}(E\cap (0,\tau)\times K)=\int_0^\tau \mathcal{H}^N(E^s)\de s\leq \tau \mathcal{H}^N(E^0)\leq \tilde C\mathrm{Per}(E),
\end{equation}
{for $\tilde C=\tilde C(\haus^N(K),k,N,v_0)>0$.} By summing together \eqref{eqn:Finally2} and \eqref{eqn:lalala2} we finally get the sought claim.
\end{proof}

For results related to the forthcoming Lemma, see \cite{DuzaarSteffen96, Castro22}. In \cite{DuzaarSteffen96} \cref{lem:IsoperimetricLargeVolumesProduct} was first proved for compact Riemannian manifolds $K$, and in \cite{Castro22} a version of \cref{lem:IsoperimetricLargeVolumesProduct} was proved for smooth cylinders with density constant along the real factor.
\begin{lemma}\label{lem:IsoperimetricLargeVolumesProduct}
Let $X=(\mathbb R\times K,\dist_{\mathrm{eu}}\otimes \dist_K,\mathcal{H}^{N+1})$ be an $(N+1)$-dimensional {$\RCD(k,N+1)$ metric measure space such that $\inf_{x \in X}\haus^{N+1}(B_1(x)) \ge v_0>0$, with $K$ compact. There exists a volume $V_0:=V_0(\mathcal{H}^N(K),k,N,v_0)>0$} such that for every $V\geq V_0$ the unique, up to translation along the factor $\mathbb R$, isoperimetric region of volume $V$ is $(0,a)\times K $ for suitable $a>0$.
\end{lemma}

\begin{proof}
Let $\pi:\mathbb R\times K\to K$ be the projection on $K$. Notice that $\pi$ is $1$-Lipschitz. Let $f:\mathbb R\times K\to \mathbb R$ be the projection on $\mathbb R$.

Let us take $V_0:=4\mathcal{H}^N(K)/\ell$, where $\ell$ is chosen as in \cref{prop:IsopRinorm}. Since every pGH limit at infinity of $\mathbb R\times K$ is isometric to $\mathbb R\times K$, we get that, by exploiting the asymptotic mass decomposition result in \cref{thm:MassDecompositionINTRO}, for every $V>0$ the isoperimetric problem has at least one solution for the volume $V$ and it is bounded, cf. the first part of \cref{thm:Regularity}. Let $E$ be one isoperimetric set with volume $V\geq V_0$. By the topological regularity of $E$, see the first part of \cref{thm:Regularity}, up to take a representative of $E$, $E$ is open, and the essential boundary of $E$ coincides with its topological boundary $\partial E$.

By \cref{lemma:Symmetrization} we get that $\mathrm{sym}(E)$ is a normalized isoperimetric set and 
\begin{equation}\label{eqn:Slab}
\mathrm{Per}(E)=\mathrm{Per}(\mathrm{sym}(E))\leq 2\mathcal{H}^N(K),
\end{equation}
by comparing with a slab with the same volume. From \cref{prop:IsopRinorm} and from the choice of $V_0$, we get that $\mathcal{H}^N(K\setminus (\mathrm{sym}(E))^0)=0$. Thus $\mathcal{H}^N(K\setminus \pi(E))=0$, since $\pi(E)=\pi(\mathrm{sym}(E))=(\mathrm{sym}(E))^0$. Thus for $\mathcal{H}^N$-almost every $x\in K$, $|\pi^{-1}(x)\cap E|\geq 1$, and then $|\pi^{-1}(x)\cap \partial E|\geq 1$, since $E$ is bounded. 

{
We now claim that for $\mathcal{H}^N$-almost every $x\in K$, $\mathcal{H}^1(\pi^{-1}(x)\cap  E)> 0$. Assume it is not the case, and let us call $X':=\{x\in K:\mathcal{H}^1(\pi^{-1}(x)\cap  E)= 0\}$. We have that $X'$ is measurable and $\mathcal{H}^{N}(X')>0$. By the rectifiability results for sets of locally finite perimeter recalled at the beginning of \cref{sec:RCD} we have that for $\mathcal{H}^N$-almost every $p\in \partial E$, the following holds
\begin{equation}\label{eqn:DensityHalf}
\lim_{r\to 0^+}\frac{\mathcal{H}^{N+1}(B_r^X(p)\cap E)}{\omega_{N+1}r^{N+1}}=\frac{1}{2}.
\end{equation}
Thus, since $\pi$ is Lipschitz, there exists $x\in X'$ a regular point of $\mathcal{H}^N$-density $1$ in $X'$ such that there is $p_x\in\partial E$ for which $\pi(p_x)=x$, and $p_x$ satisfies \eqref{eqn:DensityHalf}. Let us now fix this $x$, and let $\vartheta>0$ be small enough such that 
\begin{equation}\label{eqn:EndEnd}
    \mathcal{H}^{N+1}(B_\vartheta^X(p_x)\cap E)\geq \frac{1}{4}\omega_{N+1}\vartheta^{N+1}, \qquad {\mathcal{H}^N(B_\vartheta^X(x)\setminus X')\leq \frac{1}{16}\omega_{N+1}\vartheta^N},
\end{equation}
and let us find a contradiction. First notice that, by Fubini,
\[
\mathcal{H}^{N+1}(E\cap\pi^{-1}(B_\vartheta^K(x)\cap X'))=\int_{B_\vartheta^K(x)\cap X'}\mathcal{H}^1(\pi^{-1}(y)\cap E)\de \mathcal{H}^N(y)=0.
\]
Thus, by also exploiting the previous equality, \eqref{eqn:EndEnd}, and the coarea formula, we get
\begin{equation}\label{eqn:Chain}
    \begin{split}
    \frac{1}{4}\omega_{N+1}\vartheta^{N+1}&\leq \mathcal{H}^{N+1}(B^X_\vartheta(p_x)\cap E)=\mathcal{H}^{N+1}(B^X_\vartheta(p_x)\cap E\cap \pi^{-1}(B^K_\vartheta(x)))\\
    &= \mathcal{H}^{N+1}(B^X_\vartheta(p_x)\cap E\cap \pi^{-1}(B^K_\vartheta(x)\setminus X'))\\
    &\leq \int_{f(p_x)-\vartheta}^{f(p_x)+\vartheta}\mathcal{H}^{N}(\{f=t\}\cap E\cap \pi^{-1}(B^K_\vartheta(x)\setminus X'))\\
    &\leq 2\vartheta\mathcal{H}^{N}(B^K_\vartheta(x)\setminus X') \leq \frac{1}{8}\omega_{N+1}\vartheta^{N+1},
    \end{split}
\end{equation}
which is the sought contradiction. We remark that in the second-last inequality we are using that 
\[
\pi(\{f=t\}\cap E\cap \pi^{-1}(B^K_\vartheta(x)\setminus X'))\subset B^K_\vartheta(x)\setminus X',
\]
and then 
\[
\mathcal{H}^N(\{f=t\}\cap E\cap \pi^{-1}(B^K_\vartheta(x)\setminus X'))=\mathcal{H}^N(\pi(\{f=t\}\cap E\cap \pi^{-1}(B^K_\vartheta(x)\setminus X')))\leq \mathcal{H}^N(B^K_\vartheta(x)\setminus X').
\]
}

We thus have proved that for $\mathcal{H}^N$-almost every $x\in K$, $\mathcal{H}^1(\pi^{-1}(x)\cap  E)> 0$. This now readily implies that for $\mathcal{H}^N$-almost every $x\in K$, $|\pi^{-1}(x)\cap \partial E|\geq 2$, again using that $E$ is bounded. We now claim that for $\mathcal{H}^N$-almost every $x\in K$, $|\pi^{-1}(x)\cap \partial E|= 2$. Indeed, if not, there is $A\subset K$ with $\mathcal{H}^N(A)>0$ such that $|\pi^{-1}(x)\cap \partial E|\geq 3$ for all $x\in A$. By using the rectifiability of the boundary, the fact that $\mathrm{Per}(E,\cdot)=\mathcal{H}^{N-1}\res \partial E$ and the area formula for metric spaces \cite[Theorem 8.2]{Ambrosio2000RectifiableSpaces} applied to the projection map $\pi$, which has Lipschitz constant $=1$, we infer
\[
\mathrm{Per}(E,\mathbb R\times A)\geq 3\mathcal{H}^N(A).
\]
Moreover, by the same reasoning, we always have $\mathrm{Per}(E,\mathbb R\times (K\setminus A))\geq 2\mathcal{H}^N(K\setminus A)$, since for $\mathcal{H}^N$-almost every $x\in K$, $|\pi^{-1}(x)\cap \partial E|\geq 2$. Thus we infer, by summing the previous two inequalities, that $\mathrm{Per}(E)>2\mathcal{H}^N(K)$, which is a contradiction with \eqref{eqn:Slab}.

Hence $|\pi^{-1}(x)\cap \partial E|=2$, and $\mathcal{H}^1(\pi^{-1}(x)\cap E)>0$ for $\mathcal{H}^N$-a.e. $x\in K$. Thus, for $\mathcal{H}^N$-a.e. $x\in K$, $E_x$ is an open and bounded interval. Thus $E$ is the intersection of the subgraph of a function $v\in\mathrm{BV}_{\mathrm{loc}}(K)$ with the complement of the subgraph of a function $u\in\mathrm{BV}_{\mathrm{loc}}(K)$, with $u<v$, see \cite[Theorem 3]{AntonelliBrenaPasqualetto2}. By using \cite[Item (b) of Theorem 5.1]{APS15}, we can finally write
\[
2\mathcal{H}^N(K)\geq \mathrm{Per}(E) = \left(\sqrt{1+g_u^2}+\sqrt{1+g_v^2}\right)\mathcal{H}^N(K)+(|Du|^s+|Dv|^s)(K)\geq 2\mathcal{H}^N(K).
\]
From the last chain of inequalities we infer that $u,v$ are constant and thus $E$ is a slab.
\end{proof}
}

We conclude with a characterization of perimeter minimizers with compact boundary in $\RCD$ cylinders.

\begin{lemma}\label{lem:PerimeterMinimizer}
    Let $X=(\mathbb R\times K,\dist_{\mathrm{eu}}\otimes \dist_K,\mathcal{H}^{N+1})$ be an $(N+1)$-dimensional $\RCD(k,N+1)$ metric measure space, with $K$ compact. Let $E\subset X$ be a perimeter minimizer, i.e., a set of finite perimeter such that $\Per(E)\le \Per(F)$ whenever $E\Delta F$ is essentially bounded in $X$. Up to choice of the representative $E=E^{(1)}$, assume that $\partial E$ is bounded and that $\haus^{N+1}(E)=+\infty$.
    
    Then, for some $t_0\in\mathbb R$, either $E=(-\infty,t_0)\times K$ or $E=(t_0, +\infty)\times K$.
\end{lemma}

\begin{proof}
Like in the case of isoperimetric sets, it follows from the regularity theory in \cite{ShanmugalingamQuasiminimizers} (see also \cite{MoS21}) and from the representation for the perimeter measure in \cite{BPSGaussGreen} that $E^{(1)}$ and $E^{(0)}$ are open and that $\Per(E,\cdot)= \haus^N\res \partial E$.
Without loss of generality, assume that $\partial E \subset (0,t_1)\times K$ for some $t_1>0$. Hence either $(-\infty,0)\times K \subset E$ or $(-\infty,0)\times K \subset X\setminus\overline{E}$; similarly, either $(t_1,+\infty)\times K \subset E$ or $(t_1,+\infty)\times K \subset X\setminus\overline{E}$. Without loss of generality assume that $(-\infty,0)\times K \subset E$. Hence $(t_1,+\infty)\times K \subset X\setminus\overline{E}$, otherwise $E\cup (-1,t_1+1)\times K$ would be the whole space, contradicting the minimality of $E$. Moreover $P(E) \le P(E \cup  (-1,t_1+1)\times K) = \haus^N(K)$.\\
Let $\pi:\R\times K \to \{t_1\}\times K$ be the projection on the slice $\{t_1\}\times K$. Hence $\pi|_E$ is surjective over $\{t_1\}\times K$. Arguing as in the proof of \cref{lem:IsoperimetricLargeVolumesProduct}, the fact that $P(E) \le \haus^N(K)$ implies that for $\haus^N$-a.e. point $(t_1,k) \in \{t_1\}\times K$ there exists at most one point in $\pi^{-1}((t_1,k)) \cap \partial E$. Hence \cite[Theorem 3]{AntonelliBrenaPasqualetto2} implies that $E$ is represented as the subgraph of a $BV_{\rm loc}(K)$ function $f:K\to \R$. Therefore $P(E) \le \haus^N(K)$ and the inequality in \cite[Theorem 5.1, item (b)]{APS15} implies that the total variation of $f$ is zero, and thus $f$ is constant.
\end{proof}

\section{Geodesics in products and cylinders}\label{AppendixB}

The next observation characterizes rays in cylinders over compact spaces, and its proof, which involves classical arguments in Metric Geometry, is left to the reader.

\begin{proposition}\label{prop:GeodeticheProdotto}
Let $(X_1,\dist_1)$ and $(X_2,\dist_2)$ be metric spaces and let $(X,\dist)$ be the product metric space $X=X_1\times X_2$ endowed with the product distance $\dist^2((x_1,x_2),(y_1,y_2))\eqdef \dist_1^2(x_1,y_1) +\dist_2^2(x_2,y_2)$ for any $(x_1,x_2),(y_1,y_2) \in X$. Assume that $(X,\dist)$ is $\CBB(0)$. Let $\gamma=(\gamma_1,\gamma_2):[0,T]\to X$ be a geodesic.
Then the following hold.
\begin{itemize}
    \item $L(\gamma_i)=\dist_i(\gamma_i(0),\gamma_i(T))$ for $i=1,2$.
    
    \item If $|\gamma'|(t)=1$ for a.e. $t$, then there exist $c_1,c_2\ge0$ such that $|\gamma'_i|(t)=c_i$ for a.e. $t$ for $i=1,2$
    
    \item If $X=\R\times K$, $K$ is compact and $\sigma:[0,+\infty)\to X$ is a ray, then there is $(t_0,k_0) \in X$ such that $\sigma(t)=(t_0+ t,k_0)$ or $\sigma(t)=(t_0- t,k_0)$. 
\end{itemize}
\end{proposition}

\section{Concavity of distance from geodesics in dimension $2$}\label{AppendixC}

We shall prove that the distance function from a geodesic on a $2$-dimensional $\CBB(0)$ space is concave. As pointed us to the authors by A. Lytchak, the forthcoming \cref{prop:ConcavityDistance} can also be obtained as a consequence of \cite[Corollary 2.8]{WornerPhD}, which in turn relies on \cite{Perelman91}. Indeed, one can consider the region bounded by the curves $\sigma$ and $\gamma$ in \cref{prop:ConcavityDistance} and by geodesics joining their endpoints. The assumptions imply that such set is an Alexandrov space for which $\gamma$ is a boundary stratum, and thus \cite[Corollary 2.8]{WornerPhD} applies. However, we include a proof based on classical comparison arguments.

\begin{proposition}\label{prop:ConcavityDistance}
Let $(X,\dist)$ be a $2$-dimensional $\CBB(0)$ metric space. Let $\sigma:(a,b)\to X $ and $\gamma:[0,T]\to X$ be geodesics such that $0<\dist(\sigma(t),\gamma)<\min\{\dist(\sigma(t),\gamma(0)), \dist(\sigma(t),\gamma(T))\}$, for any $t \in (a,b)$. Then $t\mapsto \dist(\sigma(t),\gamma)$ is concave.
\end{proposition}

{
We need to recall some basic facts about the first variation formula on Alexandrov spaces. Let $X, \sigma, \gamma$ be as in \cref{prop:ConcavityDistance}, and assume for simplicity that $(a,b)=(-a_0,a_0)$. Let $q=\sigma(0)$ and $p \in \gamma$ such that $\dist(\sigma(0),\gamma)=\dist(\sigma(0),p)$, and denote $\tilde\sigma(r)\eqdef\sigma(-r)$ for $r \in [0,a_0)$. Let $s_i<0<t_i$ be sequences $s_i\to0^-$, $t_i\to 0^+$ such that $[\sigma(t_i)p]$ and $[\sigma(s_i)p]$ converge to some geodesics $\gamma_+$ and $\gamma_-$, respectively. By \cite[Theorem 4.5.6, Corollary 4.5.7]{BuragoBuragoIvanovBook}, it follows that
\begin{equation*}
\begin{split}
        -\cos\angle\gamma_+'(0)\sigma'(0) &= \lim_{t_i\to0^+} \frac{\dist(\sigma(t_i),p) - \dist(q,p)}{t_i} = - \cos \theta_+, \\
        -\cos\angle\gamma_-'(0)\tilde\sigma'(0) &= \lim_{s_i\to0^-} \frac{\dist(\sigma(s_i),p) - \dist(q,p)}{-s_i} = - \cos \theta_-,
\end{split}
\end{equation*}
where
\[
\begin{split}
    \theta_+ &\eqdef \min\{ \angle \tilde\gamma'(0)\sigma'(0) \st \tilde\gamma \text{ geodesic from $q$ to $p$} \}, \\
    \theta_- &\eqdef \min\{ \angle \tilde\gamma'(0)\tilde \sigma'(0)x \st \tilde\gamma \text{ geodesic from $q$ to $p$} \}.
\end{split}
\]
In particular $\angle\gamma_+'(0)\sigma'(0)= \theta_+$, $\angle \gamma_-'(0)\tilde\sigma'(0)=\theta_-$, and $\theta_-\le \pi-\theta_+$. Hence
\begin{equation}\label{eq:TetaMenoTetaPiu}
    \cos \theta_- \ge - \cos \theta_+.
\end{equation}

The proof of \cref{prop:ConcavityDistance} is inspired by \cite[Proposition 3.7]{AlexanderBishop}, as the geodesic $\gamma$ in \cref{prop:ConcavityDistance} plays the role of a boundary for the $2$-dimensional space $X$. The argument is originally due to Perelman \cite{Perelman91}, see also \cite[Theorem 3.3.1]{PetruninSurvey} for a different proof.

\begin{proof}[Proof of \cref{prop:ConcavityDistance}]
We prove that $$\limsup_{h\to0^+} \Big(\dist(\sigma(t_0+h),\gamma) + \dist(\sigma(t_0-h),\gamma) -2 \dist(\sigma(t_0),\gamma) \Big)/h^2 \le 0,$$
for any $t_0 \in (a,b)$. This readily implies that $\dist(\sigma(t),\gamma)$ is concave. We prove such claim by constructing suitable upper touching concave functions.

We assume for simplicity that $(a,b)=(-a_0,a_0)$, and we prove the claim for $t_0=0$. Let $s_i,t_i$, $\gamma_+, \gamma_-$, $\theta_+,\theta_-$ be as above. Denote $q_i^-\eqdef \sigma(s_i)$, $q_i^+\eqdef \sigma(t_i)$. Let $\triangle \bar p \bar q \bar q_i^+$ and $\triangle \bar p \bar q \bar q_i^-$ be comparison triangles in $\R^2$ sharing the same edge $[\bar p \bar q]$. Let $\bar\gamma$ be the line in $\R^2$ passing through $\bar p$ and orthogonal to $[\bar p \bar q]$; assume that $\bar \gamma(0)=\bar p$ and that for $t>0$ the point $\bar \gamma(t)$ lies on the same side of $\triangle \bar p \bar q \bar q_i^+$ with respcet to the axis $[\bar p \bar q]$. Without loss of generality we can assume that $\bar p =0$ is the origin and that $\bar q \in \{0\}\times (0,+\infty)\subset \R^2$. Finally let $\bar \sigma:\R\to \R^2$ be the piecewise linear curve defined by
\[
\bar \sigma(t) \eqdef
\begin{cases}
\bar q + (\sin \theta_+, -\cos \theta_+) t & t\ge 0, \\
\bar q + (\sin \theta_-, \cos \theta_-) t & t< 0.
\end{cases}
\]
We consider the function $\bar h(t)\eqdef |\bar \sigma(t) \bar \gamma|$, that is
\[
\bar h(t) = \begin{cases}
|\bar p \bar q| - (\cos\theta_+)t & t\ge0, \\
|\bar p \bar q| + (\cos\theta_-)t & t<0.
\end{cases}
\]
Hence $\bar h (0) = |\bar p \bar q |=\dist(\sigma(0),\gamma)$ and $\bar h$ is concave by \eqref{eq:TetaMenoTetaPiu}.

We want to show that $\dist(\sigma(r_i),\gamma) \le \bar h(r_i)$ for $r_i \in \{t_i,s_i\}$, for any $i$. We show the latter for the times $t_i$. In the following, $[pq]$ denotes the inverse parametrization of $\gamma_+$.

By angle condition, $\angle \bar q_{\bar p}^{\bar q_i^+} \le \angle \gamma_+'(0)\sigma'(0)= \angle \bar q_{\bar p}^{\bar \sigma(1)}$, hence
\begin{equation}\label{eq:zzhti}
    |\bar q_i^+\bar \gamma| \le \bar h(t_i),
\end{equation}
for any $i$. Denote by $\bar r_i \in \bar \gamma$ the projection of $\bar q_i^+$ on $\bar \gamma$, so that $|\bar q_i^+\bar \gamma|=|\bar q_i^+\bar r_i|$.\\
Let $p_i \in \gamma\setminus\{p\}$ such that $\dist(p,p_i)=|\bar p \bar r_i|$ and $\angle p^{q_i^+}_{p_i} \le \tfrac\pi2$. By \cref{rem:AngleContinuity}, we know that, up to subsequence, $\lim_i \angle p^{q_i^+}_{p_i} = \angle [pq]'(0)\gamma'(0)=\tfrac\pi2$. 
Since $X$ has dimension $2$, then
\begin{equation}\label{eq:zaz}
\angle p^{q}_{q_i^+} +  \angle p^{q_i^+}_{p_i} = \frac\pi2
\qquad\text{or}\qquad
\angle p^{q}_{q_i^+} =  \angle p^{q_i^+}_{p_i} + \frac\pi2.
\end{equation}
Indeed, the previous alternative can be argued by passing to the tangent cone at $p$, which is either the half-plane or $\mathbb R^2$, and distinguishing the cases $\angle p^{q}_{q_i^+}<\pi/2$ or $\angle p^{q}_{q_i^+}>\pi/2$.

Suppose by contradiction that the second alternative in \eqref{eq:zaz} occurs along some (non-relabeled) subsequence, then $\exists\, \lim_i \angle p^{q}_{q_i^+}= \pi$.
Observe that $[pq_i^+]$ cannot intersect $\gamma$ for positive times; indeed, if $[pq_i^+](t_0)= \gamma(r_0)$ for some $t_0>0$ and $r_0$, then $r_0=r_1+t_0$, where $p=\gamma(r_1)$, and since $t\mapsto \dist([pq_i^+](t), \gamma(r_1+t))/t$ is nonincreasing, then $[pq_i^+](t)$ would coincide with $\gamma(r_1+t)$ for $t\ge t_0$, and then $q_i^+ \in \gamma$, which contradicts the assumption $0<\dist(\sigma(t_i),\gamma)$. For the same reason, $[pq]$ cannot intersect $\gamma$ for positive times.\\
For large $i$, we have $\angle p^{q_i^+}_{p_i}>\pi/4$ and $ \angle p^{q}_{q_i^+}>3\pi/4$. Hence, for some $$\rho\in (0,\min\{|pq|/2,|p\gamma(0)|/2,|p\gamma(T)|/2\})$$ small enough, the ball $B_\rho(p)$ is homeomorphic to an open disk in the plane; indeed, the inequalities satisfied by the angles $\angle p^{q_i^+}_{p_i}, \angle p^{q}_{q_i^+}$ imply that the tangent cone of $X$ at $p$ is $\R^2$ and Perelman's local structure theorem applies, see \cite[Section 13.2, Theorem b)]{BuragoGromovPerelman}. Hence $\gamma$ divides $B_\rho(p)$ in two open connected components $A_1,A_2$, and $[pq]$ intersects $A_1$ only, up to renaming. Since $\angle p^{q_i^+}_{p_i}>\pi/4$ and $ \angle p^{q}_{q_i^+}>3\pi/4$, then $[pq_i^+](t) \in A_2$ for $t \in (0,\eps_i)$; indeed if otherwise there is $\tau_j\to0$ with $[pq_i^+](\tau_j) \in A_1$, rescaling the distance by $\tau_j^{-1}$ and passing to the tangent cone at $p$ we would get that $\angle p^q_{q_i^+}< \pi/2$. Therefore, since $[pq_i^+]$ does not intersect $\gamma$ for positive times, then $[pq_i^+](0,\rho)\subset A_2$, for large $i$. But $[pq_i^+]$ converges to $[pq]$ uniformly, that gives a contradiction.

Therefore $\angle p^{q}_{q_i^+} +  \angle p^{q_i^+}_{p_i} = \tfrac\pi2$, for large $i$. Hence, letting $\triangle \bar p \bar q_i^+ \bar p_i$ be a comparison triangle for $\triangle  p  q_i^+ p_i$ which extends $\triangle \bar p \bar q \bar q_i^+$ to a quadrilateral, the angle condition implies
\[
\angle \bar p^{\bar q_i^+}_{\bar p_i} \le \angle  p^{q_i^+}_{p_i} = \frac\pi2 - \angle p^{q}_{q_i^+} \le \frac\pi2 - \angle \bar p^{\bar q}_{\bar q_i^+} = \angle \bar p^{\bar q_i^+}_{\bar r_i}. 
\]
The previous inequality together, recalling that $|\bar p \bar p_i|= \dist(p,p_i) = |\bar p \bar r_i|$, implies
\[
\begin{split}
\dist(\sigma(t_i),\gamma) &\le \dist(q_i^+, p_i) = |\bar q_i^+ \bar p_i| \le |\bar q_i^+ \bar r_i| \\
&\overset{\eqref{eq:zzhti}}{\le}  \bar h(t_i),
\end{split}
\]
for large $i$. The analogous argument can be employed to show that $\dist(\sigma(s_i),\gamma) \le \bar h(s_i)$, for large $i$.

Now $t_i,s_i$ are arbitrary infinitesimal sequences, so we can choose $s_i=-t_i$ with $t_i$ such that
\[
\lim_i \frac{\dist(\sigma(t_i), \gamma) + \dist(\sigma(-t_i), \gamma) - 2 \dist(\sigma(0),\gamma)}{t_i^2} = 
\limsup_{t\to0^+} \frac{\dist(\sigma(t), \gamma) + \dist(\sigma(-t), \gamma) - 2 \dist(\sigma(0),\gamma)}{t^2}.
\]
Hence
\[
\limsup_{t\to0^+} \frac{\dist(\sigma(t), \gamma) + \dist(\sigma(-t), \gamma) - 2 \dist(\sigma(0),\gamma)}{t^2}
\le \limsup_i \frac{\bar h(t_i) +\bar h(-t_i) - 2 \bar h(0)}{t_i^2} \le 0.
\]

\end{proof}
}

\printbibliography[title={References}]

\typeout{get arXiv to do 4 passes: Label(s) may have changed. Rerun} %Questo comando pare che sia per far funzionare poi la compilazione su arXiv

\end{document}